\documentclass[10pt, a4paper]{amsart}
\usepackage{amsfonts,amsmath,latexsym,amssymb}
\usepackage{graphicx}
\usepackage{longtable}
\allowdisplaybreaks[4]

\def\C{\mathbb{C}}
\def\N{\mathbb{N}}
\def\P{\mathbb{P}}%
\def\R{\mathbb{R}}
\def\Z{\mathbb{Z}}
\def\Bs{\mathop{\mathrm{Bs}}\nolimits}
\def\Cls{\mathop{\mathrm{Cls}}\nolimits}
\def\disc{\mathop{\mathrm{disc}}\nolimits}
\def\Disc{\mathop{\mathrm{Disc}}\nolimits}

\def\Int{\mathop{\mathrm{Int}}\nolimits}
\def\Item#1{\par\hangindent\parindent\indent\llap{#1\enspace}\ignorespaces}
\def\Ker{\mathop{\mathrm{Ker}}\nolimits}
\def\Krulldim{\mathop{\mathrm{Krulldim}}\nolimits}

\def\lto{\longrightarrow}

\def\mapr#1{\allowbreak\displaystyle \mathrel{\mathop{\longrightarrow}\limits^{#1}}\allowbreak}

\def\rank{\mathop{\mathrm{rank}}\nolimits}
\def\Rat{\mathop{\mathrm{Rat}}\nolimits}
\def\Reg{\mathop{\mathrm{Reg}}\nolimits}

\def\sign{\mathop{\mathrm{sign}}\nolimits}
\def\Sing{\mathop{\mathrm{Sing}}\nolimits}
\def\Spec{\mathop{\mathrm{Spec}}\nolimits}

\def\subsetne{\subsetneq}

\def\Zar{\mathop{\mathrm{Zar}}\nolimits}
\def\cC{\mathcal{C}}
\def\cD{\mathcal{D}}
\def\cE{\mathcal{E}}
\def\cF{\mathcal{F}}
\def\cH{\mathcal{H}}
\def\cI{\mathcal{I}}
\def\cJ{\mathcal{J}}
\def\cK{\mathcal{K}}
\def\cL{\mathcal{L}}
\def\cO{\mathcal{O}}
\def\cP{\mathcal{P}}
\def\cR{\mathcal{R}}
\def\cZ{\mathcal{Z}}
\def\frC{\mathfrak{C}}
\def\fre{\mathfrak{e}}
\def\frf{\mathfrak{f}}
\def\frg{\mathfrak{g}}
\def\frh{\mathfrak{h}}
\def\frp{\mathfrak{p}}
\def\frm{\mathfrak{m}}
\def\frn{\mathfrak{n}}
\def\frq{\mathfrak{q}}
\def\frS{\mathfrak{S}}

\catcode`\@=11
\outer\def\proclaim#1{%
  \removelastskip\penalty-400\vskip0.8em plus0.3em minus0.3em
  {\bf#1.}}%
\def\endproclaim{\par\penalty-400\vskip0.8em plus0.3em minus0.3em}%
\def\niceskip{\removelastskip\penalty-400\vskip0.8em plus0.3em minus0.3em}%
\def\Proof{\removelastskip\par\begin{proof}}%

\catcode`\@=11
\newdimen\htqed \newdimen\wdqed \newdimen\dpqed
\def\hidehrule#1#2{\kern-#1 \hrule height#1 depth#2 \kern-#2 }
\def\hidevrule#1#2{\kern-#1{\dimen0=#1
    \advance\dimen0 by#2\vrule width\dimen0}\kern-#2 }
\def\makeblankbox#1#2{\hbox{\lower\dpqed\vbox{\hidehrule{#1}{#2}\kern-#1 %
    \hbox to \wdqed{\hidevrule{#1}{#2}%
    \raise\htqed\vbox to #1{}%
    \lower\dpqed\vtop to #1{}%
    \hfil\hidevrule{#2}{#1}}%
    \kern-#1\hidehrule{#2}{#1}}}}
\def\QED{\htqed=6.7pt \dpqed=0pt \wdqed=6.7pt
    {\unskip\nobreak\hfil\penalty50\quad\null\nobreak\hfil
    {\hbox{\makeblankbox{0.17pt}{0.17pt}}}
    \parfillskip0pt\finalhyphendemerits0\par\medskip}}

\catcode`\@=\active
\catcode`\@=11

\def\Tqaa{Theorem 1.1}
\def\Tqab{Theorem 1.2}
\def\Tqabn{1.2}
\def\Tqabc{Corollary 1.3}
\def\Tqac{Theorem 1.4}
\def\Tqacn{1.4}
\def\Tqade{Proposition 1.5}
\def\Tqadef{Proposition 1.6}
\def\Tqadeg{Proposition 1.7}
\def\Tqad{Theorem 1.7}
\def\Tqadn{1.8}
\def\Tqae{Theorem 1.9}
\def\Tqah{Theorem 1.10}
\def\Tqai{Proposition 1.11}
\def\Tqba{Definition 2.1}
\def\Tqbb{Definition 2.3}
\def\Tqbdb{Definition 2.6}
\def\Tqbd{Theorem 2.7}
\def\Tqbae{Theorem 2.8} 
\def\Tqbbz{Definition 2.9}
\def\Tqbby{Theorem 2.10}
\def\Tqbbb{Example 2.11} 
\def\Tqbab{Corollary 2.13} 
\def\Tqbag{Proposition 2.14} 
\def\Tqbaiz{Theorem 2.15} 
\def\Tqbaf{Theorem 2.16} 
\def\Tqccn{Example 3.2} 
\def\Tqca{Proposition 3.3}
\def\Tqcb{Theorem 3.4}
\def\Tqccea{Lemma 3.5}%
\def\Tqcceb{Lemma 3.6}
\def\Tqch{Lemma 3.7}
\def\Tqci{Theorem 3.8}
\def\Tqcin{3.8}%
\def\Tqcn{Lemma 3.16}
\def\Tqcbd{Theorem 4.1}%
\def\Tqcbbb{Proposition 4.2}%
\def\Tqcbbd{Proposition 4.3}%
\def\Tqcbb{Lemma 4.4}%
\def\Tqcbf{Lemma 4.5}%
\def\Tqcbg{Lemma 4.6}%
\def\Tqcbh{Lemma 4.7}%
\def\Tqcbi{Lemma 4.8}%
\def\Tqcbj{Lemma 4.10}%
\def\Tqcca{Lemma 4.11}%
\def\Tqccbdm{Proposition 4.12}%
\def\Tqccbd{Theorem 4.13}%
\def\Tqccb{Theorem 4.15}%
\def\Tqccd{Proposition 4.16}%
\def\Tpga{Definition 5.1} %
\def\Tpgb{Definition 5.2} %
\def\Tpgc{Definition 5.5} %
\def\Tpgk{Proposition 5.10} %
\def\Tpgla{Theorem 5.11} %
\def\Tpglc{Corollary 5.13} %
\def\Tpgm{Corollary 5.14} %
\def\Tpgn{Theorem 5.15} %

\begin{document}
\title{Some Cubic and Quartic Inequalities of Four Variables}

\vskip2mm

\author{Tetsuya ANDO}

\address{
Department of Mathematics and Informatics, 
Chiba University, \\
Yayoi-cho 1-33, Inage-ku, 
Chiba 263-8522, JAPAN \\}
\email{ando@math.s.chiba-u.ac.jp}

\date{23.09.2026}  

\subjclass[2010]{26D05, 14P10, 14Q05}

\keywords{Algebraic inequalities, Positive Semidefinite Cone, Semialgebraic variety, Positive semidefinite forms.}

\begin{abstract}
Let $\cH \subset \cH_{n,d} := \R[x_1$,$\ldots$, $x_n]_d$ be a 
vector space, and $A$ be a compact semialgebraic subset of $\P_{\R}^{n-1}$. 
We shall study some PSD cones 
$\cP = \cP(A$, $\cH) := \big\{f \in \cH$ $\big|$ 
 $f(a) \geq 0$ ($\forall a \in A$)$\big\}$. 
Our interests are (1) to determine the extremal elements of $\cP$, 
(2) to determine discriminants of $\cP$, 
(3) to describe $\cP$ as a union of basic semialgebraic subsets, 
and (4) to find a nice test set when $\dim \cH$ is low. 
In this article, we present (1), (2), (3) and (4) for 
$\cP(\R^4$, $\cH_{4,4}^{s0})$ and 
$\cP(\R_+^4$, $\cH_{4,4}^{s0})$, where 
$\cH_{n,d}^{s0} := \big\{f \in \cH_{n,d}$ 
 $\big|$ $f$ is symmetric and $f(1,\ldots,1)=0 \big\}$. 
We also provide (1)---(4) for 
$\cP(\R_+^4$, $\cH_{4,3}^{c0})$, where 
$\cH_{n,d}^{c0} := \big\{f \in \cH_{n,d}$ 
 $\big|$ $f$ is cyclic and $f(1,\ldots,1)=0 \big\}$. 
\end{abstract}

\maketitle




\section{Introduction}
Let $\cH_{n,d} := \R[x_1$,$\ldots$, $x_n]_d$ (the part of degree $d$), 
and $\cH \subset \cH_{n,d}$ be a vector subspace. 
For a semialgebraic subset $A$ of $\R^n$, 
\[\cP(A, \, \cH) 
 := \big\{f \in \cH \; \big| \; 
      \hbox{$f(a) \geq 0$ for all $a \in A$}\big\}\]
is called the PSD cone on $A$ in $\cH$. 
Our interests are: \par
{\parindent=30pt
\Item{(I1)} To determine all the extremal elements of 
$\cP := \cP(A$, $\cH)$. 
\Item{(I2)} To determine all the discriminants of $\cP$ (see \Tqbdb). 
\Item{(I3)} To describe $\cP$ as a union of 
basic semialgebraic subsets using some inequalities. 
\Item{(I4)} Find a nice test set for $(A$, $\cH)$ when 
$\dim \cH$ is low (see \Tqbbz). 

}
In this article, we present (I1), (I2), (I3) and (I4) for PSD cones 
$\cP_{4,4}^{s0}$, $\cP_{4,4}^{s0+}$ and $\cP_{4,3}^{c0+}$. 
We also treat some SOS problems relating these PSD cones. 
We shall explain these symbols. 
Let 
\begin{align*}
 & \cH_{n,d}^c := \big\{ f \in \cH_{n,d} \; \big| \; \hbox{
     $f(x_2,\ldots,x_n,x_1) = f(x_1,\ldots,x_n)$}\big\}, \\
 & \cH_{n,d}^s := \big\{ f \in \cH_{n,d} \; \big| \; \hbox{
     $f(x_{\sigma(1)},\ldots,x_{\sigma(n)}) = f(x_1,\ldots,x_n)$ for 
      all $\sigma \in \frS_n$}\big\}, \\
 & \cH_{n,d}^0 := \big\{ f \in \cH_{n,d} \; \big| \; \hbox{
     $f(a,a,\ldots,a) = 0$ for all $a \in \R$}\big\}, \\
 & \cE(\cP) := \big\{ f \in \cP \; \big| \; 
     \hbox{$f$ is a extremal element of $\cP$} \big\}, \\
 & \R_+ := \big\{ x \in \R \; \big| \; \hbox{$x \geq 0$} \big\}, 
\end{align*}
and $\cH_{n,d}^{c0} := \cH_{n,d}^c \cap \cH_{n,d}^0$, 
$\cH_{n,d}^{s0} := \cH_{n,d}^s \cap \cH_{n,d}^0$. 
We denote 
$\cP_{n,d} := \cP(\R^n$, $\cH_{n,d})$, 
$\cP_{n,d}^+ := \cP(\R_+^n$, $\cH_{n,d})$, 
$\cP_{n,d}^s := \cP(\R^n$, $\cH_{n,d}^s)$, 
$\cP_{n,d}^{s+} := \cP(\R_+^n$, $\cH_{n,d}^s)$, 
$\cP_{n,d}^{s0} := \cP(\R^n$, $\cH_{n,d}^{s0})$, 
$\cP_{n,d}^{s0+} := \cP(\R_+^n$, $\cH_{n,d}^{s0})$, 
$\cP_{n,d}^c := \cP(\R^n$, $\cH_{n,d}^c)$, 
$\cP_{n,d}^{c+} := \cP(\R_+^n$, $\cH_{n,d}^c)$, 
$\cP_{n,d}^{c0} := \cP(\R^n$, $\cH_{n,d}^{c0})$, 
and $\cP_{n,d}^{c0+} := \cP(\R_+^n$, $\cH_{n,d}^{c0})$. 
The rule of indexing will be clear. 
``c'' means cyclic, ``s'' means symmetric, 
``0'' means an equality condition $f(a$,$\ldots$, $a)=0$, 
and ``$+$'' means $A = \R_+^n$. 

We have already completed (I1), (I2) and (I3) for the PSD cones 
$\cP_{3,3}^{c+}$, $\cP_{3,3}^{c0+}$, 
$\cP_{3,4}^{c0}$, $\cP_{3,4}^{c0+}$, 
$\cP_{3,4}^s$ and $\cP_{3,5}^{s0+}$. 
See \cite{RefAc}, \cite{RefAb} and \cite{RefAa}. 
For $\cP_{3,4}^{c0}$, see also \cite{RefCa} and \cite{RefMM}. 
(I4) for $\cP_{3,3}^{c+}$ is provided in \Tqbbb. 
(I1) for $\cP_{3,3}^+$ is given in \cite{RefAd}. 

In \S 3, we study $\cP_{4,4}^{s0}$ and $\cP_{4,4}^{s0+}$. 
(I1)---(I4) for $\cP_{4,4}^{s0}$ are given in \Tqcb, 
and these for $\cP_{4,4}^{s0+}$ are given in \Tqci. 
Here, we present (I3) for $\cP_{4,4}^{s0}$ and $\cP_{4,4}^{s0}$ 
slightly different style from \Tqcb \ and \Tqcin. 

\def\Tqaa{Theorem 1.1}
\proclaim{Theorem 1.1} 
{\sl Let $\sigma_1 := a_0+a_1+a_2+a_3$, 
\[\sigma_2 := \sum_{0 \leq i<j \leq 3} a_i a_j, \hskip15pt 
  \sigma_3 := \sum_{0 \leq i<j<k \leq 3} a_i a_j a_k,\] 
and $\sigma_4 := a_0 a_1 a_2 a_3$. 
Consider a family of quartic symmetric polynomials 
\[f(a_0,a_1,a_2,a_3)
  = \sigma_1^4 + p_1 \sigma_1^2 \sigma_2 + p_2 \sigma_2^2 
     + p_3 \sigma_1 \sigma_3 - (256 + 96p_1 + 36p_2 + 16p_3) \sigma_4
  \in \cH_{4,4}^{s0}\]
{\rm ($p_1$, $p_2$, $p_3 \in \R$)}. Then \par}
{\parindent=20pt
\Item{(1)} {\sl $f(a_0,a_1,a_2,a_3) \geq 0$ for 
all $a_0$,$\ldots$, $a_3 \in \R$ if and only if 
$16+6p_1+2p_2+p_3 \geq 0$ and $9p_1^2 \leq 128 + 24p_1 + 36p_2 + 12p_3$. }
\Item{(2)} {\sl $f(a_0,a_1,a_2,a_3) \geq 0$ for 
all $a_0 \geq 0$,$\ldots$, $a_3 \geq 0$ if and only if 
``{\rm (i)} or {\rm (ii)}'' and ``{\rm (iii)} or {\rm (iv)}'' hold: \par}
{\parindent=50pt
\Item{(i)} {\sl $p_1 \leq -8$ and $p_1^2 \leq 4p_2$.} 
\Item{(ii)} {\sl $p_1 \geq -8$ and $4p_1+p_2+16 \geq 0$.} 
\Item{(iii)} {\sl $p_1 \leq -14/3$ and 
 $9p_1^2 \leq 128 + 24p_1 + 36p_2 + 12p_3$.}
\Item{(iv)} {\sl $p_1 \geq -14/3$ and $27 + 9p_1 + 3p_2 + p_3 \geq 0$.}

}}
\endproclaim

Next, we present (I1). 

\def\Tqab{Theorem 1.2}
\proclaim{Theorem 1.2} 
{\sl All the extremal elements of $\cP_{4,4}^{s0}$ are 
positive multiples of the following polynomials: 
\begin{align*}
 & \frg_t(a_0,a_1,a_2,a_3) := \frac{1}{3}\Big(
       3 \sigma_1^4 - 2(t+7) \sigma_1^2 \sigma_2 + (t+3)^2 \sigma_2^2 \\
 & \hskip100pt  - 2(t^2-9) \sigma_1 \sigma_3 -4(t+3)^2 \sigma_4 \Big), \\
 & \frg_{\infty}(a_0,a_1,a_2,a_3) 
     := \sigma_2^2 - 2 \sigma_1 \sigma_3 - 4 \sigma_4, \\
 & \frp(a_0,a_1,a_2,a_3) 
     := \sigma_2^2 - 3 \sigma_1 \sigma_3 + 12 \sigma_4. 
\end{align*}
Here, $t \in \R$. 
Conversely, these are extremal elements of $\cP_{4,4}^{s0}$. 

$\frg_t$ ($t \ne 1$, $-3$) is characterized by the equality conditions 
$\frg_t(t$, $1$, $1$, $1) = \frg_t(-1$, $-1$, $1$, $1) = 0$. 
$\frg_1$ is characterized by the equality conditions 
$\frg_1(x$, $x$, $1$, $1) = 0$ for all $x \in \P_{\R}^1$. 
$\frg_{-3}$ is characterized by the equality conditions 
$\frg_{-3}(a$, $b$, $c$, $-a-b-c) = 0$ for all $a$, $b$, $c \in \R$. 
$\frg_{\infty}$ is characterized by the equality conditions 
$\frg_{\infty}(0$, $0$, $0$, $1) 
 = \frg_{\infty}(-1$, $-1$, $1$, $1) = 0$. 

$\frp$ is characterized by the equality conditions 
$\frp(0$, $0$, $0$, $1) = 1$ and 
$\frp(s$, $1$, $1$, $1) = 0$ for all $s \in \R$. 

}
\endproclaim

We say $f$ is {\it characterized by the equality conditions} 
$f({\bf x}_{\lambda}) = 0$ ($\lambda \in \Lambda$) if 
\[\R_+ \cdot f := \big\{ g \in \cP \; \big| \; 
   \hbox{$g({\bf x}_{\lambda}) = 0$ for all $\lambda \in \Lambda$}\big\}.\]
Note that if $f \in \cP$ is characterized by certain equality conditions, 
then $f$ is extremal. 
About the converse, please read \cite{RefAd}. 

\niceskip

An elements $f \in \cP_{n,2d}$ is called {\it SOS}, 
if there exists $r \in \N$ and $g_1$,$\ldots$, $g_r \in \cP_{n,d}$ 
such that $f = g_1^2 + \cdots + g_r^2$. 
The set of all the SOS elements in $\cP_{n,2d}$ are 
written by the symbol $\Sigma_{n,2d}$, and is called a {\it SOS cone}. 
In this case, 
$\frg_t$, $\frg_{\infty}$, $\frp \in \Sigma_{4,4}$, since 
\begin{align*}
 3\frg_t(a,b,c,d) & = \big(a^2+b^2-c^2-d^2 + (t+1)(c d-a b)\big)^2 \\
   & \hskip20pt + \big(a^2-b^2+c^2-d^2 + (t+1)(b d-a c)\big)^2 \\
   & \hskip20pt + \big(a^2-b^2-c^2+d^2 + (t+1)(b c-a d)\big)^2 \\
   & = \frac{1}{16} \sum_{\tau \in \frS_4}
     \big(a_{\tau(0)} - a_{\tau(1)}\big)^2
     \big(2(a_{\tau(0)} + a_{\tau(1)}) 
                    - (t+1)(a_{\tau(2)}+a_{\tau(3)})\big)^2, \\
 \frg_{\infty}(a,b,c,d) 
   & = (a b-c d)^2 + (a c-b d)^2 + (a d-b c)^2, \\
 \frp(a,b,c,d) 
   & = (1/2)\big((a-b)^2(c-d)^2 + (a-c)^2(b-d)^2 + (a-d)^2(b-c)^2\big). 
\end{align*}
Here $(a_0$, $a_1$, $a_2$, $a_3)=(a$, $b$, $c$, $d)$. 
Moreover, $\frg_t$, $\frp \notin \cE(\cP_{4,4}^{s0+})$.
Thus we obtain: 

\def\Tqabc{Corollary 1.3}
\proclaim{Corollary 1.3} 
{\sl $\cP_{4,4}^{s0} \subset \Sigma_{4,4}$, and 
$\cE(\cP_{4,4}^{s0}) \cap \cE(\cP_{4,4}) = \emptyset$. }
\endproclaim

Remember that 
$\cE(\cP_{3,4}^{c0}) \cap \cE(\cP_{3,4}) = \emptyset$, 
because $f \in \cE(\cP_{3,4}^{c0})$ is not 
a square of a quadric polynomial (see \cite{RefCa}). 
The following theorem provides extremal elements which do not 
appear in \cite{RefTc}. 

\def\Tqac{Theorem 1.4}
\proclaim{Theorem 1.4} 
{\sl All the extremal elements of $\cP_{4,4}^{s0+}$ are 
positive multiples of the following polynomials: 
\begin{align*}
 \frf_t^{ab}(a_0,a_1,a_2,a_3) 
 & := (1/3)\Big(3\sigma_1^4 - 2(t+7) \sigma_1^2 \sigma_2 + 8(t+1) \sigma_2^2 \\
 & \hskip20pt + (t^2-6t+21) \sigma_1 \sigma_3 - 16(t^2+3) \sigma_4 \Big)
      \quad \hbox{\rm ($0 \leq t \leq 5$)}, \\
 \frf_t^c(a_0,a_1,a_2,a_3) 
 & := (1/9)\Big(9 \sigma_1^4 - 6(t+7) \sigma_1^2 \sigma_2 
        + (t+7)^2 \sigma_2^2 \\
 & \hskip20pt + 12(t-1)\sigma_1 \sigma_3 -12(t-1)(3t+13) \sigma_4 \Big)
     \quad \hbox{\rm ($t \geq 5$)}, \\
 \frp(a_0,a_1,a_2,a_3) 
 & := \sigma_2^2 - 3 \sigma_1 \sigma_3 + 12 \sigma_4, \\
 \frq_1(a_0,a_1,a_2,a_3) 
 & := \sigma_1^2 \sigma_2 - 4 \sigma_2^2 + 3 \sigma_1 \sigma_3
    = \sum_{i<j} a_i a_j(a_i-a_j)^2, \\
 \frq_2(a_0,a_1,a_2,a_3) 
 & := \sigma_1 \sigma_3 - 16 \sigma_4 
   = \frac{1}{4} \sum_{\tau \in \frS_4} 
          a_{\tau(0)} a_{\tau(1)}\big(a_{\tau(2)}-a_{\tau(3)}\big)^2. 
\end{align*}
Conversely, these are extremal elements of $\cP_{4,4}^{s0+}$. 

$\frf_t^{ab}$ ($0 \leq t<1$ or $1<t \leq 5$) is 
characterized by the equality conditions 
\[\frf_t^{ab}(t,1,1,1) = \frf_t^{ab}(0,0,1,1) = 0.\]
$\frf_1^{ab}$ is characterized by the equality conditions 
$\frf_1^{ab}(t,t,1,1) = 0$ for all $t \geq 0$ and 
$\displaystyle 
  \frac{\partial^2}{\partial a_0^2}\frf_1^{ab}(1,1,1,1) = 0$. 
$\frf_t^c$ ($t>5$) is characterized by the equality conditions 
\[\frf_t^c(t,1,1,1) = \frf_t^c(0,0,u,1) = 0,\]
where $u \in \R_+$ is any root of $3u^2-(t+1)u+3=0$. 
$\frp$ is characterized by the equality conditions 
\[\frp(0,0,0,1) = \frp_a(0,0,0,1) = \frp(x,1,1,1) = 0\]
for all $x \geq 0$. 
$\frq_1$ is characterized by the equality conditions 
\[\frq_1(1,1,1,0) = \frq_1(1,1,0,0) = \frq_1(1,0,0,0) = 0.\]
$\frq_2$ is characterized by the equality conditions 
$\frq_2(s,1,0,0) = 0$ for all $s \geq 0$. 
}
\endproclaim

By the above representation, we have $\frp(a^2$, $b^2$, $c^2$, $d^2)$, 
$\frq_i(a^2$, $b^2$, $c^2$, $d^2) \in \Sigma_{4,8}$ ($i=1$, $2$). 
But for $f=\frf_t^{ab}$ and $\frf_t^c$, we obtain: 

\def\Tqade{Proposition 1.5}
\proclaim{Proposition 1.5} 
{\sl If $0 < t \leq 5$ and $t \ne 1$, 
then $\frf_t^{ab}(a^2$, $b^2$, $c^2$, $d^2) 
   \notin \Sigma_{4,8}$. 
If $t > 5$, 
then $\frf_t^c(a^2$, $b^2$, $c^2$, $d^2) \notin \Sigma_{4,8}$. }
\endproclaim

It is clear that $\frp$, $\frq_1$, $\frq_2 
 \notin \cE(\cP_{4,4}^+)$. But we have: 

\def\Tqadef{Proposition 1.6}
\proclaim{Proposition 1.6} 
{\sl If $t > 5$, then 
$\frf_t^c \in 
  \cE(\cP_{4,4}^{s0+}) \cap \cE(\cP_{4,4}^+)$. }
\endproclaim

Remember that if $f \in \cE(\cP_{3,4}^{s0})$, 
$f$ can be written as $f = g \overline{g}$, 
where $g$ is an imaginal quadric polynomial. 

\def\Tqadeg{Proposition 1.7}
\proclaim{Proposition 1.7} 
(1) {\sl If $t \ne -3$, then $\frg_t$ is irreducible in 
$\C[a,b,c,d]$.} \par
{\parindent=20pt
\Item{(2)} {\sl If $0 \leq t \leq 5$, 
then $\frf_t^{ab}$ is irreducible in $\C[a,b,c,d]$.} 
\Item{(3)} {\sl If $t > 5$, 
then $\frf_t^c$ is irreducible in $\C[a,b,c,d]$.} 

}
\endproclaim

%
%

We should explain about the discriminants of $\cP = \cP(A$, $\cH)$. 
Let $s_0$, $s_1$,$\ldots$, $s_N$ be a basis of the vector space $\cH$, 
and let $\Phi_\cH \colon A \to\cdots \P_{\R}^N$ be the rational map 
defined by $\Phi_\cH({\bf a}) 
  = \big(s_0({\bf a})\colon \cdots \colon s_N({\bf a})\big)$. 
$X := \Phi_\cH(A)$ is called the {\it characteristic variety}. 
Let $\Delta(X) = \big\{D_1$,$\ldots$, $D_r\big\}$ be the critical 
decomposition of $X$ (see \Tqbb). 
Each $D \in \Delta(X)$ is a smooth semialgebraic variety, 
and $D$ has its dual variety $D^{\vee}$. 
Let $\disc(D)$ be the defining equation of the Zariski 
closure of $D^{\vee}$ in $\cH$, 
and let $V_\cH(\disc(D))$ be 
the zero locus of $\disc(D)$ in $\cH$. 
If $\dim \big(V_{\cH}(\disc(D)) \cap \partial \cP\big) 
  = \dim \cP - 1$, we say $\disc(D)$ is a {\it discriminant} of $\cP$. 
For any $f \in \partial \cP$, there exists $D \in \Delta(X)$ such 
that $f \in V_{\cH}(\disc(D))$. 
Assume that a subset $B \subset A$ satisfies $\Phi_\cH(B) = D$. 
Then, for each $f \in V_\cH(\disc(D)) \cap \partial \cP$, 
there exists a point ${\bf a} \in B$ such that $f({\bf a}) = 0$. 
In this case, we shall say that $\disc(D)$ is a discriminant 
{\it corresponding to} $B$. 

\def\Tqad{Theorem 1.8}
\proclaim{Theorem 1.8} 
{\sl Let's denote the elements of $\cH_{4,4}^{s0}$ as 
\[f(a_0,a_1,a_2,a_3)
  = p_0 \sigma_1^4 + p_1 \sigma_1^2 \sigma_2 + p_2 \sigma_2^2 
     + p_3 \sigma_1 \sigma_3 - (256p_0 + 96p_1 + 36p_2 + 16p_3) \sigma_4,\]
and use $(p_0$,$\ldots$, $p_3)$ as 
a coordinate system of $\cH_{4,4}^{s0}$. \par}
{\parindent=20pt
\Item{(1)} {\sl $\cP_{4,4}^{s0}$ has the following two discriminants: 
\[d_1 := 128p_0^2 + 24p_0p_1 + 36p_0p_2 + 12p_0p_3 - 9p_1^2, \quad 
  d_2 := 16p_0+6p_1+2p_2+p_3.\]
$d_1$ corresponds to $\big\{(t,1,1,1) \in \R^4$ $\big|$ 
   $t \in \R$, $t \ne -3$, $1\big\}$, 
and $d_2$ corresponds to a point $(1$, $1$, $-1$, $-1)$. }
\Item{(2)} {\sl $\cP_{4,4}^{s0+}$ has the following five discriminants: 
\begin{align*}
 & d_1 := 128p_0^2 + 24p_0p_1 + 36p_0p_2 + 12p_0p_3 - 9p_1^2, \quad 
   d_3 := 4p_0p_2-p_1^2, \\
 & d_4 := 27p_0 + 9p_1 + 3p_2 + p_3, \quad
   d_5 := 16p_0 + 4p_1+p_2, \quad 
   d_6 := p_0. 
\end{align*}
$d_3$ corresponds to $\big\{(0,0,t,1) \in \R^4$ $\big|$ $0<t<1\big\}$. 
$d_4$, $d_5$, $d_6$ corresponds to points $(1$, $1$, $1$, $0)$, 
$(1$, $1$, $0$, $0)$, $(1$, $0$, $0$, $0)$ respectively. }
\endproclaim

We explain about (I4). 
Let $r_0 := \max \{2$, $\lfloor d/2 \rfloor\}$. 
For general $f \in \cH_{n,d}^s$, Riener, Timofte and Harris proved that 
$f \in \cP_{n,d}^s$ if $f(x) \geq 0$ for all 
$x \in \big\{(x_1,\ldots,x_n) \in \R^n$ $\big|$ 
$\# \{x_1$,$\ldots$, $x_n\} \leq r_0\big\}$. 
Moreover, $f \in \cP_{n,d}^{s+}$ if $f(x) \geq 0$ for all 
$x \in \big\{(x_1,\ldots,x_n) \in \R^n$ $\big|$ 
 $\# \big(\{x_1$,$\ldots$, $x_n\}-\{0\}\big) \leq r_0\big\}$. 
(See \cite[Corollary 1.3]{RefRie}, \cite[Corollary 2.1]{RefTa}. 
See also \cite{RefTb}, \cite{RefTc}.) 

In the case $f \in \cH_{4,4}^s$, 
if $f(t$, $t$, $1$, $1) \geq 0$ and $f(t$, $1$, $1$, $1) \geq 0$ 
for all $t \in \R$ then $f \in \cP_{4,4}^s$. 
If $f(t$, $t$, $1$, $1) \geq 0$, $f(t$, $1$, $1$, $1) \geq 0$, 
$f(0$, $t$, $1$, $1) \geq 0$ and $f(0$, $0$, $t$, $1) \geq 0$ 
for all $t \geq 0$ then $f \in \cP_{4,4}^{s+}$. 

We prove that the number of test conditions can be decreased as 
the following theorem in the cases of $\cP_{4,3}^{s0}$ and $\cP_{4,3}^{s0+}$. 

\def\Tqae{Theorem 1.9}
\proclaim{Theorem 1.9} 
{\rm (1)} {\sl If $f \in \cH_{4,4}^{s0}$ satisfies 
$f(-1$, $-1$, $1$, $1) \geq 0$ and $f(t$, $1$, $1$, $1) \geq 0$ for 
all $t \in \R$, then $f(a$, $b$, $c$, $d) \geq 0$ 
for all $a$, $b$, $c$, $d \in \R$. }

{\rm (2)} {\sl If $f \in \cH_{4,4}^{s0}$ satisfies 
$f(t$, $1$, $1$, $1) \geq 0$ and $f(0$, $0$, $t$, $1) \geq 0$ for 
all $t \geq 0$, then $f(a$, $b$, $c$, $d) \geq 0$ 
for all $a$, $b$, $c$, $d \in \R_+$. }
\endproclaim

In \S 4, we study the PSD cone of cyclic cubic polynomials $\cP_{4,3}^{c0+}$. 
(I2) and (I3) for $\cP_{4,3}^{c0+}$ are given in \Tqcbd. 
(I1) is presented in \Tqccb. 
$\cP_{4,3}^{c0+}$ has 4 discriminants. 
Since one of them is very complicated polynomial, 
the structure of $\cP_{4,3}^{c0+}$ is not simple. 
We also need somewhat strange algebraic numbers to state (I3). 
This is completely different from cases of $\cP_{3,3}^{c0+}$ and 
$\cP_{3,3}^{c+}$. 
Please read \Tqcbd \ and \Tqccb. 
(I1) and (I4) for $\cP_{4,3}^{c0+}$ are as the following: 

\def\Tqah{Theorem 1.10}
\proclaim{Theorem 1.10} 
{\rm (1)} {\sl All the elements of $\cE(\cP_{4,3}^{c0+})$ is 
the positive multiple of $\fre_{u,v,w}^h$ 
{\rm ($(u \colon v \colon w) \allowbreak \in D_e^h$)} or 
$\fre_t^{P_2}$ {\rm ($t \in \P_{\R}^1$).}} 

{\rm (2)} {\sl If $f \in \cH_{4,3}^c$ satisfies 
$f(1$, $1$, $1$, $1) \geq 0$ and 
$f(0$, $s$, $t$, $1) \geq 0$ for all $s$, $t \in \R_+$, 
then $f(a$, $b$, $c$, $d) \geq 0$ for all $a$, $b$, $c$, $d \in \R_+$. }
\endproclaim

Definitions of $\fre_{u,v,w}^h$, $\fre_t^{P_2}$ and 
$D_e^h$ are given in \Tqcbbd, \Tqcbh \ and \Tqccbd \ respectively. 
(1) is proved in \Tqccb, and (2) is proved in \S 4.2. 

In \cite{RefAd}, we have proved that 
$\cE(\cP_{3,3}^{c0+}) \subset \cE(\cP_{3,3}^{c+}) 
  \subset \cE(\cP_{3,3}^+)$. 
But $\cE(\cP_{4,3}^{c0+}) \not\subset \cE(\cP_{4,3}^+)$. 
An extremal element $\fre_{u,v,w}^h \in \cE(\cP_{3,3}^{c0+})$ 
is somewhat complicated to present here. 
The definition of $\fre_{u,v,w}^h$ is given at \Tqcbbd, 
and $D_e^h$ is described in \Tqccbd. 
Relating SOS problem, $\fre_{u,v,w}^h$ satisfies: 

\def\Tqai{Proposition 1.11}
\proclaim{Proposition 1.11} 
{\sl Assume that $(u \colon v \colon w) \in D_e^h$, 
$u>0$, $v>0$, $w>0$ and $v \ne u+w$. 
Then, $\fre_{u,v,w}^h(a^2$, $b^2$, $c^2$, $d^2) 
\in \cP_{4.6} - \Sigma_{4,6}$. }
\endproclaim


In \S 5, we will give an exact definition of semialgebraic varieties, 
and prove some basic general theorems. 
In this article, we use $\P_{\R}^3/\frS_4$ and $\P_+^3/\frS_4$. 
These are not real algebraic variety. 
$\P_{\R}^3/\frS_4$ does not agree with a real weighted projective space 
$\P_{\R}(1,2,3,4)$. 
But we need to treat these with certain variety structure, i.e. 
semialgebraic varieties. 
So, the author think it will be better to give an exact definition of 
semialgebraic variety. 
For example, there exists continuous rational map 
which is not holomorphic (see \Tqccea). 
Such maps do not exist in complex algebraic geometry. 
Some results will be useful for studies of real algebraic varieties. 
Especially, \Tpgla \ and \Tpgn \ show that semialgebraic geometry 
is very different from complex algebraic geometry. 
In our theory of algebraic inequalities in this article, 
a phenomenon of \Tpgn \ occurs. For example, 
$\Phi_\cH \colon A \cdots\to X$ may include some exceptional set 
even if $A = \P_{\R}^3$ (see \Tqccea). 

\smallskip

We shall explain a short history of study of PSD cones. 
Originally, $\cP_{n,2d}$ is called a PSD cone. 
Hilbert proved, $\cP_{n,2d}=\Sigma_{n,2d}$ if and 
only if $n \leq 2$ or $2d=2$ or $(n$, $2d)=(3$, $4)$ (\cite{RefHil}). 
History of studies till 1991 are written in \cite[\S 6.6]{RefBCR}. 
So we don't explain them again. 
Choi and Lam found some extremal forms of $\cP_{n,2d}$ which 
don't belong to $\Sigma_{n,2d}$ in \cite{RefCL}.
In \cite{RefR}, Reznick studied the condition 
that $f \in \cP_{n,2d}$ is included in $\Sigma_{n,2d}$. 
He also studied the condition that $f \in \cP_{n,2d}$ is extremal. 
See also \cite{RefCLRb}. 
They implies that if $f \in \cE(\cP_{n,2d})$, 
then $V_{\R}(f)$ is larger set. 
This fact is formalized in \cite[Theorem 2.11, Proposition 2.13]{RefAd}. 

An element $f \in \cH_{n,2d}$ is called even, 
if $f \in \R[x_1^2$,$\ldots$, $x_n^2]$. 
Choi, Lam and Reznick studied $\cP_{n,2d}^{es} := 
\cP_{n,2d} \cap \R[x_1^2$,$\ldots$, $x_n^2]$ in \cite{RefCLR}. 
They studied the condition for $\cP_{n,2d}^{es} \subset \Sigma_{n,2d}$. 
Note that $\cP_{n,2d}^{es} \cong \cP_{n,d}^{s+}$, 
as is stated in \cite{RefCLRc}. 
Harris proved $\cP_{3,8}^{es} \subset \Sigma_{3,8}$ in \cite{RefHs}. 
But $\cE(\cP_{3,3}^+) \cong \cE(\cP_{3,6}^e) \subset \cE(\cP_{3,6})$ and 
$\cE(\cP_{3,6}^e) \not\subset \Sigma_{3,6}$ (see \cite{RefAd}). 

The relations $\cP_{n,2d}^s$ and $\Sigma_{n,2d}^s$ are 
studied by Goel, Kuhlmann and Reznick in \cite{RefGKR}. 
A related study can be found in \cite{RefGKBR}. 
Our study of $\cE(\cP_{4,4}^{s0+})$ 
and $\cE(\cP_{4,3}^{c0+})$ will give a small contribution for it. 

About discriminants of $\cP(A$, $\cH)$, 
Nie shown some interesting results in \cite{RefN}. 
He treated the case that $A$ is an affine real algebraic variety. 
In this article, we only treat the cases that $A$ is a compact 
semialgebraic variety. 
But they have very close relation. 
\cite{RefB} provides many nice ideas to treat algebraic inequalities 
using complex algebraic geometry. 

About $\cP_{3,6}$, $\Sigma_{3,6}$, 
$\cP_{4,4}$ and $\Sigma_{4,4}$, very important results 
are obtained in \cite{RefBHORS}. 
It provides relation with theory of K3 surfaces. 

$\dim \cH_{4,3} =  20$ and $\dim \cH_{4,4} =  35$ are 
somewhat large to proceed precise numerical analysis. 
It will not be insignificant to study some lower dimensional 
subspaces $\cH \subset \cH_{n,d}$. 

\smallskip

To check many calculations in this article, 
we will need Mathematica or a similar tool. 
The author provides a file for Mathematica 
in the authors WEB and in arXiv's anc folder. 
It will be useful for experimentation of inequalities. 

\removelastskip\penalty-400\vskip2.5em plus0.3em minus0.3em
{\bf \S 2. General theories}
\par\penalty1000\vskip0.8em plus0.2em minus0.2em
{\bf 2.1. Known results.}
\par\penalty1000\vskip0.4em plus0.1em minus0.1em
By studies in \cite{RefAc}, we have better to use $\P_{\R}^{n-1}$ and 
$\P_+^{n-1}$ instead of $\R^n$ and $\R_+^n$ where 
\[\P_+^n := (\R_+^{n+1} - \{0\})/\R_+^{\times}
 = \big\{(x_0\colon\cdots\colon x_n) \in \P_{\R}^n \; \big| \; 
      \hbox{$x_0 \geq 0$,$\ldots$, $x_n \geq 0$} \big\}.\]
The merits are that $\P_{\R}^{n-1}$ is compact and $\dim \P_{\R}^{n-1} < 
\dim \R^n$. 
But $f \in \cH_{n,d}$ is not a function on $\P_{\R}^{n-1}$. 
So, we must treat $\cH_{n,d}$ as a signed linear system on $\P_{\R}^{n-1}$. 
We need some more generalizations. 
About the exact definition of a semialgebraic variety, please see \S 5. 
We may understand here that a semialgebraic variety 
$(A$, $\cR_A)$ is a locally ringed space with semialgebraic set $A$ 
and a sheaf of rings $\cR_A$ which represent real holomorphic 
functions on open subsets of $A$. 
We only use $\cR_A$ to define singularities of $A$, 
regular maps between semialgebraic varieties, and signed linear systems. 
The author apologizes that \cite[Definition 1.7]{RefAc} \ must be 
corrected as the following: 

\def\Tqba{Definition 2.1}
\proclaim{Definition 2.1} 
{\rm Let $(A$, $\cR_A)$ be a semialgebraic variety, 
and $\cC^0_A$ be the sheaf of germs of real continuous functions on A. \par
{\parindent=20pt
\Item{(1)} Let $\cI$ be an invertible $\cR_A$-sheaf. 
$\cI$ is called a {\it signed invertible sheaf} on $A$ if 
{\parindent=40pt
\Item{(i)} there exists $\cC^0_A$-invertible sheaf $\cJ$ 
such that $\cI \otimes_{\cR_A} \cC^0_A = \cJ \otimes_{\cC^0_A} \cJ$, and 
\Item{(ii)} there exists $e \in \cJ(A)$ such 
that $e^2 \in \cI(A)$ and $\cI(A) = \cR_A(A) \cdot e^2$. 

}
Then, for $f \in H^0(A$, $\cI)$, 
there exists $g \in H^0(A$, $\cR_A)$ such that $f = g e^2$. 
We define $\sign(f(P)) \in \{0$, $\pm 1\}$ by $\sign(f(P)) = \sign(g(P))$ 
for $P \in A$. 
\Item{(2)} Let $\cI$ be a signed invertible $\cR_A$-sheaf. 
A finite dimensional vector subspace $\cH \subset H^0(A$, $\cI)$ is 
called a {\it signed linear system} on $A$. 
For $f \in \cH$, we say $f$ is {\it PSD} on $A$ if 
$f(P) \geq 0$ for all $P \in A$. 
\Item{(3)} The cone $\cP = \cP(A, \, \cH) 
 := \big\{f \in \cH \; \big| \; 
      \hbox{$f(P) \geq 0$ for all $P \in A$}\big\}$ 
is called the {\it PSD cone} on $A$ in $\cH$. 
Note that $\cP_{n,d} = \cP(\P_{\R}^{n-1}$, $\cH_{n,d})$ and 
$\cP_{n,d}^+ = \cP(\P_+^{n-1}$, $\cH_{n,d})$ and so on. 
\Item{(4)} The set $\Bs \cH := \big\{P \in A \; \big| \;
 \hbox{$f(P) = 0$ for all $f \in \cH$} \big\}$ 
is called the {\it base locus} of $\cH$. 
When $\cP$ is non-degenerate in $\cH$, 
we define $\Bs \cP := \Bs \cH$. 

}}
If $\dim \Bs \cP < \dim A$, we can define a rational map 
$\Phi_\cH : A \cdots\to \P_{\R}(\cH^{\vee})$, 
using a basis of $\cH$. 
$X = X(A$, $\cH) := \Cls(\Phi_\cH(A))$ (Euclidian closure) is called 
the {\it characteristic variety} of $A$. 
\endproclaim

For example, 
\[\cH_{n+1,d} := \big\{ f(x_0,\ldots, x_n) \; \big| \; 
 \hbox{$f$ is a homogeneous polynomial of degree $d$} \big\} \cup \{0\}\]
is a signed linear system on $\P_+^n$. 
For $f \in \cH_{n+1,d}$ and $P \in \P_+^n$, 
we cannot define the value $f(P)$ but can define $\sign(f(P))$. 
If $d$ is even, $\cH_{n+1,d}$ is also 
a signed linear system on $\P_{\R}^n$. 

\proclaim{Proposition 2.2} 
{\sl Let $X := X(A$, $\cH)$, and let $Y$ be the convex closure 
of $X$ in $\P(\cH^{\vee})$. Then
\[\cP(A, \, \cH) = \cP(X, \, \cH_{N+1,1}) 
    = \cP(Y, \, \cH_{N+1,1}),\]
where $\cH_{N+1,1}$ is the set of linear polynomials 
on $\P(\cH^{\vee})$. }
\endproclaim

\Proof
$\cP(A$, $\cH) = \cP(X$, $\cH_{N+1,1})$ is 
proved at \cite[Proposition 1.13]{RefAc}. 
$\cP(X$, $\cH_{N+1,1}) \allowbreak = \cP(Y$, $\cH_{N+1,1})$ is 
clear since every element of $\cH_{N+1,1}$ is linear. 
\end{proof}

\niceskip

Assume that a semialgebraic set $B$ is a subset of a complete 
real algebraic variety $V$. 
The minimal closed algebraic subset which contains $B$ is called 
the Zariski closure of $B$ and is denoted by $\Zar_V(B)$. 
We denote the Euclidian closure 
of $B$ in $V$ by $\Cls_V(B)$ or $\overline{B}$. 
Assume that $\Zar_V(B) = V$. 
The interior of $B$ is defined by $\Int(B) := V - \Cls_V(V - B)$. 
The boundary of $B$ is defined by $\partial B := B - \Int(B)$. 
Do not confuse with $\partial_V B := \Cls_V(B) - \Int(B)$. 
Note that $\Int(B)$ and $\partial B$ do not depend on the choice of $V$. 
But $\Cls_V(B)$ and $\partial_V B$ depend on $V$. 

\def\Tqbb{Definition 2.3}
\proclaim{Definition 2.3} 
{\rm (Critical decomposition. See \cite[Definition 1.5]{RefAc}) 
Let $A$ be a reduced semialgebraic variety with $\dim A = n$. 
We shall define $\Delta^i(A)$ ($i=0$,$\ldots$, $n$) by induction on $n$. 
If $\dim A = 0$, then $A = \{P_1$,$\ldots$, $P_m\}$ where $P_i$ are points. 
In this case we put $\Delta^0(A) = \{P_1$,$\ldots$, $P_m\}$, 
and put $\Delta^i(A) = \emptyset$ for $i \ne 0$. 

Assume that $n = \dim A \geq 1$. 
Let $Z_1$,$\ldots$, $Z_r$ be all the irreducible components of $A$ with 
$\dim Z_i = n$. 
Put $A_i := \Int(Z_i - \Sing(A)\big)$, 
and $\Delta^n(A) := \big\{A_1$,$\ldots$, $A_r\big\}$. 
Note that $Z_i \cap Z_j \cap \Int(A) \subset \Sing(A)$ for $i \ne j$. 

Let $Y_1$,$\ldots$, $Y_k$ be all the irreducible components of $A$ with 
$\dim Y_j \leq n-1$, and 
let $B_j := Y_j - (A_1 \cup \cdots \cup A_r)$. Put 
\[B := \Sing(A) \cup \partial A \cup B_1 \cup \cdots \cup B_k.\]
Then, we can regard $B$ to be a semialgebraic subvariety of $A$ with 
the reduced structure. 
Note that $\dim B < \dim A$. 
Thus we put $\Delta^i(A) := \Delta^i(B)$ for $i \ne n$. 

We denote $\Delta(A) := \Delta^0(A) \cup \Delta^1(A) \cup 
  \cdots \cup \Delta^n(A)$, 
and is called a {\it critical decomposition} of $A$. 
Each element $D \in \Delta(A)$ is called a {\it critical set} of $A$. 
Note that $D$ is a non-singular semialgebraic 
variety with $\partial D = \emptyset$. }
\endproclaim

\proclaim{Example 2.4} 
{\rm Consider the case $A = \P_+^2$. This is homeomorphic to a triangle. 
Let $P_x := (1 \colon 0 \colon 0)$, $P_y := (0 \colon 1 \colon 0)$, 
and $P_z := (0 \colon 0 \colon 1)$. 
For two points $P$, $Q \in \P_+^2$, we denote the open line segment 
connecting $P$ and $Q$ as $(PQ)$. 
Then, the critical decomposition of $\P_+^2$ is 
$\Delta^0(\P_+^2) = \big\{P_x$, $P_y$, $P_z \big\}$, 
$\Delta^1(\P_+^2) = \big\{(P_xP_y)$, $(P_yP_z)$, $(P_zP_x) \big\}$, 
$\Delta^2(\P_+^2) = \big\{\Int(\P_+^2)\big\}$. 

On the other hand, if $A = \P_{\R}^n$, 
then $\Delta^n(\P_{\R}^n) = \big\{\P_{\R}^n\big\}$, 
and $\Delta^r(\P_{\R}^n) = \emptyset$ for $r \ne n$. }
\endproclaim

\proclaim{Definition 2.5} 
{\rm (1) Let $X$ be a subset of $\R^n$ or $\P_{\R}^n$. 
$e \in X$ is said to be {\it extremal} in $X$, 
if $a>0$, $b>0$ and $x$, $y \in X$ satisfy $e = ax+by$ then $x = y = e$. 
Let $\cP$ be a closed convex cone which contain no lines. 
$0 \ne f \in \cP$ is called {\it extremal} in $\cP$, 
if $g$, $h \in X$ satisfy $f = g+h$ then 
$g$ and $h$ are multiples of $f$. 
For both cases $Y=X$ and $Y=\cP$, we denote that 
\[\cE(Y) := \big\{ y \in Y \; \big| \; 
  \hbox{$y$ is extremal in $Y$} \big\}.\]

(2) For a semialgebraic variety $A$ and $a \in A - \Bs \cH$ and 
a signed linear system $\cH$ on $A$, we put 
\[\cH_a := \big\{ f \in \cH \; \big| \; 
   \hbox{$f(a) = 0$}\big\}, \quad 
  \cP_a := \cP \cap \cH_a = \cP(A, \, \cH_a).\]
$\cP_a$ is called the {\it local cone} of $\cP$ at $a$. }
\endproclaim

Even if $a \in \Bs \cH$, we can define $\cP_a$ as 
\cite[Definition 2.6]{RefAd}. 
But we don't use it in this article. 

\def\Tqbdb{Definition 2.6}
\proclaim{Definition 2.6}
(See \cite[Definition 1.15, 1.17]{RefAc}) \par
{\rm (1) Let $\P = \P_{\R}^N$ and $\P^{\vee}$ 
be the set of all the hyperplanes in $\P$. 
Assume that $D \subset \P$ is a non-singular semialgebraic variety 
with $\partial D = \emptyset$ (i.e. $\Delta(D) = \{D\}$). 
For $x \in D$, let $T_{D,x} := T_{\Zar(D),x} \subset \P$ be the tangent space 
of $\Zar(D)$ at $x$. Then, 
\[D^{\vee} := \big\{ H \in \P^{\vee} \; \big| \; 
  \hbox{$H \supset T_{D,x}$ for a certain $x \in D$}\big\}\]
is called the {\it dual variety} of $D$. 
Since $D$ is irreducible and non-singular, $D^{\vee}$ is irreducible. 
Thus $D^{\vee}$ is a semialgebraic variety. 

(2) Under the same notation with \Tqba, let 
$\pi : (\cH - \{0\}) \to \P(\cH)$ be the natural surjection. 
For $D \in \Delta(X)$, we denote 
\[\cF(D) 
 := \Cls_\cH(\pi^{-1}(D^{\vee}) \cap \partial\cP).\]
If $\dim \cF(D) = \dim (\partial \cP)$, 
then $\cF(D)$ is called a {\it face component} 
of $\cP$ or of $\partial \cP$, 
and an irreducible defining equation 
of the Zariski closure $\Zar(\cF(D))$ is 
called a {\it discriminant} of $\cP$, 
and denoted by $\disc_D$ or $\disc(D)$. 

Especially, if $D \in \Delta^{\dim X}(X)$ and $\cF(D)$ is a face 
component, then $\cF(D)$ is called a {\it main component} of $\cP$, 
and $\disc(D)$ is called a {\it main discriminant} of $\cP$. }
\endproclaim


For example, if $X \cong \P_{\R}^n = A$, 
then $\cP$ has unique discriminant which is a main discriminant. 

In the case $D \in \Delta^0(X)$, $\disc(D)$ is linear. 
That is, if $\Phi_\cH$ is defined by basis $\{s_0$,$\ldots$, $s_N\}$ of 
$\cH$, and if we represent $f \in \cH$ as 
$f = p_0 s_0 + \cdots + p_N s_N$, 
and $D = (b_0 \colon \cdots \colon b_N) \in \P(\cH^{\vee})$, 
then $\disc(D) = b_0 p_0 + \cdots + b_N p_N$. 

\def\Tqbd{Theorem 2.7}
\proclaim{Theorem 2.7} (\cite[Theorem 1.18]{RefAc})} 
{\sl We use the same notation as \Tqba \ and the above.} \par
{\parindent=20pt
\Item{\rm (1)} {Let
\[\cD := \big\{ D \in \Delta(X) \; \big| \; 
  \hbox{$\cF(D)$ is a face component of $\cP$}\big\}.\]
Then 
$\displaystyle \partial\cP = \bigcup_{D \in \cD} \cF(D)$. }
\Item{\rm (2)} 
{\sl For $D \in \Delta(X)$, take a subset $B \subset A$ such that 
$\Phi_\cH(B) \subset D$ and $\Cls_D\big(\Phi_\cH(B)\big) 
\allowbreak = D$. 
Put $B_0 := B - \Bs \Phi_{\cH}$. Then, }
\[\cF(D) = \Cls_\cH
           \left(\bigcup_{a \in B_0} \cP_a\right).\]
\Item{\rm (3)} 
{\sl Assume that $\cP := \cP(X$, $\cH_{N+1,1})$ is 
non-degenerate in $\cH_{N+1,1}$. 
Take $x \in D \in \Delta^r(X)$. Then $\dim \cP_x \leq N-r$.}

}
\endproclaim

The author should apologize for that 
\cite[Proposition 1.27]{RefAc} is not correct. 
It should be corrected as (3) of the above theorem. 
We present a corrected proof of (3). 

\Proof 
(3) For $f \in \cH$, let $H_f$ be the hyperplane 
in $\P(\cH^{\vee})$ defined by $f = 0$. 
Since $\cP$ is non-degenerate, $\dim (U \cap \cP) = N+1$ for 
any Euclidean open neighborhood $U$ of $x$. 
Let $\cL := \big\{ f \in \cH$ $\big|$ $T_{D,x} \subset H_f \big\}$. %
Note that $\dim T_{D,x} = \dim D = r \leq N+1$, since $D$ is non-singular. 
The condition $T_{D,x} \subset H_f$ means that 
$f$ passes through independent $r+1$ points. 
Thus, $\dim \cL = \dim \cH - (r+1) = N-r$. 
Since $\cP_x = \cP \cap \cL$, 
we have $\dim \cP_x \leq N-r$. 
\end{proof}

Even if we determine all the discriminants of $\cP$, 
the signature of $\disc(D)$ may not be constant in $\Int(\cP)$. 
To describe $\cP$ as a union of basic semialgebraic sets 
of $\cH$ using some inequalities, 
we need some more inequalities to cut off extra parts or 
to avoid the interior zero locus $\Int(\cP) \cap V_\cH(\disc(D))$. 
Such inequalities are called {\it separators}. 
Note that discriminants are unique up to multiplication by non-zero constant, 
but there may be many possibilities of the choice of separators. 

About extremality of $f \in \cP$, 
the following theorem is useful. 
About the definition of infinitesimal local cone, 
please see \cite[Definition 2.6, 2.12]{RefAd}. 

\def\Tqbae{Theorem 2.8} 
\proclaim{Theorem 2.8} (\cite[Theorem 2.11, Proposition 2.13]{RefAd}) 
{\sl Let $\cP = \cP(A$, $\cH)$. 
Assume that $\dim \cP \geq 2$. }\par
{\parindent=20pt
\Item{\rm(1)} {\sl If $f \in \cE(\cP)$, 
then there exists local cones or infinitesimal local cones 
$\cP_1$,$\ldots$, $\cP_r \subset \cP$ with respect to $f$ 
which satisfy $\cP_1 \cap \cdots \cap \cP_r = \R_+ \cdot f$. }
\Item{\rm (2)} {\sl Let $f \in \cP$. 
If there exists local cones or infinitesimal local cones 
$\cP_1$,$\ldots$, $\cP_r \subset \cP$ such that 
$\cP_1 \cap \cdots \cap \cP_r = \R_+ \cdot f$. 
Then, $f \in \cE(\cP)$. }

}
\endproclaim

In the above theorem, infinitesimal local cones appear for 
special $f \in \cE(\cP)$. 
In ordinary case, 
there exists points ${\bf a}_1$,$\ldots$, ${\bf a}_r \in A$ such that 
\[\R_+ \cdot f = \big\{ g \in \cP \; \big| \; 
 \hbox{$g({\bf a}_1) = \cdots = g({\bf a}_r) = 0$}\big\}.\]
We can choose each ${\bf a}_i$ so that $\Phi_\cH({\bf a}_i) 
 \in \cE(X)$. 
Infinitesimal local cones appears when not less than two zero points of $f$ 
become infinitely near points. 

\def\Tqbbz{Definition 2.9} 
\proclaim{Definition 2.9} 
{\rm Let $\cH$ be a signed linear system on a semialgebraic variety $A$. 
A subset $\Omega \subset A$ is called a {\it test set} for $(A$, $\cH)$, 
if $f({\bf a}) \geq 0$ for all ${\bf a} \in \Omega$, 
then $f({\bf a}) \geq 0$ for all ${\bf a} \in A$. }
\endproclaim

The following theorem will be trivial. 

\def\Tqbby{Theorem 2.10} 
\proclaim{Theorem 2.10} 
{\sl Let $\cH$ be a signed linear system on a 
compact semialgebraic variety $A$ with $\dim \cH \geq 3$, 
and let $X := \Cls(\Phi_{\cH}(A))$ be the characteristic variety. 
Take a subset $\Omega \subset A$. 
If $\cE(X) \subset \Cls(\Phi_{\cH}(\Omega))$, 
then $\Omega$ is a test set for $\cH$. }
\endproclaim

Some artices add the following condition for a test set: \par
(Additional condition) For any ${\bf a} \in \Omega$, 
there exixts $f \in \cH$ such that $f({\bf a}) = 0$. \par
\noindent Under this definition, 
$\cE(X) \subset \Cls(\Phi_{\cH}(\Omega))$ must be 
replaced by $\cE(X) = \Phi_{\cH}(\Omega)$. 

\def\Tqbbb{Example 2.11} 
\proclaim{Example 2.11} 
{\rm Consider the case $A = \P_+^2$, $\cH = \cH_{3,3}^c$. Then 
\[\Omega := \{(1 \colon 1 \colon 1)\} \cup 
   \big\{(0 \colon t \colon 1) \in \P_+^2 \; \big| \; t \geq 0 \big\}\]
is a test set for $\cH_{3,3}^c$ (see \cite[Theorem 3.1]{RefAc}). 
Thus, if $f \in \cH_{3,3}^c$ satisfies 
$f(1$, $1$, $1) \geq 0$ and $f(0$, $t$, $1) \geq 0$ for all $t \geq 0$, 
then $f(a$, $b$, $c) \geq 0$ for all $a$, $b$, $c \in \R_+$. }
\endproclaim

\bigbreak
{\bf 2.2. Some more general theorems.}%
\par\penalty1000\vskip0.4em plus0.1em minus0.1em
Let $V$ and $W$ be non-singular semialgebraic varieties 
with $\dim V = n$, $\dim W = m$, 
and $\varphi \colon V \to W$ be a regular map. 
Take a point $a \in V$ and put $b := \varphi(a)$. 
We can take open neighborhoods $a \in U_V \subset V$ 
and $b \in U_W \subset W$ such that $\varphi(U_V) \subset U_W$ 
and that $U_V$, $U_W$ have local coordinate systems $(x_1$,$\ldots$, $x_n)$ 
and $(y_1$,$\ldots$, $y_m)$ whose origins are $a$, $b$. 
$\varphi$ can be represented by 
functions $y_j = \varphi_j(x_1$,$\ldots$, $x_n)$ ($j=1$,$\ldots$, $m)$. 
Let $\displaystyle J_a 
 := \left.\left(\frac{\partial y_j}{\partial x_i}\right)
\right|_{(x_1,\ldots,x_n)=a}$ be the Jacobian matrix of $\varphi$ at $a$. 
Note that $\rank J_a$ does not depend on the choice 
of $(x_1$,$\ldots$, $x_n)$ and $(y_1$,$\ldots$, $y_m)$. 
We denote 
\[\Sing(\varphi) := \big\{ a \in V \; \big| \; 
\hbox{$\rank J_a < \dim \varphi(V)$}\big\}.\]

\proclaim{Proposition 2.12} 
{\sl If $V$ is a non-singular complete real algebraic variety, 
then $\partial\big(\varphi(V)\big) \allowbreak \subset 
  \varphi\big(\Sing(\varphi)\big)$. }
\endproclaim

\Proof 
Put $r := \dim \varphi(V)$, and assume that $\rank J_a = r$. 
We may assume that 
\[\det \left(\frac{\partial y_j}{\partial x_i}
   \right)_{1 \leq i \leq r, \, 1 \leq j \leq r} \ne 0\]
at $a$. Let $U' := \big\{(x_1$,$\ldots$, $x_n) \in U_V$ $\big|$ 
$x_{r+1} = \cdots = x_n = 0 \big\}$. 
If $U_V$ is sufficiently small Euclidean open set, 
$\varphi\big|_{U'} \colon U' \lto \varphi(U')$ is an isomorphism. 
Thus $b \notin \partial \big(\varphi(V)\big)$. 
\end{proof}


When $V$ and $W$ are open subsets of $\P_{\R}^r$, and $\varphi$ is given by 
$y_j = \varphi_j(x_0 \colon \cdots \colon x_r)$ ($j=0$,$\ldots$, $r)$ 
using homogeneous coordinate system, 
the condition $\rank J_a = r$ can be replaced by 
\[\det \left(\frac{\partial y_j}{\partial x_i}
   \right)_{0 \leq i \leq r, \, 0 \leq j \leq r} \ne 0.\]
When $V$ has singularities, we put $\Sing(\varphi) 
  := \Sing\big(\varphi\big|_{\Reg(V)}\big)$. 

\def\Tqbab{Corollary 2.13} 
\proclaim{Corollary 2.13} 
{\sl Assume that $A$ is a compact semialgebraic variety, then,} 
\[\partial\big(\varphi(A)\big) \subset 
   \varphi\big(\Sing(\varphi) \cup \Sing(A) \cup \partial A\big).\]
\endproclaim

For $f \in \R[x_1$,$\ldots$, $x_n]_d$ and $K = \R$ or $\C$, we denote 
\[V_K(f) := \big\{ a \in \P_K^n \; \big| \; f(a) = 0 \big\}, \quad
  V_+(f) := V_{\R}(f) \cap \P_+^n.\]
In some articles, $V_K(f)$ are also denoted by $\cZ(f)$. 
The symbol $V_K(f)$ is rather popular in algebraic geometry. 

\def\Tqbag{Proposition 2.14} 
\proclaim{Proposition 2.14} 
{\sl Let $X_{3,d}^{s+ }:= X(\P_+^2$, $\cH_{3,d}^s)$. 
If $d \geq 4$, then $X_{3,d}^{s+} \cong \P_+^2/\frS_3$. }
\endproclaim

\Proof 
We denote the coordinate system of $\P_+^2$ by $(a \colon b \colon c)$, 
and put $S_1 := a+b+c$. 
$\Phi_{3,d} := \Phi_{\cH_{3,d}^s} : \P_+^2 \to X_{3,d}^{s+}$ is 
decomposed as 
$\Phi_{3,d} \colon \P_+^2 \mapr{\sigma} \P_+^2/\frS_3 
    \mapr{\Psi_{3,d}} X_{3,d}^{s+}$. 
By \cite[Proposition 2.13, 2.14 and \S 4.5]{RefAc}, 
$\Psi_{3,4} \colon \P_+^2/\frS_3 \lto X_{3,4}^{s+}$ is an isomorphism. 
Since $V_+(S_1) \cap \P_+^2 = \emptyset$, the multiplication map 
$\times S_1 \colon \cH_{s,d}^s \lto \cH_{s,d+1}^s$ induces 
an isomorphism $X_{3,d+1}^{s+} \to X_{3,d}^{s+}$. 
\end{proof}

\bigskip

In the cyclic case $X_{n,d}^{c+} := X(\P_+^{n-1}$, $\cH_{n,d}^c)$, 
we know that $X_{n,d}^{c+} \cong \P_+^{n-1}/\frC_n$ if $d \geq n$, 
here $\frC_n = \Z/n\Z$ (see \cite[Proposition 1.36]{RefAc}). 
When $n=3$, $\Delta^1(X_{3,d}^{c+})$ has a unique element 
$C_{3,d}^{c+} := \big\{ \Phi_{3,d}^c(0 \colon s \colon 1)$ 
 $\big|$ $s>0 \big\}.$ 
We call $\disc(C_{3,d}^{c+})$ to be the {\it edge discriminant} 
of $\cP_{3,d}^{c+}$ (see \cite[Definition 2.7]{RefAc}). 
The following Theorem is a replacement of 
\cite[Proposition 2.10, Theorem 5.9, 6.8]{RefAc}. 

We denote the discriminant of $c_nx^n+c_{n-1}x^{n-1}+\cdots+c_1x+c_0$ by 
\[\Disc_n(c_n,c_{n-1},\ldots,c_1,c_0).\]

\def\Tqbaiz{Theorem 2.15} 
\proclaim{Theorem 2.15} 
{\sl Let's denote the coordinate system 
of $\P_+^2$ by $(a \colon b \colon c)$, 
and put $S_{m,n} = S_{m,n}(a,b,c) := a^m b^n + b^m c^n + c^m a^n$, 
$S_n = S_n(a,b,c) := a^n+b^n+c^n$, and $U = U(a,b,c) := abc$. 
Take the basis of $\cH_{3,d}^c$ so that $s_0 = S_d$, 
$s_1 = S_{d-1,1}$, $s_2 = S_{d-2,2}$,$\ldots$, $s_{d-1} = S_{1,d-1}$,$\ldots$. 
Here, if $i \geq d$, then $s_i$ is a multiple of $abc$. 
We represent $f \in \cH_{3,d}^c$ as $f = \sum p_i s_i$. 
Then, the edge discriminant of $\cP_{3,d}^{c+}$ agrees 
with $\Disc_d(p_0,p_1,\ldots,p_{d-1},p_0)$. }
\endproclaim

\Proof 
Let $\cL_{0,t}^{c+}$ be the local cone of $\cP_{3,s}^{c+}$ 
at $(0 \colon t \colon 1) \in \P_+^2$. 
Take $f \in \cL_{0,t}^{c+} \subset \cF(C_{n,d}^+)$ 
($p_0>0$ and $t>0$). 
Then $f(0,t,1)=0$. 
Since $f(0,x,1) \geq 0$ for all $x>0$, 
the equation $f(0,x,1)=0$ has a multiple root at $x=t$. 
Thus, the discriminant of $f$ is equal to $0$. 
Since $S_{i,d-1}(0,x,1) = x^i$ ($1 \leq i \leq d-1$), 
$S_d(0,x,1) = x^d + 1$ and $U(0,x,1)=0$, we have 
$f(0,x,1) = p_0 x^d + p_1 x^{d-1} + \cdots + p_{d-1} x + p_0$. 

Since $\Disc_d$ and $\disc(C_{3,d}^{c+})$ are irreducible, 
we have the conclusion. 
\end{proof}

\def\Tqbaf{Theorem 2.16} 
\proclaim{Theorem 2.16} 
{\sl Consider the cases $A = \P_{\R}^{n-1}$ or $\P_+^{n-1}$, 
and $\cH = \cH_{n,d}^s$ or $\cH_{n,d}^{s0}$. 
Let $\cP := \cP(A$, $\cH)$, $X := X(A$, $\cH_{n,d})$, 
and $\Phi := \Phi_{\cH} \colon A \cdots\to X$. 
Let $\sigma \colon \P_{\R}^{n-1} \lto \P_{\R}^{n-1}/\frS_n 
\subset \P_{\R}(1,2,\ldots,n)$ be the natural surjection, 
and $\Psi \colon \P_{\R}^{n-1}/\frS_n \cdots\to X$ be the 
rational map such that $\Psi \circ \pi = \Psi$. 
Assume that $\Psi$ is a birational map. 
Let $D \in \Delta^r(X)$ with $r \geq \max \{2$, $\lfloor d/2 \rfloor\}$. 
Then $\cF(D)$ is not a face component of $\cP$. }
\endproclaim

\Proof 
Let $r_0 := \max \{2$, $\lfloor d/2 \rfloor\}$, and take 
$D \in \Delta^r(X)$ with $r_0 \leq r \leq n-1$. 
Assume that $\cF(D)$ is a face component of $\cP$. 
Then $\dim \cF(D) = n-1$. 

\smallskip

(1) Consider the case $A = \P_{\R}^{n-1}$. 

Let $\Omega := \big\{(x_1 \colon \cdots \colon x_n) \in \P_{\R}^{n-1}$ 
$\big|$ $\# \{x_1$,$\ldots$, $x_n\} \leq r_0\big\}$. 
Here $\# \{x_1$,$\ldots$, $x_n\} \leq r_0$ means that 
at most $r_0$ members of $x_1$,$\ldots$, $x_n$ are distinct. 
$\Omega$ is a test set by \cite{RefRie}. 

$\Omega$ is included in a union of some $(r_0-1)$-dimensional 
linear subspace of $\P_{\R}^{n-1}$. 
Take general $f \in \cF(D)$. 
There exists a semialgebraic subset $E \subset A$ such that $\Phi(E) = D$, 
and ${\bf a} \in E$ such that $f({\bf a}) = 0$. 
Since $\cF(D)$ is a face component, 
we may assume that the hyperplane $H_f \subset \P(\cH^{\vee})$ 
corresponding to $f$, tangents to $X$ only 
at the unique point $\Phi({\bf a})$. 
This means that if ${\bf b} \in A - \Bs \cH$ satisfies $f({\bf b}) = 0$, 
then $\Phi({\bf b}) = \Phi({\bf a})$, if ${\bf a} \in E$ is a general point. 
We can choose such $f$ and ${\bf a}$. 

By \cite[Corollary 1.3]{RefRie} \ or \cite[Corollary 2.1]{RefTa}, 
there exists ${\bf b} \in \Omega$ such that $f({\bf b}) = 0$. 
We denote this ${\bf b}$ by ${\bf b}({\bf a})$. 
${\bf a}$ can move a certain $r$-dimensional subset of $E$. 
But $\dim \Omega = r_0-1 < r$. 
Thus, there exists ${\bf a} \in E$ such that 
$\Phi({\bf b}({\bf a})) \ne \Phi({\bf a})$. 
A contradiction. 
Thus $\cF(D)$ is not a face component of $\cP$. 

\smallskip

(2) Consider the case $A = \P_+^{n-1}$. 

Let $\Omega' := \big\{(x_1 \colon \cdots \colon x_n) \in \P_+^{n-1}$ 
$\big|$ $\# \big(\{x_1$,$\ldots$, $x_n\}-\{0\}\big) \leq r_0\big\}$. 
If $f \in \cH_{n,d}^s$ satisfies $f({\bf a}) \geq 0$ 
for all ${\bf a} \in \Omega'$, then $f \in \cP_{n,d}^{s+}$ by \cite{RefRie}. 
$\Omega'$ is also included in a union of some $(r_0-1)$-dimensional 
linear subspace of $\P_{\R}^{n-1}$. 

The left part is same as (1). 
\end{proof}

\bigskip

If $\cF(D)$ is not a face component, then, 
for each $f \in \cF(D)$, 
there exist $D_1$,$\ldots$, $D_r \in \Delta(X) - \{D\}$ such that 
$f \in \cF(D_1) \cap \cdots \cap \cF(D_r)$, 
and that all $\cF(D_i)$ are face components. 

\removelastskip\penalty-400\vskip2.5em plus0.3em minus0.3em
{\bf Section 3. Quartic Inequalities of Four Variables} 
\par\penalty1000\vskip0.8em plus0.2em minus0.2em
In this section, we shall study 
$\cP_{4,4}^{s0}$ and $\cP_{4,4}^{s0+}$. 
We write the homogeneous coordinate system of $A = \P_{\R}^3$ or $A = \P_+^3$ 
by $(a \colon b \colon c \colon d)$ 
or $(a_0 \colon a_1 \colon a_2 \colon a_3)$. 
We regard $a_{4n+i} = a_i$ for $n \in \Z$. 
We denote 
\begin{align*}
 & \displaystyle S_d := \sum_{i=0}^3 a_i^d, \quad  
   \displaystyle T_{p,q} := \sum_{i=0}^3 a_i^p(a_{i+1}^q+a_{i+2}^q+a_{i+3}^q),
   \quad
   \displaystyle S_{p,p} := \sum_{0 \leq i < j \leq 3} a_i^p a_j^p, \\
 & \displaystyle T_{p,q,q} := \sum_{i=0}^3 a_i^p(a_{i+1}^q a_{i+2}^q
     + a_{i+1}^q a_{i+3}^q + a_{i+2}^q a_{i+3}^q), \quad
   U := a_0a_1a_2a_3. 
\end{align*}
A polynomial $f \in \cH_{n,d}^s$ or $\cH_{n,d}^c$ is 
called {\it monic}, if the coefficient of 
$S_d = a_0^d + \cdots + a_{n-1}^d$ is equal to $1$. 
For a subset $V \subset \cH_{n,d}^c$, we denote 
\[\breve{V} := \big\{ f \in V \; \big| \; \hbox{$f$ is monic}\big\}.\]
We denote as $\P_{\R}^n:(a_0 \colon \cdots \colon a_n)$ 
when we treat $\P_{\R}^n$ with a 
homogeneous coordinate system $(a_0 \colon \cdots \colon a_n)$. 
Similarly we denote as $\R^n:(x_1,\ldots,x_n)$ when 
we study $\R^n$ with a coordinate system $(x_1$,$\ldots$, $x_n)$. 

\par\penalty1000\vskip1.0em plus0.3em minus0.3em
{\bf 3.1. Structure of $\P_{\R}^3/\frS_4$}
\par\penalty1000\vskip0.4em plus0.1em minus0.1em
Let $(a_0 \colon \cdots \colon a_n)$ be the homogeneous coordinate system 
of $\P_{\R}^n$, and $\sigma_k \allowbreak
= \sigma_k(a_0$,$\ldots$, $a_n)$ be the $k$-th 
symmetric function of $a_0$,$\ldots$, $a_n$ ($0 \leq k \leq n+1$). 
The sequence of functions $(\sigma_1$,$\ldots$, $\sigma_{n+1})$ defines 
the regular map $\sigma \colon \P_{\P}^n \lto \P_{\R}(1,2,\ldots,n+1)$, 
where $\P_{\R}(1,2,\ldots,n+1)$ is the real weighted projective space 
which is defined as the real part of the complex weighted projective 
space $\P_{\C}(1,2,\ldots$, $n+1)$. 
The image $\sigma(\P_{\R}^n)$ is isomorphic to $\P_{\R}^n/\frS_{n+1}$ 
as semialgebraic varieties. 
Note that $\P_{\C}^n/\frS_{n+1} \cong \P_{\C}(1,2,\ldots,n+1)$, 
but $\P_{\R}^n/\frS_{n+1} \not\cong \P_{\R}(1,2,\ldots,n+1)$. 
In general, for two points $P$, $Q \in \P_{\R}^n$, 
$(PQ)$ represents an open line segment, 
$[PQ] := (PQ) \cup \{P$, $Q\}$ represents a closed line segment, 
and $PQ$ represents a line. 

\proclaim{Definition 3.1} 
{\rm Assume that a finite group $G$ acts on a semialgebraic variety $A$. 
Let $\sigma \colon A \to A/G$ be the natural surjection. 
A closed semialgebraic subset $A_0 \subset A$ is called a 
fundamental domain of $A/G$, 
if $\sigma(A_0) = A/G$ and $\sigma : \Int(A_0) \lto \sigma(\Int(A_0)) 
 \subset A/G$ is an isomorphism. }
\endproclaim

\def\Tqccn{Example 3.2} 
\proclaim{Lemma 3.2} 
{\rm (1) Let $A = \P_{\R}^n$ and $G = \Z/(n+1)\Z$. 
Then $(\P_{\R}^n)^G = \{{\bf 1}\}$, 
and $\Sing(\P_{\R}^n/G) = \sigma\big((\P_{\R}^n)^G\big) 
= \{\sigma({\bf 1})\}$, 
here ${\bf 1} = (1 \colon 1 \colon \cdots \colon 1) \in A$. 
The following $A_c$ is a fundamental domain. 
\[A_c := \left\{ (a_0 \colon \cdots \colon a_{n-1} \colon 1) \in \P_{\R}^n \;
   \left| \; 
   \vcenter{\hbox{$a_0+a_1+\cdots+a_{n-1}+1 \geq 0$,}
            \hbox{$a_0 \leq 1$, $a_1 \leq 1$,$\ldots$, $a_{n-1} \leq 1$}}\ 
   \right.\right\}.\]

(2) Let $A = \P_+^n$ and $G = \Z/(n+1)\Z$. 
Then $(\P_{\R}^n)^G = \{{\bf 1}\}$, and 
\[A_c^+ := \big\{ (a_0 \colon \cdots \colon a_{n-1} \colon 1) 
 \in \P_{\R}^n \; \big| \; 
 \hbox{$0 \leq a_0 \leq 1$,$\ldots$, $0 \leq a_{n-1} \leq 1$} \big\}\]
is a fundamental domain. 

(3) Let $A = \P_{\R}^n$ and $G = \frS_{n+1}$. Then 
\[A_s := \big\{ (a_0 \colon  \cdots \colon a_{n-1} \colon 1) \in A_c \; \big| 
 \; \hbox{$a_0 \leq a_1 \leq \cdots \leq a_{n-1}$, 
     $a_0 + \cdots + a_{n-1} + 1 \geq 0$}\big\}\]
is a fundamental domain. 

(4) Let $A = \P_+^n$ and $G = \frS_{n+1}$. Then 
\[A_s^+ := \big\{ (a_0 \colon \cdots \colon a_{n-1} \colon 1) \in \P_{\R}^n 
 \; \big| \; 
 \hbox{$0 \leq a_0 \leq a_1 \leq \cdots \leq a_{n-1} \leq 1$} \big\}\]
is a fundamental domain. }

\bigskip

Note that $\P_{\C}^3/\frS_4 \cong \P_{\C}(1,2,3,4)$ has 
cyclic quotients singularities 
at $\tilde{P}_0:=(0 \colon 1 \colon 0 \colon 0)$, 
$\tilde{P}_0':=(0 \colon 0 \colon 1 \colon 0)$ 
and $\tilde{P}_0'':=(0 \colon 0 \colon 0 \colon 1)$. 

\def\Tqca{Proposition 3.3}
\proclaim{Proposition 3.3} 
{\sl About the structures of $\P_{\R}^3/\frS_4$ 
and $\P_+^3/\frS_4$, we have the followings: }
{\parindent=20pt
\Item{\rm (1)} {\sl Let $\sigma \colon \P_{\R}^3 \lto \P_{\R}^3/\frS_4 
 \mapr{\subset} \P_{\R}(1,2,3,4)$ be the natural map. 
Then $\sigma^{-1}(\tilde{P}_0') = \emptyset$, 
$\sigma^{-1}(\tilde{P}_0'') = \emptyset$, and }
$\sigma(-1$, $0$, $0$, $1) = \tilde{P}_0$. 
\Item{\rm (2)} {\sl 
$\Delta^2(\P_{\R}^3/\frS_4) = \big\{\tilde{D}_1\big\}$, 
$\Delta^1(\P_{\R}^3/\frS_4) = \big\{\tilde{C}_1$, $\tilde{C}_2 \big\}$, 
and $\Delta^0(\P_{\R}^3/\frS_4) \allowbreak
   = \big\{\tilde{P}_0$, $\tilde{P}_1$, $\tilde{P}_2 \big\}$, 
where $\tilde{D}_1$, $\tilde{C}_i$ and $\tilde{P}_i$ are as follows:} 
\begin{align*}
 & \tilde{D}_1 := \big\{ \sigma(s \colon t \colon 1 \colon 1) 
    \in \P_{\R}(1,2,3,4) \; \big| \; 
        \hbox{$s+t+2 > 0$, $s < t$, $s \ne 1$, $t \ne 1$} \big\}, \\
 & \tilde{C}_1 := \big\{ \sigma(s \colon 1 \colon 1 \colon 1) 
    \in \P_{\R}(1,2,3,4) \; \big| \; 
        \hbox{$s \in \P_{\R}^1$, $s \ne -3$, $1$} \big\}, \\
 & \tilde{C}_2 := \big\{ \sigma(s \colon s \colon 1 \colon 1) 
    \in \P_{\R}(1,2,3,4) \; \big| \; \hbox{$-1<s<1$} \big\}, \\
 & \tilde{P}_1 := \sigma(1 \colon 1 \colon 1 \colon 1) 
    = (4 \colon 6 \colon 4 \colon 1) \in \P_{\R}(1,2,3,4), \\
 & \tilde{P}_2 := \sigma(-1 \colon -1 \colon 1 \colon 1) 
    = (0 \colon -2 \colon 0 \colon 1) \in \P_{\R}(1,2,3,4). 
\end{align*}
\Item{\rm (3)} {\sl $\Delta^2(\P_+^3/\frS_4)
   = \big\{\tilde{D}_1^+$, $\tilde{D}_0\big\}$, 
$\Delta^1(\P_+^3/\frS_4) 
  = \big\{\tilde{C}_1^+$, $\tilde{C}_2^+$, $\tilde{C}_3$, $\tilde{C}_4\big\}$, 
and $\Delta^0(\P_+^3/\frS_4) \allowbreak
   = \big\{\tilde{P}_1$, $\tilde{P}_3$, $\tilde{P}_4$, $\tilde{P}_5 \big\}$, 
where $\tilde{D}_1^+$, $\tilde{D}_0$, $\tilde{C}_1'$, $\tilde{C}_i$ 
and $\tilde{P}_i$ are as follows:} 
\begin{align*}
 & \tilde{D}_1^+ 
     := \big\{ \sigma(s \colon t \colon 1 \colon 1) \in \P_{\R}(1,2,3,4) 
        \; \big| \; \hbox{$0 < s < t$, $s \ne 1$, $t \ne 1$} \big\}, \\
 & \tilde{D}_0 
     := \big\{ \sigma(0 \colon s \colon t \colon 1) \in \P_{\R}(1,2,3,4) 
        \; \big| \; \hbox{$0<s<t<1$} \big\}, \\
 & \tilde{C}_1^+ 
     := \big\{ \sigma(s \colon 1 \colon 1 \colon 1) \in \P_{\R}(1,2,3,4) 
        \; \big| \; \hbox{$0<s<1$ or $s>1$} \big\}, \\
 & \tilde{C}_2^+ 
     := \big\{ \sigma(s \colon s \colon 1 \colon 1) \in \P_{\R}(1,2,3,4) 
        \; \big| \; \hbox{$0<s<1$} \big\}, \\
 & \tilde{C}_3 
     := \big\{ \sigma(0 \colon s \colon 1 \colon 1) \in \P_{\R}(1,2,3,4) 
        \; \big| \; \hbox{$0<s<1$ or $1<s$} \big\}, \\
 & \tilde{C}_4 
     := \big\{ \sigma(0 \colon 0 \colon s \colon 1) \in \P_{\R}(1,2,3,4) 
        \; \big| \; \hbox{$0<s<1$} \big\}, \\
 & \tilde{P}_3 := \sigma(0 \colon 1 \colon 1 \colon 1) 
      = (3 \colon 3 \colon 1 \colon 0) \in \P_{\R}(1,2,3,4), \\
 & \tilde{P}_4 := \sigma(0 \colon 0 \colon 1 \colon 1) 
      = (2 \colon 1 \colon 0 \colon 0) \in \P_{\R}(1,2,3,4), \\
 & \tilde{P}_5 := \sigma(0 \colon 0 \colon 0 \colon 1) 
      = (1 \colon 0 \colon 0 \colon 0) \in \P_{\R}(1,2,3,4). 
\end{align*}
\Item{\rm (4)} {\sl $\disc(\tilde{D}_1)=\Disc_4$, 
and $\tilde{C}_1 \cup \tilde{C}_2 \subset \Sing(V(\Disc_4))$, 
here $V(f)$ is the zero locus of $f$ in $\P_{\R}(1,2,3,4)$. }
\Item{\rm (5)} {\sl $\Cls{\tilde{C}_1}$ is isomorphic to a cubic curve 
on $\P_{\R}^2$ with a cusp at $\tilde{P}_1$.} 
\Item{\rm (6)} {\sl $\tilde{C}_2 = (\tilde{P}_1\tilde{P}_2)$ is isomorphic to 
an open line segment with ends $\tilde{P}_1$ and $\tilde{P}_2$. } 
\Item{\rm (7)} {\sl $\P_{\R}^3/\frS_4$ is the semialgebraic subset 
of $\P_{\R}(1,2,3,4)$ defined by 
$\Disc_4(1$, $\sigma_1$, $\sigma_2$, $\sigma_3$, $\sigma_4) \geq 0$, 
$8\sigma_2 \leq 3\sigma_1^2$, and 
$64\sigma_4 - 16\sigma_2^2 + 16\sigma_1^2\sigma_2 
  - 16\sigma_1\sigma_3 - 3\sigma_1^4 \leq 0$. 
Here, $\sigma_i$ is the elementary symmetric polynomials of $a_0$, $a_1$, 
$a_2$, $a_3$ of degree $i$. }

}
\endproclaim

\Proof 
(1) is clear. 

(2) and (3) follow from the critical decompositions of 
fundamental domains $A_s$ and $A_s^+$ in the above example. 

(4) This follows from conditions that a quartic equation has a double root, 
a triple root or two double roots. 

(5) Eliminate $t$ from 
$\displaystyle x = \frac{\sigma_2(t,1,1,1)}{\sigma_1(t,1,1,1)^2}$, 
$\displaystyle y = \frac{\sigma_3(t,1,1,1)}{\sigma_1(t,1,1,1)^3}$, 
$\displaystyle z = \frac{\sigma_4(t,1,1,1)}{\sigma_1(t,1,1,1)^4}$, 
then we obtain 
\[32\left(x-\frac{3}{8}\right)^3+27\left(x-\frac{3}{8}\right)^2
   -108\left(x-\frac{3}{8}\right)\left(y-\frac{1}{16}\right)
   +108\left(y-\frac{1}{16}\right)^2 = 0\]
and $x^2 = 3 y - 12 z$. 
This cuve is isomorphic to a cubic curve on $\P_{\R}^2$, and have a cusp at 
$(x$, $y$, $z) = (3/8$, $1/16$, $1/256) = \tilde{P}_1$. 

(6) Eliminate $t$ from 
$\displaystyle x = \frac{\sigma_2(t,t,1,1)}{\sigma_1(t,t,1,1)^2}$, 
$\displaystyle y = \frac{\sigma_3(t,t,1,1)}{\sigma_1(t,t,1,1)^3}$, 
$\displaystyle z = \frac{\sigma_4(t,t,1,1)}{\sigma_1(t,t,1,1)^4}$, 
then we obtain 
$4x - 8y = 1$ and $y^2 = z$. This is a non-singular rational curve. 

(7) This follows from theory of quartic equations. 
$g(a,b,c,d) : = 64\sigma_4 - 16\sigma_2^2 + 16\sigma_1^2\sigma_2 
  - 16\sigma_1\sigma_3 - 3\sigma_1^4$ is a separator. Note that 
\begin{align*}
 & g(a,a,c,d) = -(c-d)^2 (8a^2 - 8a c + 3c^2 - 8a d + 2c d + 3d^2), \\
 & g(a,a,a,d) = -3(a-d)^4. 
\end{align*}
Thus, $V(g)$ pass through  $\tilde{C}_2$. 
\end{proof}

\bigbreak
{\bf 3.2 The PSD cone $\cP_{4,4}^{s0}$}
\par\penalty1000\vskip0.4em plus0.1em minus0.1em
In this subsection, we shall study 
$\cP_{4,4}^{s0} := \cP(\P_{\R}^3$, $\cH_{4,4}^{s0})$. 
We choose 
\[s_0 := S_4-4U, \quad s_1 := T_{3,1}-12U, \quad s_2 := S_{2,2}-6U, \quad 
s_3 := T_{2,1,1} - 12U\]
as a basis of $\cH_{4,4}^{s0}$. 
We define $\Phi_{4,4}^{s0} : \P_{\R}^3 \cdots\to \P_{\R}^3$ by 
$\Phi_{4,4}^{s0}(a) = \big(s_0(a):s_1(a):s_2(s):s_3(a)\big)$. 
Let 
\[X_{4,4}^{s0} := \Phi_{4,4}^{s0}(\P_{\R}^3)
  = X(\P_{\R}^3, \, \cH_{4,4}^{s0}) 
    \subset \P\big((\cH_{4,4}^{s0})^{\vee}\big)),\]
and let $\Psi \colon \P_{\R}^3/\frS_4 \cdots\to X_{4,4}^{s0}$ be 
the rational map such that $\Phi_{4,4}^{s0} = \Psi \circ \sigma$. 
Moreover, let 
\begin{align*}
 & C_1 := \Cls(\Psi(\tilde{C}_1))
       = \big\{\Phi_{4,4}^{s0}(t \colon 1 \colon 1 \colon 1) 
            \; \big| \; \hbox{$t \in \P_{\R}^1$} \big\}, \\
 & C_2 := \Psi(\tilde{C}_2)
       = \big\{\Phi_{4,4}^{s0}(t \colon t \colon 1 \colon 1) 
            \; \big| \; \hbox{$-1 < t < 1$} \big\}, \\
 & P_0 := \Psi(\tilde{P}_0) = \Phi_{4,4}^{s0}(-1 \colon 0 \colon 0 \colon 1)
    = (2 \colon -2 \colon 1 \colon 0), \\
 & P_2 := \Psi(\tilde{P}_2) = \Phi_{4,4}^{s0}(-1 \colon -1 \colon 1 \colon 1)
    = (0 \colon 1 \colon 0 \colon 1), \\
 & P_{-3} := \Phi_{4,4}^{s0}(-3 \colon 1 \colon 1 \colon 1) 
    = (2 \colon -1 \colon 1 \colon 1), \\
 & E_0 := \big\{ (a \colon b \colon c \colon d) \in \P_{\R}^4 
     \; \big| \; \hbox{$a$, $b$, $c$, $d \in \R$, $a+b+c+d=0$} \big\}, \\
 & L_0 := \Phi_{4,4}^{s0}(E_0), \\
 & D_1' := \big\{ (a \colon b \colon c \colon c) \in \P_{\R}^4 
     \; \big| \; \hbox{$a$, $b$, $c \in \R$} \big\}, \\
 & D_1 := \Psi(\tilde{D}_1) = \Phi_{4,4}^{s0}(D_1').  
\end{align*}

Note that $\Bs \cH_{4,4}^{s0} = \{(1 \colon 1 \colon 1 \colon 1)\}$, 
and $\Psi$ is undefinite at $\tilde{P}_1$. 

The aim of this subsection is to prove the following theorem. 

\def\Tqcb{Theorem 3.4}
\proclaim{Theorem 3.4} 
{\rm (1)} {\sl For a monic $f = s_0 + p s_1 + q s_2 + r s_3 
  \in \breve{\cH}_{4,4}^{s0}$, 
$f({\bf a}) \geq 0$ for all ${\bf a} \in \R^4$ if and only if }
\[p + r \geq 0 \quad \hbox{and} \quad
  - 9p^2 + 12p + 12q + 12r + 8 \geq 0.\]
{\parindent=20pt
\Item{{\rm (2)}} {\sl All the extremal elements of $\cP_{4,4}^{s0}$ are 
positive multiples of $\frg_t$ 
{\rm ($t \in \P_{\R}^1 = \R \cup \{\infty\}$)} or $\frp$.}
\Item{{\rm (3)}} {\sl All the discriminants of $\cP_{4,4}^{s0}$ are 
$\disc_{C_1} = 9p^2 + 12p + 12q + 12r + 8$ and $\disc_{P_2} = p + r$.}
\Item{{\rm (4)}} {\sl 
$\big\{(t \colon 1 \colon 1 \colon 1) \in \P_+^3 \; \big| \; 
   \hbox{$t \geq 0$}\big\} \cup \{(-1 \colon -1 \colon 1 \colon 1)\}$ 
is a test set for $\cP_{4,4}^{s0}$. }

}
\endproclaim

This theorem will be proved after \Tqch. 

\def\Tqccea{Lemma 3.5}%
\proclaim{Lemma 3.5} 
{\parindent=20pt
\Item{(1)} {\sl $\Psi : \P_{\R}^3/\frS_4 \lto X_{4,4}^{s0}$ is continuous 
map and $\Psi : (\P_{\R}^3/\frS_4 - \{\tilde{P}_1\}) \lto X_{4,4}^{s0}$ is 
a birational morphism. 
All the exceptinal set of $\Phi_{4,4}^{s0} : \P_{\R}^3 \cdots\to \P_{\R}^3$ 
is $E_0$. }
\Item{(2)} {\sl $\partial \cP_{4,4}^{s0} = \cF(C_1) \cup \cF(P_2)$, 
and $\cE(X_{4,4}^{s0}) \subset C_1 \cup \{P_2\}$. }

}
\endproclaim

\Proof 
We denote the coordinate system of $\P\big((\cH_{4,4}^s)^{\vee}\big)
 = \P_{\R}^3$ by $(x_0 \colon x_1 \colon x_2 \colon x_3)$. 
$\Phi_{4,4}^{s0}$ is defined by $x_i = s_i({\bf a})$. 

(0) Let $P_1 := (2 \colon 3 \colon 1 \colon 1)$. 
When $a$, $b$, $c \to 0$, 
\begin{align*}
\Phi_{4,4}^{s0}(1 \colon 1+a \colon 1+b \colon 1+c)
 &  = (3a^2-2ab+3b^2-2ac-2bc+3c^2)(2 \colon 3 \colon 1 \colon 1) \\
 & \hskip50pt + (\hbox{higher terms of $a$, $b$, $c$}).
\end{align*}
Thus $\Phi_{4,4}^{s0}(1 \colon 1 \colon 1 \colon 1) 
= \Psi(\tilde{P_1}) = P_1$, and $\Psi$ is continuoius at $\tilde{P_1}$. 

\smallskip

(1) We take $A_s$ as \Tqccn (3). 
It is easy to see that $\Phi_{4,4}^{s0} \otimes_{\R} \C : \P_{\C}^3 
  \lto \P_{\C}^3$ is a generically finite rational map of degree $24$. 
Thus $\Phi_{4,4}^{s0} \colon A_s \lto X_{4,4}^{s0}$ is generically one to one. 
Using PC, we have 
\[J_P := \det \left(
  \frac{\partial s_i(a_0,a_1,a_2,a_3)}{\partial a_j}\right)_{0\leq i,j \leq 3}
  = 16 S_1^2 (3S_2 - 2S_{1,1})^2 \prod_{i<j} (a_i-a_j).\]
$J_p \ne 0$ on $\Int(A_s) - \{(1 \colon 1 \colon 1 \colon 1)\}$. 
Thus $\Phi_{4,4}^{s0} \colon \Int(A_s) \lto X_{4,4}^{s0}$ is injective. 
Since $\displaystyle \partial A_s 
  \subset E_0 \cup \bigcup_{\tau \in \frS_4} \tau(D_1')$, 
we have $\partial X_{4,4}^{s0} = L_0 \cup D_1$. 
So, $\Psi : (\P_{\R}^3/\frS_4 - \{\tilde{P}_1\}) 
\lto X_{4,4}^{s0}$ is a birational morphism. 
Let 
\begin{align*}
 & f_{4,4}^{s0}(x_0,x_1,x_2,x_3) := 
    -3 x_0 x_1^4 + 4 x_1^5 + 6 x_0^2 x_1^2 x_2 - 24 x_0 x_1^3 x_2 
    + 14 x_1^4 x_2 - 3 x_0^3 x_2^2 \\
 & \hskip20pt + 20 x_0^2 x_1 x_2^2 - 48 x_0 x_1^2 x_2^2 + 16 x_1^3 x_2^2 
    + 34 x_0^2 x_2^3 + 16 x_0 x_1 x_2^3 + 8 x_1^2 x_2^3 + 44 x_0 x_2^4 \\
 & \hskip20pt 
    - 48 x_1 x_2^4 - 72 x_2^5 + 12 x_0^2 x_1^2 x_3 + 12 x_0 x_1^3 x_3 
    - 36 x_1^4 x_3 - 12 x_0^3 x_2 x_3 + 20 x_0^2 x_1 x_2 x_3 \\
 & \hskip20pt 
     + 120 x_0 x_1^2 x_2 x_3 - 56 x_1^3 x_2 x_3 - 76 x_0^2 x_2^2 x_3 
     - 32 x_0 x_1 x_2^2 x_3 - 64 x_1^2 x_2^2 x_3 \\
 & \hskip20pt 
    - 32 x_0 x_2^3 x_3 + 112 x_1 x_2^3 x_3 + 144 x_2^4 x_3 - 12 x_0^3 x_3^2 
    - 40 x_0^2 x_1 x_3^2 - 112 x_1 x_2^3 x_3 \\
 & \hskip20pt 
    + 144 x_2^4 x_3 - 12 x_0^3 x_3^2 - 40 x_0^2 x_1 x_3^2 
    - 18 x_0 x_1^2 x_3^2 + 104 x_1^3 x_3^2 + 14 x_0^2 x_2 x_3^2 \\
 & \hskip20pt 
    - 104 x_0 x_1 x_2 x_3^2  + 84 x_1^2 x_2 x_3^2 + 64 x_0 x_2^2 x_3^2 
    + 16 x_1 x_2^2 x_3^2 - 152 x_2^3 x_3^2 + 28 x_0^2 x_3^3 \\
 & \hskip20pt 
    + 12 x_0 x_1 x_3^3 - 136 x_1^2 x_3^3 + 8 x_0 x_2 x_3^3 
    - 56 x_1 x_2 x_3^3 + 32 x_2^2 x_3^3 - 3 x_0 x_3^4 + 84 x_1 x_3^4 \\
 & \hskip20pt 
    + 14 x_2 x_3^4 - 20 x_3^5. 
\end{align*}
Since 
\[f_{4,4}^{s0}(s_0,s_1,s_2,s_3) 
    = 16 (a_0+a_1+a_2+a_3)^4 \left(\prod_{i<j} (a_i - a_j)^2 \right) 
       \left(\sum_{i<j} (a_i-a_j)^2\right)^2,\]
we have $\partial X_{4,4}^{s0} = \Phi_{4,4}^{s0}(E_0 \cup D_1')
   \subset V_{\R}(f_{4,4}^{s0}) \subset \P_{\R}^3$ by \Tqbab. 
Since $f_{4,4}^{s0}$ is irreducible, 
we have $\Zar(\partial X_{4,4}^{s0}) = V_{\R}(f_{4,4}^{s0})$. 
Note that $f_{4,4}^{s0} \geq 0$ on $X_{4,4}^{s0}$. 

It is easy to see that $L_0$ is a closed line segment $[P_2P_{-3}]$ defined 
by $x_0 = 2 x_2$, $x_0 - x_1 + x_3 = 0$ and $x_1/x_0 \leq -1/2$. 
This also means that $E_0$ is an exceptinal set of $\Phi_{4,4}^{s0}$. 

Similarlym $C_2$ is an open line segment $(P_1P_2)$ defined 
by $x_0 = 2 x_2$, $x_0 - x_1 + x_3 = 0$ and $x_1/x_0 < 3/2$. 
Note that $L_0$, $C_2 \subset \partial X_{4,4}^{s0}$. 

Next we consider $C_1$. 
Let 
\[g_2(x_0, x_1, x_2, x_3) 
   := (x_1-x_3)^2 + 2 x_2^2 - 3 x_2 x_0.\]
Then $C_1$ is the conic defined 
by $x_2=x_3$ and $g_2(x_0$, $x_1$, $x_2$, $x_3) = 0$. 
Note that $x_2 - x_3 \geq 0$ on $X_{4,4}^{s0}$, 
because $s_2-s_3 = \frp \geq 0$ on $A_s$. 

Let $B$ be the ellise domain on the plane $x_2=x_3$ defined 
by $g_2(1, x_1, x_2, x_2) \leq 0$, 
and let $Y$ be the cone with the base $B$ and the virtex $P_2$. 

\smallskip

(2) We shall show that $Y$ is the convex closure of $X_{4,4}^{s0}$. 

A point on $C_1$ can be written as 
\[P(t) = \Phi_{4,4}^{s0}(t,1,1,1) 
   = (t^2+2t+3 \colon 3(t+2) \colon 3 \colon 3)\]
where $t \in \P_{\R}^1$. 
$P(1) = P_1$ and $P(-3) = P_{-3}$. 
Let $L(t) := (P_2 P(t))$ be an open line segment. 
Note that $L(-3) = (P_2P_{-3}) \subset L_0$, and $L(1) = (P_2P_1) = C_2$. 
A point on $L(t)$ can be written as 
$P(t,s) = P(t) + s P_2$ by $s > 0$. 
Using PC, we have 
\[f_{4,4}^{s0}(P(t,s)) = -12 s^2(s-1)^2(t+3)^4.\]
This implies $L(t) \cap X_{4,4}^{s0} = \emptyset$, if $t \ne -3$, $1$. 
This means $X_{4,4}^{s0} \subset Y$. 
Since $C_1 \cup \{P_2\} \subset X_{4,4}^{s0} \subset Y$, 
we conclude that $Y$ is the convex closure of $X_{4,4}^{s0}$. 
This also implies $\cE(X_{4,4}^{s0}) \subset C_1 \cup \{P_2\}$. 
Since $X_{4,4}^{s0} \cap \partial Y = C_1 \cup C_2 \cup L_0 \cup \{P_2\}$, 
we have $\partial \cP_{4,d}^{s0}
  = \cF(C_1) \cup \cF(C_2) \cup \cF((P_2P_{-3})) \cup \cF(P_2)$ by \Tqbd(1). 
But $\cF(C_2)$ and $\cF((P_2P_{-3}))$ are not face components, 
because dual varieties of $\Zar(C_2)$ and $\Zar((P_2P_{-3}))$ are 
linear subspaces of $cH_{4,4}^{s0}$ of codimension $2$. 
Thus all the face components of $\cP_{4,4}^{s0}$ are $\cF(C_1)$ and $\cF(P_2)$. 
Therefore $\partial \cP_{4,4}^{s0} = \cF(C_1) \cup \cF(P_2)$. 
\end{proof}

\bigskip

{\it Proof of \Tqae(1).} 
Put $\Omega := \{(-1 \colon -1 \colon 1 \colon 1)\} \cup 
 \big\{(t \colon 1 \colon 1 \colon 1) \in \P_{\R}^3$ 
  $\big|$ $t \in \R \big\}$. 
By \Tqbby, it is enough to show that $\Phi_{4,4}^{s0}(\Omega) 
  \supset C_1 \cup \{P_2\} \supset \cE(X_{4,4}^{s0})$. 
But this is clear. \QED

\bigskip

We regard $\cH_{4,3}^{s0} = \R^4$, by identifying 
$\displaystyle f = \sum_{i=0}^3 p_i s_i \in \cH_{4,3}^{s0}$ 
and $(p_0$, $p_1$, $p_2$, $p_3) \in \R^4$. 
We also use $(p_0$, $p_1$, $p_2$, $p_3)$ as a coordinate system of 
$\cH_{4,3}^{s0} = \R^4$. 
We denote the local cone of $\cP_{4,4}^{s0}$ at 
$(t \colon 1 \colon 1 \colon 1) \in \P_{\R}^3$ by $\cL_t^{s0}$. 
Note that if $f \in \cF(C_1)$, there exists $t \in \R$ 
such that $f(t$, $1$, $1$, $1) = 0$. 
Thus $f \in \cL_t^{s0}$. 
For $t=\infty \in \P_{\R}^1$, 
we denote the local cone of $\cP_{4,4}^{s0}$ at 
$(1 \colon 0 \colon 0 \colon 0) \in \P_{\R}^3$ by $\cL_{\infty}^{s0}$. 

\bigskip

We shall observe $\frg_t$, $\frg_{\infty}$ and 
$\frp \in \cP_{4,4}^{s0}$. Note that 
\begin{align*}
 3\frg_t(a,b,c,d) 
 & = 3 s_0 - 2(t+1) (s_1-s_3) + (t^2+2t-1) s_2 \\
 & = \big(a^2+b^2-c^2-d^2 + (t+1)(c d-a b)\big)^2 \\
 & \hskip20pt + \big(a^2-b^2+c^2-d^2 + (t+1)(b d-a c)\big)^2 \\
 & \hskip20pt + \big(a^2-b^2-c^2+d^2 + (t+1)(b c-a d)\big)^2, \\
 \frg_{\infty}(a,b,c,d) 
 & = s_2 = (a b-c d)^2 + (a c-b d)^2 + (a d-b c)^2, \\
 \frp 
 & = s_2 - s_3 
   = (a-b)^2(c-d)^2 + (a-c)^2(b-d)^2 + (a-d)^2(b-c)^2. 
\end{align*}
Espacially, $\frg_t$, $\frg_{\infty}$, $\frp \in \Sigma_{4,4}$. 

If $f \in \cE(\cP_{4,4}) \cap \Sigma_{4,4}$, 
then there exists $g \in \cP_{2,4}$ such that $f = g^2$. 
Therefore $\frg_t$, $\frg_{\infty}$, $\frp \notin \cE(\cP_{4,4})$. 
But $\frg_t$, $\frg_{\infty}$, $\frp \notin \cE(\cP_{4,4}^{s0})$ 
as the following Lemma. 

For $f(a,b,c,d) \in \R[a,b,c,d]$, 
we dnote $\displaystyle \frac{\partial}{\partial a}f$ by $f_a$, 
$\displaystyle \frac{\partial^2}{\partial a^2}f$ by $f_{aa}$, and so on. 

\def\Tqcceb{Lemma 3.6}
\proclaim{Lemma 3.6} 
{\sl $\frg_t \in \cE(\cP_{4,4}^{s0})$ for 
all $t \in \P_{\R}^1$, and $\frp \in \cE(\cP_{4,4}^{s0})$. 
These are characterized as follows: \par}
{\parindent=20pt
\Item{\rm (1)} {\sl Let $t \in \R - \{1$, $-3\}$. 
If $f \in \cP_{4,4}^{s0}$ satisfies 
$f(t$, $1$, $1$, $1) = 0$ and $f(-1$, $-1$, $1$, $1) = 0$, 
then there exists $\alpha \geq 0$ such that $f = \alpha \frg_t$.} 
\Item{\rm (2)} {\sl If $f \in \cP_{4,4}^{s0}$ satisfies 
$f(x$, $x$, $1$, $1) = 0$ for all $x \in \R$, 
then there exists $\alpha \geq 0$ such that $f = \alpha \frg_1$.} 
\Item{\rm (3)} {\sl If $f \in \cP_{4,4}^{s0}$ satisfies 
$f(x$, $y$, $z$, $-x-y-z) = 0$ for all $x$, $y$, $z \in \R$, 
then there exists $\alpha \geq 0$ such that $f = \alpha \frg_{-3}$.} 
\Item{\rm (4)} {\sl If $f \in \cP_{4,4}^{s0}$ satisfies 
$f(0$, $0$, $0$, $1) = 0$ and $f(-1$, $-1$, $1$, $1) = 0$, 
then there exists $\alpha \geq 0$ such 
that $f = \alpha \frg_{\infty}$.} 
\Item{\rm (5)} {\sl If $f \in \cP_{4,4}^{s0}$ satisfies 
$f(0$, $0$, $0$, $1) = 0$ and $f(x$, $1$, $1$, $1) = 0$ for all $x \in \R$, 
then there exists $\alpha \in \R_+$ such that $f = \alpha \frp$.} 

}
\endproclaim

\Proof 
Note that of $f \in \cP_{4,4}^{s0}$ satisfies 
$f(a,b,c,d) = 0$, then $f_a(a,b,c,d) = 0$. 
Similarly, if $f_{aa}(a,b,c,d) = 0$, then $f_{aaa}(a,b,c,d) = 0$. 
Otherwise, $f$ will be negative at a certain point near $(a$, $b$, $c$, $d)$. 
$f \in \cH_{4,4}^s$ can be written as 
$f = p_0s_0 + p_1s_1 + p_2s_2 + p_3s_3$ by 
$p_0$, $p_1$, $p_2$, $p_3 \in \R$. 

\smallskip

(1) Take $t \in \R - \{1$, $-3\}$. Let's consider a system of equations 
\[f(t,1,1,1) = 0, \quad f_a(t,1,1,1) = 0, \quad f(-1,-1,1,1) = 0. \eqno (*)\]
Let $a_{0,j} := s_j(t,1,1,1)$, $a_{1,j} := (s_j)_a(t,1,1,1)$, 
$a_{2,j} := s_j(1,1,-1,-1)$, and $A := (a_{i,j}) \in M_{3,4}(\R)$. 
Then, ($*$) is equivalent to $A {\bf p} = {\bf 0}$. 
That is 
\[\left(\begin{matrix}
  (t-1)^2(t^2+2t+3) & 3(t-1)^2(t+2) & 3(t-1)^2 & 3(t-1)^2 \\
  4(t^3-1) & 9(t^2-1) & 6(t-1) & 6(t-1) \\
  0 & -16 & 0 & -16 \end{matrix}\right)\left(
\begin{matrix} p_0 \\ p_1 \\ p_2 \\ p_3 \end{matrix}\right) 
= \left(\begin{matrix} 0 \\ 0 \\ 0 \end{matrix}\right).\]
Using Mathmatica, we can soon check that $\Ker A = \R \cdot \frg_t$. 
If $f \in \cP_{4,4}^{s0}$ satisfies $f(t,1,1,1) = 0$, 
then $f_a(t,1,1,1) = 0$ always holds. 
Thus, if $f \in \cP_{4,4}^{s0}$ satisfies $f(t,1,1,1) = 0$ 
and $f(-1,-1,1,1) = 0$, 
then $f = \alpha \frg_t$ by a certain $\alpha > 0$. 

\smallskip

(2) Consider a system of equations $f(0,0,1,1)=0$, $f(2,2,1,1) = 0$ 
instead of ($*$). Then $\dim \Ker A = 2$, and 
$\frg_1$ and $g := s_1-2s_2-s_3$ is a basis of $\Ker A$. 
Note that $g$ is not PSD, and $\frg_1 = \frf_1^{ab}$. 
Since $\frg_1(x,1,1,1) + c g(x,1,1,1) 
 = (x-1)^3(x-1+3c)$, $\frg_1 + c g$ is PSD only if $c=0$. 

\smallskip

(3) Consider $f(1,2,3,-6) = 0$, $f_a(1,2,3,-6)=0$, $f(1,2,4,-7) = 0$. 

(4) Consider $f(0,0,0,1) = 0$, $f_a(0,0,0,1)=0$, $f(-1,-1,1,1) = 0$. 

(5) Consider $f(2,1,1,1) = 0$, $f(0,0,0,1) = 0$, $f_a(0,0,0,1)=0$. 

Each $A$ of the cases (2)---(5) are as follows: 
\begin{align*}
& (2) \hskip10pt A = \left(\begin{matrix}
  2 & 2 & 1& 0 \\ 18 & 26 & 9 & 8 \end{matrix}\right), 
\hbox{} \hskip20pt (3) \hskip10pt A = \left(\begin{matrix}
  1538 & -962 & 769 & 576 \\ 148 & 248 & 314 & 516 \\
      2898 & -2002 & 1449 \end{matrix}
\right), \\
& (4) \hskip10pt A = \left(\begin{matrix}
  1 & 0 & 0 & 0 \\ 0 & 1 & 0 & 0 \\ 0 & 16 & 0 &-16 
\end{matrix}\right),
\hbox{} \hskip20pt (5) \hskip10pt A = \left(\begin{matrix}
  11 & 12 & 3 &3 \\ 1 & 0 & 0 & 0 \\ 0 & 1 & 0 & 0 
\end{matrix}\right). 
\end{align*}
\end{proof}

\bigskip

$\frg_t$ ($t \in \P_{\R}^1$) degenerates when $t=1$, $-3$. Note that 
\[\frg_{-3} = S_1^2 (3S_3 - 2T_{1,1}).\]
Thus $\cF(L_0) = \R_+ \cdot \frg_{-3}$. 

Since $\frg_1(x$, $x$, $1$, $1) = 0$ for all $x \in \P_{\R}^1$, 
we have $\cF(C_2) = \R_+ \cdot \frg_1$. 
These also implies that $\cF(L_0)$ and $\cF(C_2)$ are 
not a face component of 
$\cP_{4,4}^{s0}$, and we can omit 
$\big\{ (x \colon x \colon 1 \colon 1) \in \P_{\R}^3$ $\big|$ 
$x \in \R_+\big\}$ from the test set. 

\def\Tqch{Lemma 3.7}
\proclaim{Lemma 3.7} 
{\sl $\cL_t^{s0} 
   = \R_+ \cdot \frg_t + \R_+ \cdot \frp$, and 
the discriminant of $\cF(C_1)$ and $\cF(P_2)$ are }
\begin{align*}
 & \disc_{C_1}(p_0,p_1,p_2,p_3) 
     = 8 p_0^2 - 9p_1^2 + 12p_0p_1 + 12p_0p_2 + 12p_0p_3, \\
 & \disc_{P_2}(p_0,p_1,p_2,p_3) = p_1+p_3. 
\end{align*}
\endproclaim

\Proof 
Since $P_2 = (0 \colon 1 \colon 0 \colon 1)$, 
$\disc_{P_2}(p_0$, $p_1$, $p_2$, $p_3) = p_1+p_3$, 
by \cite[Remark 1.28]{RefAc}. 

Since $\frg_t$, $\frp \in \cL_t^{s0}$ ($t \in \P_{\R}^1$), 
we have $\dim \cL_t^{s0} \geq 2$. 
On the other hand, since $\dim \cL_t^{s0} 
 < \dim \cP_{4,4}^{s0} = 3$, we have 
$\dim \cL_t^{s0} = \dim \cL_{\infty}^{s0} = 2$ ($t \ne 1$). 
Since $\frg_t$, $\frp \in \cE(\cP_{4,4}^{s0})$, 
we have $\cL_t^{s0} = \R_+ \cdot \frg_t + \R_+ \cdot \frp$ 
for all $t \in \P_{\R}^1$. 

Using PC, we can check that $\frg_t$ ($\forall t \in \P_{\R}^1$) 
and $\frp$ exists on the hypersurface in $\cH_{4,4}^{s0}$ 
defined by $8 p_0^2 - 9p_1^2 + 12p_0p_1 + 12p_0p_2 + 12p_0p_3$. 
This equation is also the defining equation of the dual variety of $C_1$. 
So, this is $\disc_{C_1}$. 
\end{proof}

\bigskip

{\it Proof of \Tqcb.} 
By the above lemma, we have 
\begin{align*}
 & \cF(P_2) 
   = \left\{\left. \sum_{i=0}^3 p_i s_i \in \cH_{4,4}^{s0} \; \right| \; 
    \hbox{$p_1 + p_3 = 0$, $p_0 \geq 0$, 
    $-9p_1^2 + 12p_0 p_2 + 8p_0^2 \geq 0$} \right\}, \\
 & \cF(C_1) = \left\{ 
  \sum_{i=0}^3 p_i s_i \in \cH_{4,4}^{s0} \; \left| \; 
  \vcenter{\hbox{$p_1+p_3 \geq 0$, $p_0 \geq 0$,}
           \hbox{$-9p_1^2 + 12p_0p_1 + 12p_0 p_2 + 12p_0 p_3 + 8p_0^2 = 0$}}\ 
   \right.\right\}. 
\end{align*}
Thus, all the extremal elements of $\cP_{4,4}^{s0}$ are 
$\frg_t$ ($t \in \P_{\R}^1$) and $\frp$. 

Thus, for $f = s_0 + p s_1 + q s_2 + r s_3 \in \cH_{4,4}^{s0}$, 
$f(a) \geq 0$ for all $a \in \P_{\R}^3$ if 
and only if $p + r \geq 0$ and $- 9p^2 + 12p + 12q + 12r + 8 \geq 0$. 

(4) follows from 
$\partial \cP_{4,4}^{s0} = \cF(C_1) \cup \cF(P_2)$. \QED

\bigskip

{\it Proof of \Tqaa (1), \Tqabn \ and \Tqadn (1).} 
Let $t_0 := \sigma_1^4 - 256 \sigma_4$, 
$t_1 := \sigma_1^2 \sigma_2 - 96 \sigma_4$, 
$t_2 := \sigma_2^2 - 36 \sigma_4$, 
$t_3 := \sigma_1 \sigma_3 - 16 \sigma_4$. 
Then $s_0 = t_0 - 4 t_1 + 2 t_2 + 4 t_3$, 
$s_1 = t_1 - 2 t_2 - t_3$, 
$s_2 = t_2 - 2 t_3$ 
and $s_3 = t_3$. 
Using these substitution for 
$\frg_t$, $\frg_{\infty}$ and $\frp$, 
we obtain \Tqab. 

Take $f = p_0 s_0 + p_1 s_1 + p_2 s_2 + p_3 s_3 
 = q_0 t_0 + q_1 t_1 + q_2 t_2 + q_3 t_3 \in \cH_{4,4}^{s0}$. 
Since $t_0 = s_0 + 4 s_1 + 6 s_2 + 12 s_3$, 
$t_1 = s_1 + 2 s_2 + 5 s_3$, 
$t_2 = s_2 + 2 s_3$ and $t_3 = s_3$, 
we have $p_0 = q_0$, $p_1 = 4 q_0 + q_1$, $p_2 = 6 q_0 + 2 q_1 + q_2$, 
and $p_3 = 12 q_0 + 5 q_1 + 2 q_2 + q_3$. 
Substitute these for $p_i$ in $\disc_{C_1}$ and $\disc_{P_2}$ of \Tqch, 
we obtain $d_1$ and $d_2$ of \Tqad (1). 
\Tqaa (1) follows from these.  \QED

\bigskip

{\it Proof of \Tqadeg(1).} 
Let $f(x,y) := \frg_t(x,y,1,-x-y-1)/(t+3)^3$ for 
$t \in \P_{\R}^1 - \{-3\}$. 
If $\frg_t$ is reducible, then $f$ is also reducible. By 
\[\frac{\partial}{\partial x} f(x,y)
  = 2(2x+y+1)(x^2+x y+y^2+x+3y+1)\]
and so on, we have 
\[\Sing(V_{\C}(f)) = \big\{(-1 \colon -1 \colon 1), \, 
  (-1 \colon 0 \colon 1), \, (0 \colon 1 \colon 1)\big\}.\]
Moreover, these are acnodes. 
Assume that $f = gh$. 
If $\deg g = 1$, then $\# \Sing(V_{\C}(f)) \allowbreak = 4$ or 
$\# \Sing(V_{\C}(f)) \subset V_{\C}(g)$. 
This cannot occur. 
Thus, $g$ and $h$ are irreducible quadric curves 
which intersect transversally. 
Then, $\# \Sing(V_{\C}(f)) = 4$. 
Therefore, $V_{\C}(f)$ must be an irreducible rational quartic curve. \QED

\bigskip

{\it Proof of \Tqabc.} 
$\cE(\cP_{4,4}^{s0}) \subset \Sigma_{4,4}$ is already proved. 
Since, any element of $\cP_{4,4}^{s0}$ can be written as a sum 
of some elements in $\cE(\cP_{4,4}^{s0})$, 
we have $\cP_{4,4}^{s0} \subset \Sigma_{4,4}$. 

Assume that $\exists f \in \cE(\cP_{4,4}^{s0}) 
  \cap \cE(\cP_{4,4}) \ne \emptyset$. 
$f$ is SOS, since $\cE(\cP_{4,4}^{s0}) \subset \Sigma_{4,4}$. 
Since, $f \in \cE(\cP_{4,4})$, 
we have $f \in \cE(\Sigma_{4,4})$. 
Thus, there exists $g \in \cH_{4,2}$ such that $f = g^2$. 
Then $V_{\R}(g) = V_{\R}(f)$. 
Since $\# V_{\R}(\frg_t) \geq 2$ and 
$\# V_{\R}(\frp) \geq 2$, we have $\# V_{\R}(g) \geq 2$. 
Such conic $g$ satisfies $\dim_{\R} V_{\R}(g) \geq 1$. 
But, $V_{\R}(f)$ is a finite set. \QED

\bigbreak
{\bf 3.3. The PSD cone $\cP_{4,4}^{s0+}$}
\par\penalty1000\vskip0.4em plus0.1em minus0.1em
In this subsection, we shall study 
$\cP_{4,4}^{s0+} := \cP(\P_+^3$, $\cH_{4,4}^{s0})$. 
The aim of this subsection is to prove the following theorem. 

\def\Tqci{Theorem 3.8}
\def\Tqcin{3.8}%
\proclaim{Theorem 3.8} 
{\rm (I)} {\sl For a monic 
\[f = s_0 + p s_1 + q s_2 + r s_3 \in \breve{\cH}_{4,4}^{s0},\]
$f({\bf a}) \geq 0$ for all ${\bf a} \in \R_+^4$ if and only if the following 
``{\rm (1)} or {\rm (2)}'' and ``{\rm (3)} or {\rm (4)}'' hold: \par}
{\parindent=20pt
\Item{(1)} {\sl $p \leq -4$ and $p^2 \leq 4 q - 8$.}
\Item{(2)} {\sl $p \geq -4$ and $2p+q+2 \geq 0$.}
\Item{(3)} {\sl $p \leq -2/3$ and $9 p^2 \leq 12 p + 12 q + 12 r + 8$.}
\Item{(4)} {\sl $p \geq -2/3$ and $3q+3r \geq 1$.}

}
{\rm (II)} {\sl All the extremal elements of $\cP_{4,4}^{s0+}$ are 
positive multiples of 
$\frf_t^{ab}$ {\rm ($0 \leq t \leq 5$)}, 
$\frf_t^c$ {\rm ($5 < t < \infty$)}, 
$\frp = s_2-s_3$, $\frq_1 = s_1-2s_2$ or $\frq_2 = s_3$.}

{\rm (III)} {\sl The following set is a test set for 
$(\P_+^3$, $\cH_{4,4}^{s0+})$. }
\[\big\{(t \colon 1 \colon 1 \colon 1) \in \P_+^3 \; \big| \; 
      \hbox{$t \geq 0$}\big\}
   \cup \big\{(0 \colon 0 \colon t \colon 1) \in \P_+^3 \; \big| \; 
   \hbox{$t \geq 0$}\big\}.\]
\endproclaim
\niceskip

This theorem will be proved after \Tqcn. 

Essentially, we use the same symbols as the previous subsection, 
but there are some changes. 
Let $A := \P_+^3:(a_0 \colon a_1 \colon a_2 \colon a_3)$, 
$X_{4,4}^{s0+} := \Phi_{4,4}^{s0}(\P_+^3) 
 = X(\P_+^3$, $\cH_{4,4}^{s0}) 
 \subset \P\big((\cH_{4,4}^{s0})^{\vee}\big)$. 
As \S 3.2, put $D_0 := \Psi(\tilde{D}_0)$, 
$D_1^+ := \Psi(\tilde{D}_1^+) \subset D_1$, 
$P_1 := (2 \colon 3 \colon 1 \colon 1)$, 
$C_1^+ := \Psi(\tilde{C}_1^+) \cup \{P_1\} \subset C_1$, 
$C_2^+ := \Psi(\tilde{C}_2^+) \subset C_2$, 
$C_i := \Psi(\tilde{C}_i)$ for $i=3$, $4$ and 
$P_j := \Psi(\tilde{P}_j)$ for $j=3$, $4$, $5$. 
Note that 
\begin{align*}
 & P_3 = (1 \colon 2 \colon 1 \colon 1) = \Phi_{4,4}^{s0}(0,1,1,1), \\
 & P_4 = (2 \colon 2 \colon 1 \colon 0) = \Phi_{4,4}^{s0}(0,0,1,1), \\
 & P_5 = (1 \colon 0 \colon 0 \colon 0) = \Phi_{4,4}^{s0}(0,0,0,1). 
\end{align*}
%

\proclaim{Lemma 3.9} 
{\sl Let $Z := \P_+^3/\frS_4 - \tilde{C}_1^+ - \tilde{C}_2^+ 
   - \{\tilde{P}_1$, $\tilde{P}_3$, $\tilde{P}_4$, $\tilde{P}_5\}$. }
{\parindent=20pt
\Item{{\rm (1)}} {\sl $\Psi \colon \P_+^3/\frS_4 \lto X_{4,4}^{s0+}$ is 
continuous bijective map and 
$\Psi : Z \lto \Psi(Z)$ is an isomprohism.}
\Item{{\rm (2)}} 
$\Delta^0(X_{4,4}^{s0+}) = \{P_3$, $P_4$, $P_5\}$, 
$\Delta^1(X_{4,4}^{s0+}) = \{C_1^+$, $C_2^+$, $C_3$, $C_4\}$, 
$\Delta^2(X_{4,4}^{s0+}) = \{D_0$, $D_1^+\}$, 
$\Delta^3(X_{4,4}^{s0+}) = \{\Int(X_{4,4}^{s0+})\}$. 

}
\endproclaim

\Proof 
(1) We use the same symbols with the proof of \Tqccea. 
Note that $E_0 \cap A_s^+ = \emptyset$. 
So, it is enough to show that $\Phi_{4,4}^{s0}$ is injective 
on $\displaystyle A_s^+ \cap \bigcup_{\tau \in \frS_4} \tau(D_1')$. 
It is enough to show that $\Phi_{4,4}^{s0}$ is injective 
on $\P_+^2 \cap D_1'$. 
It's Jacobian is equal to 
$\displaystyle J(x_1,x_2) := \det \left( \frac{\partial h_i}{\partial x_j} 
   \right)_{1 \leq i,j \leq 2}$, 
where $h_i(x_1,x_2) := s_i(x_1,x_2,1,1)/s_0(x_1,x_2,1,1)$ ($i=1$, $2$). 
Using PC, we have 
\begin{align*}
 & J(x,y) = \frac{4(x-1)(y-1)(x-y)w(x,y)}{s_0(x,y,1,1)^3}, \\
 & w(x,y) := \frac{1}{8}\Big((x+y-2)^4(x+y+2)^2 
              + (x-y)^4(3(x+y)^2+28(x+y)+12) \\
 & \hbox{} \hskip80pt + 4(x-y)^2(x+y-2)^2((x+y)^2+6(x+y)+4)\big). 
\end{align*}
Thus $J(x$, $y) \geq 0$ on $\P_+^2 \cap D_1'$, 
and $J(x$, $y) = 0$ only at points of the form 
$(x \colon 1 \colon 1 \colon 1)$ or $(1 \colon x \colon 1 \colon 1)$ 
or $(x \colon x \colon 1 \colon 1)$. 
Thus $\Psi : Z \lto \Psi(Z)$ is an isomprohism. 

\smallskip

(2) follows from \Tqca (3) and the proof of \Tqccea. 
\end{proof}

\proclaim{Lemma 3.10} 
(1) $\partial \cP_{4,4}^{s0+} = \cF(C_1^+) \cup \cF(C_3) \cup \cF(C_4) 
   \cup \cF(P_3) \cup \cF(P_4) \cup \cF(P_5)$. 

(2) {\sl Take $f \in \cH_{4,4}^{s0}$. 
If $f(x,1,1,1) \geq 0$, $f(0,x,1,1) \geq 0$, $f(0,0,x,1) \geq 0$ 
for all $x \geq 0$, then $f \in \cH_{4,4}^{s0}$. }
\endproclaim

\Proof 
(1) $\Int(X_{4,4}^{s0+})$, $\cF(D_0)$ and $\cF(D_1^+)$ are not 
face components of $\cP_{4,4}^{s0+}$ by \Tqbaf. 
$\cF(C_2^+)$ is not also a face component of $\cP_{4,4}^{s0+}$, 
because $C_2^+$ is an open line segment $(P_1$, $P_4)$. 
Thus, we have (1). 

\smallskip

(2) Let 
\begin{align*}
 & A_1^+ := \big\{(t \colon 1 \colon 1 \colon 1) \in \P_+^3 \; \big| \; 
   \hbox{$t \geq 0$}\big\}, \\
 & A_2^+ := \big\{(t \colon t \colon 1 \colon 1) \in \P_+^3 \; \big| \; 
   \hbox{$0 \leq t \leq 1$}\big\}, \\
 & A_3 := \big\{(0 \colon t \colon 1 \colon 1) \in \P_+^3 \; \big| \; 
   \hbox{$t \geq 0$}\big\}, \\
 & A_4 := \big\{(0 \colon 0 \colon t \colon 1) \in \P_+^3 \; \big| \; 
   \hbox{$t \geq 0$}\big\}. 
\end{align*}
Note that $\Phi_{4,4}^{s0}(A_i^+) \supset C_i^+$ ($i=1$, $2$), 
and $\Phi_{4,4}^{s0}(A_j) \supset C_j$ ($j=3$, $4$). 
By \cite[Corollary 1.3]{RefRie} \ or \cite[Corollary 2.1]{RefTa}, 
we can choose $A_1^+ \cup A_2^+ \cup A_3 \cup A_4$ as a test set for 
$(\P_+^3$, $\cH_{4,4}^{s0})$. 
Since $\cF(C_2^+)$ is not a face component of $\cP_{4,4}^{s0+}$ 
and $P_1 \in C_1^+$, $P_4 \in \Cls(C_3) \cap \Cls(C_4)$, 
we can omit $A_2^+$ from the test set. 
Thus, if $f \in \cH_{4,4}^{s0}$ satisfies 
$f(x$, $1$, $1$, $1) \geq 0$, $f(0$, $x$, $1$, $1) \geq 0$ and 
$f(0$, $0$, $x$, $1) \geq 0$ for all $x \geq 0$, 
then $f \in \cP_{4,4}^{s0+}$. 
\end{proof}

In fact, $\cF(C_3)$ is not a face component, and we can omit $A_3$ from 
the test set. But it will be proved later. 
We summarize here what $C_1^+$, $C_3$ and $C_4$ are. 

\proclaim{Lemma 3.11} 
{\parindent=20pt
\Item{(1)} {\sl $\Zar(C_1^+)$ is a conic defined by 
$x_1^2 - 2x_1 x_2 - 3x_0 x_2 + 3x_2^2 = 0$, $x_2-x_3 = 0$. 
Especially, $\Zar(C_1^+)$ is nonsingular. 
The ends of $C_1^+$ are $P_3$ and $P_5$.}
\Item{(2)} {\sl $\Zar(C_3)$ has a cusp at $P_3$. 
The ends of $C_3$ are $P_4$ and $P_5$.} 
\Item{(3)} {\sl $\Zar(C_4)$ is a conic defined by 
$x_1^2 - 2 x_2^2 - x_0 x_2 = 0$ on the plane $V_{\R}(x_3)$. 
The ends of $C_4$ are $P_4$ and $P_5$.} 

}
\endproclaim

Next, we shall study $\frf_t^{ab}$ ($0 \leq t \leq 5$), 
$\frf_t^c$ ($5 \leq t < \infty$), $\frp = s_2-s_3$, 
$\frq_1 = s_1-2s_2$, and $\frq_2 = s_3$. 
Note that 
\begin{align*}
 & \frf_t^{ab} =  \frac{1}{3}
    \big(3 s_0 - 2(t+1) s_1 + 2(2t-1) s_2 + (t^2+3) s_3\big), \\
 & \frf_t^c = \frac{1}{9}
    \big(9 s_0 - 6(t+1) s_1 + (t^2+2t+19)s_2+ 2(t^2+5t-8) s_3\big), 
\end{align*}
and $\frf_5^{ab} = \frf_5^c$. 
Put $\frf_{\infty}^c := s_2 + 2 s_3$. 
Since $\frf_{\infty}^c = \frp + 3 \frq_2$, 
$\frf_{\infty}^c$ is not extremal. 
The author studied $\Phi_{4,4}^{c0}(t \colon 1 \colon 1 \colon 1) 
 \in \cF(C_1^+)$, dividing three cases (a) $0 \leq t < 1$, 
(b) $1<t \leq 5$ and (c) $t>5$. 
The symbol $\frf_t^{ab}$ stands for cases (a) and (b). 
For $u \geq 0$, let 
\[\frh_u^c := 3u^2 s_0 - 6u(u^2+1) s_1 
     + 3(u^4+4u^2+1) s_2 + 2(3u^4+3u^3+2u^2+3u+3) s_3.\]
If $t = (3u^2-u+3)/u$, then 
$\displaystyle \frh_u^c = 3u^2 \frf_t^c$. 
So, $\frh_u^c$ is not a new polynomial, 
but it is convenient to study $\cF(C_4)$ for the property 
$\frh_u^c(0$, $0$, $u$, $1) = 0$. 

We shall denote the local cone of $\cP_{4,4}^{c0+}$ at 
the point $(t \colon 1 \colon 1 \colon 1) \in \P_+^3$ by $\cL_t^{C_1}$, 
and the local cone at the point $(0 \colon 0 \colon t \colon 1)$ 
by $\cL_t^{C_4}$. 

\proclaim{Lemma 3.12} 
{\sl $\frf_t^{ab}$ {\rm ($0 \leq t \leq 5$)}, 
$\frf_t^c$ {\rm ($5 < t < \infty$)}, $\frp$, 
$\frq_1$, and $\frq_2$ are extremal elements 
of $\cP_{4,4}^{s0+}$. 
These are characterized as follows: \par}
{\parindent=20pt
\Item{\rm (1)} {\sl Let $0 \leq t < 1$ or $1 < t \leq 5$. 
If $f \in \cP_{4,4}^{s0+}$ satisfies 
$f(t$, $1$, $1$, $1) = 0$ and $f(0$, $0$, $1$, $1) = 0$, 
then there exists $\alpha \in \R_+$ such that $f = \alpha \frf_t^{ab}$.} 
\Item{\rm (2)} {\sl If $f \in \cP_{4,4}^{s0}$ satisfies 
$f_{aa}(1$, $1$, $1$, $1) = 0$ and $f(x$, $x$, $1$, $1) = 0$ for all 
$x \geq 0$, 
then there exists $\alpha \in \R_+$ such that $f = \alpha \frf_1^{ab}$.} 
\Item{\rm (3)} {\sl Assume that $t$, $u \in \R_+$ satisfy $3u^2-(t+1)u+3=0$. 
If $f \in \cP_{4,4}^{s0}$ satisfies 
$f(t$, $1$, $1$, $1) = 0$ and $f(0$, $0$, $u$, $1) = 0$, 
then there exists $\alpha \in \R_+$ such that $f = \alpha \frf_t^c$.} 
\Item{\rm (4)} {\sl If $f \in \cP_{4,4+}^{s0}$ satisfies 
$f(0$, $0$, $0$, $1) = 0$, $f_a(0$, $0$, $0$, $1) = 0$ 
and $f(x$, $1$, $1$, $1) = 0$ for all $x \geq 0$, 
then there exists $\alpha \in \R_+$ such that $f = \alpha \frp$.} 
\Item{\rm (5)} {\sl If $f \in \cP_{4,4}^{s0}$ satisfies 
$f(0$, $1$, $1$, $1) = 0$, $f(0$, $0$, $1$, $1)=0$ 
and $f(0$, $0$, $0$, $1) = 0$, then there exists $\alpha \in \R_+$ such 
that $f = \alpha \frq_1$.} 
\Item{\rm (6)} {\sl If $f \in \cP_{4,4}^{s0}$ satisfies 
$f(0$, $0$, $x$, $y) = 0$ for all $x$, $y \in \R_+$, 
then there exists $\alpha \in \R_+$ such that $f = \alpha \frq_2$.} 

}
\endproclaim

\Proof 
We shall show that $\frf_t^{ab}$ ($0 \leq 1 \leq 5$), 
$\frf_t^c$ ($t \geq 5$), $\frp$, $\frq_1$ 
and $\frq_2$ belong to $\cP_{4,4}^{s0+}$. 
Since 
\[\frf_t^{ab}(0,x,1,1) 
  = \frac{1}{3} x(x+2) \left(\left(t - \frac{2(x-1)^2}{(x+2)}\right)^2 
        + \frac{x(16-x)(x-1)^2}{(x+2)^2}\right),\]
we have $\frf_t^{ab}(0$, $x$, $1$, $1) \geq 0$ if $x \leq 16$. 
On the other hand 
\begin{align*}
 & \frf_t^{ab}(0,x,1,1) 
    = \frac{1}{3} x\Big(18(25-t)^2 + \big(t^2+120(5-t) +1575\big)(x-16) \\
 & \hskip80pt + \big(4(5-t)+120\big)(x-16)^2 + 3(x-16)^3\Big), 
\end{align*}
we have $\frf_t^{ab}(0$, $x$, $1$, $1) \geq 0$ for $x \geq 16$. 
Similarly, 
\begin{align*}
 \frf_t^{ab}(x,1,1,1) & = (x-t)^2 (x-1)^2 \geq 0, \\
 \frf_t^{ab}(0,0,x,1) 
    & = \frac{1}{3} (x-1)^2 \left(3\left(x-\frac{t-2}{3}\right)^2
        + \frac{1}{3}(5-t)(1+t)\right) \geq 0, \\
 \frf_t^c(x,1,1,1) & = (x-t)^2 (x-1)^2 \geq 0, \\
 \frf_t^c(0,x,1,1) & = \frac{(2x+1)^2}{9}
     \left(t - \frac{(x-1)^2(6x+5)}{(2x+1)^2}\right)^2 
       + \frac{24x(x-1)^2(x+2)(3x+2)}{9(2x+1)^2} \\
 & \geq 0, \\
 \frf_t^c(0,0,x,1) & = \frac{1}{9}(3x^2 - (t+1)x + 3)^2 \geq 0, \\
 \frh_u^c(0,0,x,1) & = 3(x-u)^2(u x-1)^2 \geq 0, \\
 \frq_1(x,1,1,1) & = 3x(x-1)^2 \geq 0, \\
 \frq_1(0,x,1,1) & = 2x(x-1)^2 \geq 0, \\
 \frq_1(0,0,x,1) & = x(x-1)^2 \geq 0, \\
 \frq_2(x,1,1,1) & = 3(x-1)^2 \geq 0, \\
 \frq_2(0,x,1,1) & = x(x+2) \geq 0, \\
 \frq_2(0,0,x,1) & = 0. 
\end{align*}
Thus $\frf_t^{ab}$, $\frf_t^c$, $\frq_1$, 
$\frq_2 \in \cP_{4,4}^{s0+}$. 

The left part can be proved similarly as the proof of \Tqcceb. 

\smallskip

(1) Consider a system of equations $f(t,1,1,1) = 0$, 
$f_a(t,1,1,1) = 0$, $f(0,0,1,1) \allowbreak = 0$ instead 
of ($*$) in \Tqcceb. 
Then $A{\bf p} = {\bf 0}$ become 
\[\left(
\begin{matrix}
  (t-1)^2(t^2+2t+3) & 3(t-1)^2(t+2) & 3(t-1)^2 & 3(t-1)^2 \\
  4(t^3-1) & 9(t^2-1) & 6(t-1) & 6(t-1) \\
  2 & 2 & 1 & 0 \end{matrix}\right) 
\left(\begin{matrix}  p_0 \\ p_1 \\ p_2 \\ p_3 \end{matrix}\right) 
= \left(\begin{matrix}  0 \\ 0 \\ 0 \end{matrix}\right).\]
Using Mathematica, we can check $\Ker A = \R \cdot \frf_t^{ab}$ 
if $t \ne 1$. 
$f_a(t,1,1,1) = 0$ follows from $f(t,1,1,1) = 0$ if 
$f \in \cP_{4,4}^{s0+}$. 

\smallskip

(2) Consider $f_{aaa}(1,1,1,1)=0$, $f(0,0,1,1)=0$, $f_a(0,0,1,1) = 0$. 
In this case, 
\[A = \left(\begin{matrix}
  24 & 18 & 0 & 0 \\ 2 & 2 & 1 & 0 \\ 0 & 2 & 0 & 2 
\end{matrix}\right).\]

\smallskip

(3) This case is slightly complicated. 
Let $t = (3u^2-u+3)/u$ and consider the system of 
equations $f\big((3u^2-u+3)/u ,1,1,1\big) = 0$, 
$f_a\big((3u^2-u+3)/u ,1,1,1) = 0$, $f(0,0,u,1) = 0$. 
Then $A{\bf p} = {\bf 0}$ become 
\[\left(\begin{matrix}
  (t-1)^2(t^2+2t+3) & 3(t-1)^2(t+2) & 3(t-1)^2 & 3(t-1)^2 \\
  4(t^3-1) & 9(t^2-1) & 6(t-1) & 6(t-1) \\
  u^4+1 & u(u^2+1) & u^2 & 0 \end{matrix}\right) 
\left(\begin{matrix} p_0 \\ p_1 \\ p_2 \\ p_3 \end{matrix}\right) 
= \left(\begin{matrix} 0 \\ 0 \\ 0 \end{matrix}\right).\]
Using Mathematica, we can check $\Ker A = \R \cdot \frf_t^c$. 

\smallskip

(4) Same with (5) of \Tqcceb. 

\smallskip

(5) Consider $f(0,1,1,1)=0$, $f(0,0,1,1)=0$, $f(0,0,0,1) = 0$. 
In this case, 
\[A = \left(\begin{matrix}
  3 & 6 & 3 & 3 \\ 2 & 2 & 1 &0 \\ 1 & 0 & 0 & 0 
\end{matrix}\right).\]

\smallskip

(6) Consider $f(0,0,0,1)=0$, $f(0,0,1,1)=0$, $f(0,0,1,2)=0$. 
In this case, 
\[A = \left(\begin{matrix}
  1 & 0 & 0 & 0 \\ 2 & 2 & 1 & 0 \\ 17 & 10 & 4 & 0 \end{matrix}\right).\]
\end{proof}

\proclaim{Lemma 3.13} \par
{\parindent=20pt
\Item{(1)} {\sl $\frf_t^{ab} \in \cF(C_1) \cap \cF(P_4)$ and 
$\cL_t^{C_1} = \R_+ \cdot \frf_t^{ab} + \R_+ \cdot \frp$ 
for $0 < t < 1$ or $1 < t \leq 5$. }
\Item{(2)} {\sl $\frf_t^c \in \cF(C_1) \cap \cF(C_4)$ and 
$\cL_t^{C_1} = \R_+ \cdot \frf_t^c + \R_+ \cdot \frp$ 
for $t>5$. Moreover, 
$\cL_u^{C_4} = \R_+ \cdot \frh_u^c + \R_+ \cdot \frq_2$ 
for $u \geq 0$ with $t = (3u^2-u+3)/u$. }
\Item{(3)} $\frf_0^{ab} \in \cF(C_1) \cap \cF(P_3)
 \cap \cF(P_4)$. 
\Item{(4)} $\frf_5 := \frf_5^{ab} = \frf_5^c 
   \in \cF(C_1) \cap \cF(C_4) \cap \cF(P_4)$. 
\Item{(5)} $\frf_{\infty}^c 
   \in \cF(C_1) \cap \cF(C_4) \cap \cF(P_5)$. 
\Item{(6)} $\frp \in \cF(C_1) 
   \cap \cF(P_3) \cap \cF(P_5)$. 
\Item{(7)} $\frq_1 
     \in \cF(P_3) \cap \cF(P_4) \cap \cF(P_5)$. 
\Item{(8)} $\frq_2 
     \in \cF(C_4) \cap \cF(P_4) \cap \cF(P_5)$. 

}
\endproclaim

\Proof 
If $\cF(D)$ ($D \in \Delta(X_{4,4}^{c0+})$) is a face 
component of $\cP_{4,4}^{c0+}$, 
then $\dim \cF(D) = \dim (\partial \cP_{4,4}^{c0+}) 
  = \dim \cP_{4,4}^{c0+} - 1 = 3$. 
So, if $D_1$, $D_2$, $D_3$ are distinct elements of $\Delta(X_{4,4}^{c0+})$, 
and $\cF(D_i)$ ($i=1$, $2$, $3$) are face components, 
then $\dim \big(\cF(D_1) \cap \cF(D_2)\big) = 2$ and 
$\dim \big(\cF(D_1) \cap \cF(D_2) \cap \cF(D_3)\big) = 1$. 

Now, we shall prove (1)---(8). 

\smallskip

(1) Assume that $0 \leq t < 1$ or $1 < t \leq 5$. 
By previous lemma, 
we have $\frf_t^{ab} \in \cL_t^{C_1} \cap \cF(P_4)$ 
for $0 \leq t \leq 5$. 
Since $\dim \cL_t^{C_1} = 2$, 
we have $\cL_t^{C_1} 
 = \R_+ \cdot \frf_t^{ab} + \R_+ \cdot \frp \subset \cF(C_1)$, 

\smallskip

(2) Let $u>0$ and $t = (3u^2-u+3)/u \geq 5$. 
By previous lemma, $\frf_t^c \in \cF(C_1) \cap \cF(C_4)$. 
Since $\dim \cL_t^{C_4} = 2$, 
$\cL_u^{C_4} = \R_+ \cdot \frh_u^c + \R_+ \cdot \frq_2$. 
As (1), we have 
$\cL_u^{C_1} = \R_+ \cdot \frf_t^c + \R_+ \cdot \frp$. 

\smallskip

(3)---(8) can be proved similarly. 
\end{proof}

\bigskip

Note that $\frf_1^{ab} \in \cF(C_2^+)$, 
because $\frf_1^{ab}(x$, $x$, $1$, $1) = 0$ for all $x \in \R$. 
By \Tqcceb (2), we have $\cF(C_2^+) = \R_+ \cdot \frf_1^{ab}$. 
This also implies that $\cF(C_2^+)$ is not a face component. 

Using the above lemma, we shall determine the structure of the face 
components $\cF(C_1^+)$, $\cF(C_4)$, 
$\cF(P_3)$, $\cF(P_4)$ and $\cF(P_5)$. 

\def\Fan{\mathop{\rm Fan}\nolimits}
\proclaim{Lemma 3.14} 
{\sl For $f$, $g \in \cH_{4,4}^{s0}$, 
let $\Fan(f$, $g) :=  \R_+ \cdot f + \R_+ \cdot g$ be the fan 
whose edges are $f$ and $g$. Put 
\[W^{ab} := \R_+ \cdot \big\{\frf_t^{ab} \; \big| \; 
      \hbox{$0 \leq t \leq 5$} \big\} \subset \cH_{4,4}^{s0}, \quad
  W^c := \R_+ \cdot \big\{\frf_t^c \; \big| \; 
      \hbox{$t \geq 5$} \big\} \cup \R_+ \cdot \frf_{\infty}^c.\]
Then the followings hold. } \par
{\parindent=20pt
\Item{(1)} $\partial \cF(C_1^+) = W^{ab} \cup W^c 
   \cup \Fan(\frf_{\infty}^c$, $\frp) 
   \cup \Fan(\frp$, $\frf_0^{ab})$.  
\Item{(2)} $\partial \cF(C_4) = W^c 
   \cup \Fan(\frf_5$, $\frq_2) 
   \cup \Fan(\frq_2$, $\frf_{\infty}^c)$. 
\Item{(3)} $\partial \cF(P_3) = \Fan(\frf_0^{ab}$, $\frq_1) 
   \cup \Fan(\frq_1$, $\frp) 
   \cup \Fan(\frp$, $\frf_0^{ab})$. 
{\sl That is, $\cF(P_3)$ is a triangle cone with edges 
$\frf_0^{ab}$, $\frq_1$ and $\frp$. }
\Item{(4)} $\partial \cF(P_4) = W^{ab} 
   \cup \Fan(\frf_5$, $\frq_2) 
   \cup \Fan(\frq_2$, $\frq_1)
   \cup \Fan(\frq_1$, $\frf_0^{ab})$. 
\Item{(5)} {\sl $\cF(P_5)$ is a triangle cone with edges 
$\frp$, $\frq_1$ and $\frq_2$. 
Note that $\frf_{\infty}^c \in \Fan(\frp$, $\frq_2)$, and 
$\Fan(\frp$, $\frf_{\infty}^c) = \cF(P_5) \cap \cF(C_1)$, 
$\Fan(\frf_{\infty}^c$, $\frq_2) 
    = \cF(P_5) \cap \cF(C_4)$. }

}
\endproclaim

\begin{center}
\includegraphics[width=120mm,clip]{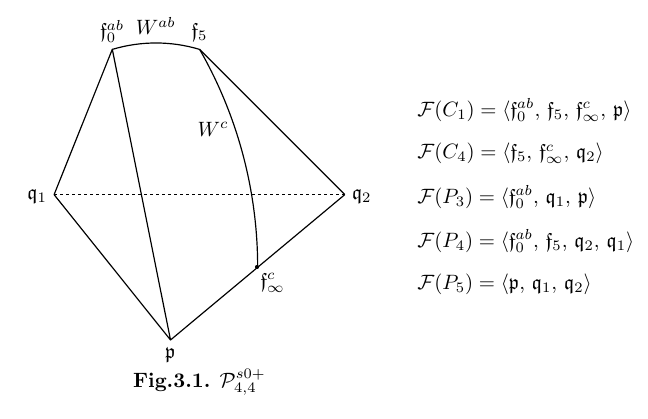}
\end{center}

By the above lemma, we know that $\partial \cP_{4,4}^{s0+}$ is 
enclosed by $\cF(C_1^+)$, $\cF(C_4)$, 
$\cF(P_3)$, $\cF(P_4)$ and $\cF(P_5)$. 
We don't need $\cF(C_3)$. See Fig.3.1. Thus, we have: 

\proclaim{Lemma 3.15} 
{\sl $\partial \cP_{4,4}^{s0+} = \cF(C_1^+) \cup \cF(C_4) 
   \cup \cF(P_3) \cup \cF(P_4) \cup \cF(P_5)$, and 
$\cE(X_{4,4}^{s0+}) \subset C_1^+ \cup C_4 
   \cup \{P_3$, $P_4$, $P_5\}$. 
Especially, $\cF(C_3)$ is not a face component of $\cP_{4,4}^{s0+}$.}
\endproclaim

{\it Proof of \Tqae(2).} 
Put $\Omega_+ := A_1^+ \cup A_4$. 
By \Tqbby, it is enough to show that $\Phi_{4,4}^{s0}(\Omega_+) 
  \supset C_1^+ \cup C_4 \cup \{P_3$, $P_4$, $P_5\}$. 
But this is clear. \QED

\bigskip

Geometrically, $C_3 - \{P_3$, $P_4$, $P_5 \}$ is included in the 
interior of the convex closure of $X_{4,4}^{s0+}$. 
So, any $f \in \cP_{4,4}^{s0+}$ cannot satisfy 
$f(0$, $x$, $1$, $1) = 0$ for $x>0$, $x \ne 1$. 

\Tqac \ is also proved from the above results. 

Finally, we shall study discriminants $\disc_D = \disc(D)$ for 
$D = C_1^+$, $C_4$, $P_3$, $P_4$ and $P_5$. 
We use $(p_0$, $p_1$, $p_2$, $p_3)$ as a coordinate system of 
$\cH_{4,4}^{s0}$ as before. 
$(p_0$, $p_1$, $p_2$, $p_3)$ corresponds to 
$\displaystyle \sum_{i=0}^3 p_i s_i \in \cH_{4,4}^{s0}$. 

\def\Tqcn{Lemma 3.16}
\proclaim{Lemma 3.16} 
\begin{align*}
 & \disc(C_1^+) = 8p_0^2 - 9 p_1^2 + 12 p_0p_1 + 12 p_0p_2 + 12 p_0p_3, \\
 & \disc(C_4) = -8 p_0^2 - p_1^2 + 4 p_0p_2, \\
 & \disc(P_3) = p_0 + 2p_1 + p_2+ p_3, \\
 & \disc(P_4) = 2p_0 + 2p_1 + p_2, \\
 & \disc(P_5) = p_0. 
\end{align*}
\endproclaim

\Proof 
$\disc(C_1^+) = \disc(C_1)$, since $\Zar(C_1^+) = \Zar(C_1)$. 

If $P = (c_0 \colon c_1 \colon c_2 \colon c_3) 
\in \Delta^0(\cP_{4,4}^{s0+})$, 
then $\displaystyle \disc(P) = \sum_{i=1}^3 c_i p_i$. 
Thus we have $\disc(P_i)$ ($i=3$, $4$, $5$). 

We shall study $\disc(C_4)$. 
Take $\displaystyle f = (1/3u^3)\frh_u^c + v \frq_2
 \in \cF(C_4)$ ($u > 0$, $v \geq 0$). 
The coefficients of $f$ are 
$p_1/p_0 = - 2(u^2+1)/u$, $p_2/p_0 = (u^4+4u^2+1)/u^2$, 
$p_3/p_0 = 2(3u^4+3u^3+2u^2+3u+3)/(3u^2) + v$. 
Eliminate $u$ and $v$ from these relations. 
Then we have $\disc(C_4) = -8 p_0^2 - p_1^2 + 4 p_0p_2 = 0$. 
\end{proof}

\bigskip

{\it Proof of \Tqci.} 
This proof is almost completed. 
What we should do is only to observe the signature of discriminants. 
Then, we find that we can use $p+4$ and $p+2/3$ as separators to 
describe $\cP_{4,4}^{c0+}$ as a union of 
basic semialgebraic sets as (1)---(4) of \Tqci(I). \QED

\bigskip

{\it Proof of \Tqaa (2), \Tqacn \ and \Tqadn (2).} 
This is same as the proof of \Tqaa (1), \Tqabn \ and \Tqadn (1). \QED

\bigskip

{\it Proof of \Tqadeg (2), (3).} 
(2-i) Consider the case $0 \leq t < 1$ or $1 < t \leq 5$. 
Let $F(x,y,z) := 3 \frf_t^{ab}(x$, $y$, $z-x-y$, $-z)$, 
and $f(x,y) := F(x,y,1)$. 
If $\frf_t^{ab}$ is reducible, then $f$ is also reducible. 
Consider the real curve $\Gamma := V_{\R}(F) \subset \P_{\R}^2$. 
Note that $f(y$, $x) = f(x$, $y)$. Since
\begin{align*}
 & f(x,0) = f(0,x) = 8(t+1)(x^2-x+1)^2 > 0, \\
 & F(1,0,z) = F(0,1,z) = 8(t+1)(z^2-z+1)^2 > 0, \\
 & f(1,1) = f(1,-1) = f(-1,1) = -16(t-1)^2 < 0, 
\end{align*}
$\Gamma$ has at least three connected components $\Gamma_1$, $\Gamma_2$, 
$\Gamma_4$ in the 1-st, 2-nd and 4-th quadrant. 
$\Gamma_1$, $\Gamma_2$, $\Gamma_4$ are all bounded. 
This implies $\Gamma$ cannot contain a line. 
Moreover, $\Gamma$ cannot be a union of two quadric curves. 
Thus $V_{\C}(F)$ must be an irreducible curve. 

\smallskip

(3) Consider the case $t > 5$. 
Let $G(x,y,z) := 9 \frf_t^c(x$, $y$, $z-x-y$, $-z)$, 
and $g(x,y) := G(x,y,1)$. Then, 
\begin{align*}
 & g(x,0) = g(0,x) = (t+7)^2(x^2-x+1)^2 > 0, \\
 & G(1,0,z) = G(0,1,z) = (t+7)^2(z^2-z+1)^2 > 0, \\
 & g(1,1) = g(1,-1) = g(-1,1) = -32(t^2+2t-11) < 0. 
\end{align*}
Thus $V_{\C}(G)$ must be an irreducible curve. 

\smallskip

(2-ii) Consider the case $t = 1$. 
Assume that $\frf_1^{ab}$ is reducible. 
Since 
\[\frf_1^{ab}(x,y,1,1) = \frac{1}{3}(x-y)^2(3x^2+2x y+3y^2-8x-8y+8),\]
$\frf_1^{ab}$ must be product of two real quadrics. 
But this is impossible. 
since $(\frf_1^{ab})_{aa}(1$, $1$, $1$, $1) = 0$. \QED

\bigskip

{\it Proof of \Tqade.} 
For $f_t = \frf_t^{ab}$ ($0 < t < 1$ or $1 < t \leq 5$) 
or $f_t = \frf_t^c$ ($t>5$), 
let $F_t(a,b,c,d) := f_t(a^2$, $b^2$, $c^2$, $d^2)$, 
and consider the zero point set $Z_t := V_{\R}(F_t) \subset \P_{\R}^3$.

Let $u$ be a positive root of $t = (3u^2-u+3)/u$ if $t > 5$, 
and $u := 1$ if $0 < t \leq 5$. 
Remember that $f_t(1,1,1,1) = f_t(t,1,1,1) = f_t(0,0,u,1) = 0$. 
Let $s := \sqrt{t}$ and $v := \sqrt{u}$. 
Then $F_t(\pm 1$, $\pm 1$, $\pm 1$, $1) = F_t(\pm s$, $\pm 1$, $\pm 1$, $1) 
 = F_t(0$, $0$, $\pm v$, $1) = 0$. 
Thus, if $0 < t < 1$ or $1 < t \leq 5$, then $\# Z_t = 52$. 
If $t>5$, then $\# Z_t = 64$. 

Assume that $F_t \in \Sigma_{4,8}$. 
Then, there exists $r \in \N$ and $g_1$,$\ldots$, $g_r \in \cH_{4,4}$ 
such that $F_t = g_1^2 + \cdots + g_r^2$. 
If ${\bf a} \in Z_t$, then $g_1({\bf a}) = \cdots = g_r({\bf a}) = 0$, 
since $F_t({\bf x}) \geq 0$ for all ${\bf x} \in \P_{\R}^3$. 
Note that $\dim \cH_{4,4} = 35$. 
So, let's find 35 points ${\bf a}_i \in Z_t$ ($1 \leq i \leq 35$) 
such that there exists no $g \in \cH_{4,4} - \{0\}$ 
which satisfy $g({\bf a}_i) = 0$ for all $1 \leq i \leq 35$. 

Let ${\bf a}_1 := (1 \colon 1 \colon -1 \colon -1)$, 
${\bf a}_2 := (1 \colon 1 \colon 1 \colon s)$, 
${\bf a}_3 := (-s \colon 1 \colon 1 \colon 1)$, 
${\bf a}_4 := (1 \colon -s \colon 1 \colon 1)$, 
${\bf a}_5 := (1 \colon 1 \colon -s \colon 1)$, 
${\bf a}_6 := (1 \colon 1 \colon 1 \colon -s)$, 
${\bf a}_7 := (s \colon -1 \colon 1 \colon 1)$, 
${\bf a}_8 := (s \colon 1 \colon -1 \colon 1)$, 
${\bf a}_9 := (s \colon 1 \colon 1 \colon -1)$, 
${\bf a}_{10} := (-1 \colon s \colon 1 \colon 1)$, 
${\bf a}_{11} := (1 \colon s \colon -1 \colon 1)$, 
${\bf a}_{12} := (1 \colon s \colon 1 \colon -1)$, 
${\bf a}_{13} := (-1 \colon 1 \colon s \colon 1)$, 
${\bf a}_{14} := (1 \colon -1 \colon s \colon 1)$, 
${\bf a}_{15} := (1 \colon 1 \colon s \colon -1)$, 
${\bf a}_{16} := (-1 \colon 1 \colon 1 \colon s)$, 
${\bf a}_{17} := (1 \colon -1 \colon 1 \colon s)$, 
${\bf a}_{18} := (1 \colon 1 \colon -1 \colon s)$, 
${\bf a}_{19} := (s \colon 1 \colon -1 \colon -1)$, 
${\bf a}_{20} := (s \colon -1 \colon 1 \colon -1)$, 
${\bf a}_{21} := (s \colon -1 \colon -1 \colon 1)$, 
${\bf a}_{22} := (-1 \colon s \colon 1 \colon -1)$, 
${\bf a}_{23} := (-1 \colon s \colon -1 \colon 1)$, 
${\bf a}_{24} := (1 \colon s \colon -1 \colon -1)$, 
${\bf a}_{25} := (-1 \colon -1 \colon s \colon 1)$, 
${\bf a}_{26} := (1 \colon -1 \colon s \colon -1)$, 
${\bf a}_{27} := (-1 \colon 1 \colon s \colon -1)$, 
${\bf a}_{28} := (1 \colon -1 \colon -1 \colon s)$, 
${\bf a}_{29} := (-1 \colon 1 \colon -1 \colon s)$, 
${\bf a}_{30} := (-1 \colon -1 \colon 1 \colon s)$, 
${\bf a}_{31} := (v \colon 1 \colon 0 \colon 0)$, 
${\bf a}_{32} := (v \colon 0 \colon -1 \colon 0)$, 
${\bf a}_{33} := (v \colon 0 \colon 0 \colon -1)$, 
${\bf a}_{34} := (0 \colon v \colon 1 \colon 0)$, 
${\bf a}_{35} := (0 \colon v \colon 0 \colon 1)$. 
Take 35 monomials $e_1$,$\ldots$, $e_{35}$ as a basis of $\cH_{4,4}$, 
and denote $g = c_1 e_1 + \cdots + c_{35} e_{35} \in \cH_{4,4}$. 
Let $A = (a_{i,j})$ be $35 \times 35$-matrix such that 
$a_{i,j} = e_j({\bf a}_i)$. Then 
\begin{align*}
 \det A
 & = \pm 549755813888 \, t^{13/2} (t-1)^{23} (t+3)^6 \\
 & \hskip20pt \times u^3 (1+t-2u)(t u + u - 2)(3u^2 - u t - u - 1). 
\end{align*}
Note that $3u^2 - u t - u - 1 = (3u^2- u t - u + 3) - 4 = -4 \ne 0$, 
$t u + u - 2 = 3u^2 + 1 > 0$ and $u>0$. 
There exist no real solutions $1+t-2u = 0$, $t = (3u^2-u+3)/u$. 
Thus $\det A \ne 0$ if $t>0$ and $t \ne 1$. 
This implies there exists no $g \in \cH_{4,4} - \{0\}$ 
which satisfy $g({\bf a}_i) = 0$ for all $1 \leq i \leq 35$. \QED

\bigskip

{\it Proof of \Tqadef.} 
Let $t > 5$. 
We shall show that $\frf_t^c \in \cE(\cP_{4,4}^+)$. 
This is equivalent to $\frh_u^c \in \cE(\cP_{4,4}^+)$ 
for all $u>0$. 

Let $e_1$,$\ldots$, $e_{35}$ be all the monomials in $\cH_{4,4}$, 
and denote $f \in \cH_{4,4}$ as 
$\displaystyle f = \sum_{i=1}^{35} c_i e_i$ ($c_i \in \R$). 
Let $t := (3u^2-u+3)/u$. 
Let $\cK$ be the subspace of all the $f \in \cH_{4,4}$ 
which satisfies the following 34 equalities: \par
\begin{center}
\begin{longtable}{llll}
  $f_a(1,1,1,1)=0$, & $f_b(1,1,1,1)=0$, 
    & $f(t,1,1,1)=0$, & $f_a(t,1,1,1)=0$, \\
  $f_b(t,1,1,1)=0$, & $f(1,t,1,1)=0$, 
    & $f_a(1,t,1,1)=0$, & $f_b(1,t,1,1)=0$, \\
  $f_c(1,t,1,1)=0$, & $f(1,1,t,1)=0$, 
    & $f_a(1,1,t,1)=0$, & $f_b(1,1,t,1)=0$, \\
  $f_c(1,1,t,1)=0$, & $f(1,1,1,t)=0$, 
    & $f_a(1,1,1,t)=0$, & $f_b(1,1,1,t)=0$, \\
  $f_c(1,1,1,t)=0$, & $f(0,0,u,1)=0$, 
    & $f_c(0,0,u,1)=0$, & $f(0,u,0,1)=0$, \\
  $f_b(0,u,0,1)=0$, & $f(0,u,1,0)=0$, 
    & $f_b(0,u,1,0)=0$, & $f(u,0,0,1)=0$, \\
  $f_a(u,0,0,1)=0$, & $f(u,0,1,0)=0$, 
    & $f(u,1,0,0)=0$, & $f_a(u,1,0,0)=0$, \\
  $f(0,0,1,u)=0$, & $f(0,1,0,u)=0$, 
    & $f(0,1,u,0)=0$, & $f(1,0,0,u)=0$, \\
  $f(1,0,u,0)=0$, & $f(1,u,0,0)=0$. & & 
\end{longtable}
\end{center}
The system of these equation can be written 
as $A{\bf c}={\bf 0}$ by a certain $34 \times 35$-matrix $A$, 
i.e. $\cK = \Ker A$. 
Add the vector $(1$, $0$,$\ldots$, $0)$ to the bottom of $A$, 
and make $35 \times 35$-matrix $B$. Then 
\[\det B = \pm t(t+3)(t-1)^{25} 
    u^{12} (u^2-1)^{12} (u^2+1)^2 (12u^4+12u^3+21u^2+10u+9) \ne 0.\]
Thus $\dim \Ker A = 1$, and $\Ker A = \R \cdot \frh_u^c$. 
This implies $\frh_u^c \in \cE(\cP_{4,4}^+)$. \QED

\bigskip

It seems that $\frf_t^{ab} \notin \cE(\cP_{4,4}^+)$ 
for $t < 5$. 
But the author does not have proof. 

\removelastskip\penalty-400\vskip2.5em plus0.3em minus0.3em
{\bf Section 4. Cubic Inequalities of Four Variables}
\par\penalty1000\vskip0.8em plus0.2em minus0.2em
{\bf 4.1. Structure of $\cP_{4,3}^{c0+}$}
\par\penalty1000\vskip0.4em plus0.1em minus0.1em
In this section, we shall study 
$\cP_{4,3}^{c0+} := \cP(\P_+^3$, $\cH_{4,3}^{c0})$. 
We use similar symbols with \S 3. 
To state the main theorem of this section we need to fix some symbols. 
Put 
\begin{align*}
 & S_3 := \sum_{i=0}^3 a_i^3, \hskip15pt 
   S_{2,1,0} := \sum_{i=0}^3 a_i^2 a_{i+1}, \hskip15pt 
   S_{2,0,1} := \sum_{i=0}^3 a_i^2 a_{i+2}, \\
 & S_{1,2,0} := \sum_{i=0}^3 a_i^2 a_{i+3}, \hskip15pt 
   S_{1,1,1} := \sum_{i=0}^3 a_ia_{i+1}a_{i+2}, 
\end{align*}
here we regard $a_{i+4}=a_i$ for all $i \in \Z$. 
We choose $s_0 := S_3 - S_{1,1,1}$, $s_1 := S_{2,1,0} - S_{1,1,1}$, 
$s_2 := S_{2,0,1} - S_{1,1,1}$, $s_3 := S_{1,2,0} - S_{1,1,1}$ as 
a basis of $\cH_{4,3}^{c0}$, and define $\Phi_{4,3}^{c0} : 
  \P_+^3 \cdots\to \P_+^3$ by $\Phi_{4,3}^{c0}({\bf a}) 
 = \big(s_0({\bf a}):s_1({\bf a}):s_2({\bf a}):s_3({\bf a})\big)$. 
The coordinate system of $A = \P_{\R}^3$ is denoted by 
$(a_0 \colon a_1 \colon a_2 \colon a_3)$ or $(a \colon b \colon c \colon d)$, 
and the coordinate system of $\P((\cH_{4,3}^{c0})^{\vee})$ is denoted by 
$(x_0 \colon x_1 \colon x_2 \colon x_3)$. 
We represent $f \in \cH_{4,3}^{c0}$ as 
$f = p_0 s_0 + \cdots + p_3 s_3$ ($p_i \in \R$), 
and the coordinate system of $\cH_{4,3}^{c0}$ is denoted by 
$(p_0$, $p_1$, $p_2$, $p_3)$. 
If $f \in \cP_{4,3}^{c0+}$, then $s_0 \geq 0$. 
When $p_0=1$, we say $f$ is ${\it monic}$. 
When $p_0=0$, we say $f$ lies {\it at infinity}. We denote 
\[\breve{\cP}_{4,3}^{c0+} 
    := \big\{ f \in \cP_{4,3}^{c0+} \; \big| \; \hbox{$f$ is monic}\big\}.\]
The characteristic variety is written by 
$X_{4,3}^{c0+} := \Phi_{4,3}^{c0}(\P_+^3)$. Let
\begin{align*}
 & A_c^+ := \big\{ (a_0 \colon a_1 \colon a_2 \colon 1) \in \P_+^3 \; \big| 
             \; \hbox{$0 \leq a_i \leq 1$ ($i=0$, $1$, $2$)} \big\}, \\
 & E_2 := \big\{ (a \colon b \colon a \colon b) \in \P_+^3 \; \big| \; 
             \hbox{$a$, $b \in \R_+$} \big\}, \\
 & E_3 := \big\{ (a_0 \colon a_1 \colon a_2 \colon a_3) 
             \in \P_+^3 \; \big| \; \hbox{$a_0+a_2 = a_1+a_3$} \big\}, \\
 & B_0 := \big\{ (0 \colon s \colon t \colon 1) \in \P_+^3 \; \big| \;  
             \hbox{$s>0$, $t>0$, $(s,t) \ne (1,0)$, $t \ne s+1$} \big\}, \\
 & \overline{B_0} := \big\{ (0 \colon s \colon t \colon u) \in \P_+^3 \; 
          \big| \; \hbox{$(s \colon t \colon u) \in \P_+^2$} \big\}, \\
 & S := \Phi_{4,3}^{c0}(B_0) \subset X_{4,3}^{c0+}, \\
 & C := \big\{\Phi_{4,3}^{c0}(0 \colon 0 \colon t \colon 1) \in \P_+^3 
             \; \big| \; \hbox{$t>0$}\big\} \subset X_{4,3}^{c0+}, \\
 & L := \big\{\Phi_{4,3}^{c0}(0 \colon t \colon 0 \colon 1) \in \P_+^3 
             \; \big| \; \hbox{$t>0$}\big\} \subset X_{4,3}^{c0+}, \\
 & P_1 := (1 \colon 0 \colon 0 \colon 0) 
        = \Phi_{4,3}^{c0}(0 \colon 0 \colon 0 \colon 1)
            \in X_{4,3}^{c0+}, \\
 & P_2 := (1 \colon 0 \colon 1 \colon 0) 
        = \big\{ \Phi_{4,3}^{c0}(a \colon b \colon a \colon b) \in \P_+^3 \; 
           \big| \; \hbox{$a$, $b \in \R_+$} \big\}
           \in X_{4,3}^{c0+}, \\
 & P_3 := (2 \colon 1 \colon 0 \colon 1) \in X_{4,3}^{c0+}. \\
\end{align*}
We denote $\cF(P_i)$, $\cF(C)$, $\cF(S)$ by 
$\cF_{P_i}$, $\cF_C$ and $\cF_S$. 
As we will prove in \Tqcbb, 
$\partial \cP_{4,3}^{c0+} = \cF_S \cup \cF_C \cup \cF_{P_1} \cup \cF_{P_2}$, 
So, we need two discriminants $\disc_C$ and $\disc_S$ which are defining 
equations of $\Zar(\cF_C)$ and $\Zar(\cF_S)$. 
Note that $\disc_{P_1} = p_0$ and $\disc_{P_2} = p_0+p_2$. 
Using PC, we obtain them. 
Since $\disc_S$ is somewhat complicated polynomial, 
we also define $d_S$ to describe $\disc_S$.
\begin{align*}
 & \disc_C(p_0,p_1,p_3) = 27 p_0^4 + 4 p_0 p_1^3 + 4 p_0 p_3^3 
        - p_1^2 p_3^2 - 18 p_0^2 p_1 p_3
    = \Disc_3(p_0,p_1,p_3,p_0), \\
 & d_S(p_0,p_2,q,r) \\
 & := (p_0 - p_2 - q)^2 
      (13 p_0^2 - 2 p_0 p_2 + p_2^2 + 2 p_0 q + 2 p_2 q)^2 \\
 & \hskip20pt  
     (104 p_0^3 + 100 p_0^2 p_2 - 4 p_0 p_2^2 + 36 p_0^2 q 
        + 36 p_0 p_2 q - p_0 q^2 - p_2 q^2 + 8 q^3) \\
 & \hskip20pt  + (17173 p_0^7 - 121 p_0^6 p_2 - 5639 p_0^5 p_2^2 
      + 7651 p_0^4 p_2^3 - 3489 p_0^3 p_2^4 + 469 p_0^2 p_2^5 \\
 & \hskip20pt  
      - 45 p_0 p_2^6  + p_2^7 + 6250 p_0^6 q + 10028 p_0^5 p_2 q 
      + 3142 p_0^4 p_2^2 q - 1368 p_0^3 p_2^3 q - 746 p_0^2 p_2^4 q \\
 & \hskip20pt 
      - 20 p_0 p_2^5 q - 6 p_2^6 q  + 898 p_0^5 q^2 +  7230 p_0^4 p_2 q^2 
      + 1748 p_0^3 p_2^2 q^2 - 1572 p_0^2 p_2^3 q^2 \\
 & \hskip20pt 
      - 86 p_0 p_2^4 q^2 - 26 p_2^5 q^2 + 2780 p_0^4 q^3 - 368 p_0^3 p_2 q^3 
      + 1448 p_0^2 p_2^2 q^3 - 496 p_0 p_2^3 q^3 \\
 & \hskip20pt 
      + 28 p_2^4 q^3 + 518 p_0^3 q^4 + 1018 p_0^2 p_2 q^4 - 190 p_0 p_2^2 q^4 
      + 78 p_2^3 q^4 + 164 p_0^2 q^5 \\
 & \hskip20pt 
      + 168 p_0 p_2 q^5 + 4 p_2^2 q^5) r^2  \\
 & \hskip20pt  + (2495 p_0^5 - 317 p_0^4 p_2 - 1886 p_0^3 p_2^2 
      + 842 p_0^2 p_2^3 - 81 p_0 p_2^4 + 3 p_2^5 + 1768 p_0^4 q \\
 & \hskip20pt 
      + 4 p_0^3 p_2 q - 988 p_0^2 p_2^2 q + 380 p_0 p_2^3 q - 12 p_2^4 q 
      + 291 p_0^3 q^2 + 897 p_0^2 p_2 q^2 - 463 p_0 p_2^2 q^2 \\
 & \hskip20pt  
      + 83 p_2^3 q^2 + 226 p_0^2 q^3 + 92 p_0 p_2 q^3 - 38 p_2^2 q^3 
      - p_0 q^4 - p_2 q^4) r^4 \\
 & \hskip20pt  + (95 p_0^3 + 65 p_0^2 p_2 - 43 p_0 p_2^2 + 3 p_2^3 
      + 98 p_0^2 q - 20 p_0 p_2 q - 6 p_2^2 q - 4 p_0 q^2) r^6 \\
 & \hskip20pt   + (-3p_0 + p_2)r^8, \\
 & \disc_S(p_0,p_1,p_2,p_3) 
     = \frac{1}{4} d_S(p_0, p_2, p_1+p_3, p_1-p_3). 
\end{align*}
Since $\disc_C(p_0,p_1,p_3)$ has an obstacle branch in the 
first quadrant $p_1/p_0>0$, $p_3/p_0>0$, we put 
\[d_C(x,z) := \begin{cases}
      \disc_C(1,x,z) & \hbox{(if $x < 0$ or $z < 0$)} \\
      1              & \hbox{(if $x \geq 0$ and $z \geq 0$)} 
\end{cases}\]
to avoid complexity. $d_C(x,z) \geq 0$ implies $\disc_C(1,x,z) \geq 0$ 
or `$x \geq 0$ and $z \geq 0$'. 
Thus, $d_C(x,z) \geq 0$ defines a convex domain, 
but $\disc_C(1,x,z) \geq 0$ does not. 
The following $\eta(x,y)$ is a nice separator 
whose property is explained in \Tqcbj. 
\begin{align*}
 & \eta(x,y) :=  61 + 62 x + 56 y + 32 x^2 + 30 x y - 6 y^2 \\
 & \hskip50pt 
    + 9 x^3 + 4 x^2 y - 6 x y^2 - 16 y^3 
    + x^4 - 4 x^2 y^2 - 6 x y^3 + y^4 - x^3 y^2. 
\end{align*}
We also need two constants $\kappa_1$, $\kappa_2$. 
Let $\kappa_1 := 0.0129074031\cdots$ be a root of 
\begin{align*}
 & 817808203 x^6 - 546807084 x^5 + 129155640 x^4 \\
 & \hskip50pt - 13342016 x^3 + 556080x^2 - 10176 x + 64 = 0, 
\end{align*}
and $\kappa_2 := 0.0318925844\cdots$ be a root of 
\[43042537 x^6 - 4514514 x^5 - 188769 x^4 - 38684 x^3 + 4119 x^2 
  - 114 x + 1 = 0.\] 
The aim of this section is to prove the following theorem. 

\def\Tqcbd{Theorem 4.1}%
\proclaim{Theorem 4.1} 
{\rm (I)} {\sl Take a monic $f = s_0 + p_1 s_1 + p_2 s_2 + p_3 s_3 
  \in \breve{\cH}_{4,3}^{c0}$. 
Then, $f({\bf a}) \geq 0$ for all ${\bf a} \in \R_+^4$, 
if and only if one of the followings holds:} \par
{\parindent=20pt
\Item{(1)} {\sl $p_2 = -1$ and $8(p_1+p_3) \geq (p_1-p_3)^2$.}
\Item{(2)} {\sl $-1 < p_2 \leq 3$, $\disc_S(1,p_1,p_2,p_3) \geq 0$ 
and $d_C(p_1,p_3) \geq 0$.}
\Item{(3)} {\sl $p_2>3$, $\kappa_1 (p_1 + p_3) + \kappa_2 p_2 \geq 1$, 
$\disc_S(1,p_1,p_2,p_3) \geq 0$, and $d_C(p_1,p_3) \geq 0$.} 
\Item{(4)} {\sl $p_2>3$, $\kappa_1 (p_1 + p_3) + \kappa_2 p_2 < 1$, 
$\eta(p_1+p_3, p_2) > 0$, $\disc_S(1,p_1,p_2,p_3) \geq 0$, 
and $d_C(p_1,p_3) \geq 0$.} 
\Item{(5)} {\sl $p_2>3$, $\kappa_1 (p_1 + p_3) + \kappa_2 p_2 < 1$, 
$\eta(p_1+p_3, p_2) \leq 0$, and $d_C(p_1,p_3) \geq 0$.} 

}
{\rm (II)} {\sl 
Let's denote $f = p_0 s_0 + p_1 s_1 + p_2 s_2 + p_3 s_3$. 
Then, all the discriminants of $\cP_{4,3}^{c0+}$ are 
$\disc_S(p_0,p_1,p_2,p_3)$, $\disc_C(p_0,p_1,p_3)$, 
$\disc_{P_1}\allowbreak  = p_0$, and $\disc_{P_2} = p_0+p_2$.}

{\rm (III)} {\sl If $f \in \cH_{4,3}^{c0}$ satisfies 
$f(0$, $s$, $t$, $1) \geq 0$ for all $s$, $t \in \R_+$, 
then $f \in \cP_{4,3}^{c0}$. }
\endproclaim

This theorem will be proved after \Tqcbi. 

\def\Tqcbbb{Proposition 4.2}%
\proclaim{Proposition 4.2} 
\begin{align*}
 & s_0(a_0,a_1,a_2,a_3) = \frac{1}{3} \sum_{i=0}^3 
     (a_i^3 + a_{i+1}^3 + a_{i+2}^3 - 3 a_i a_{i+1} a_{i+2}) \geq 0, \\
 & s_2(a_0,a_1,a_2,a_3) = (a_0-a_1+a_2-a_3)(a_0 a_2 - a_1 a_3), \\
 & s_3(a_0,a_1,a_2,a_3) = s_1(a_0,a_3,a_2,a_1), \\
 & s_0 - s_2 = \frac{1}{3}\sum_{i=0}^3 
       (a_i^3 + a_i^3 + a_{i+2}^3 - 3 a_i^2a_{i+2}) \geq 0, \\
 & s_0 + 2s_2 
   = \sum_{i=0}^3 
     (a_i^2a_{i+2} + a_{i+1}^3 + a_i a_{i+2}^2 - 3 a_i^2 a_{i+1} a_{i+2})
         \geq 0, \\
 & 2s_1 + s_2
   = \sum_{i=0}^3 (a_i^2a_{i+1} + a_{i+1}^2a_{i+2} 
      + a_{i+2}^2a_i - 3 a_i a_{i+1} a_{i+2}) \geq 0, \\
 & 2s_3 + s_2
   = \sum_{i=0}^3 (a_i a_{i+1}^2 + a_{i+1} a_{i+2}^2 
      + a_{i+2} a_i^2 - 3 a_i a_{i+1} a_{i+2}) \geq 0, \\
 & s_0 - s_1
   = \frac{1}{3}\sum_{i=0}^3 
      (a_i^3 + a_i^3 + a_{i+1}^3 - 3 a_i^2 a_{i+1}) \geq 0, \\
 & s_0 - s_3 = \frac{1}{3}\sum_{i=0}^3
      (a_i^3 + a_{i+1}^3 + a_{i+1}^3 - 3 a_i a_{i+1}^2) \geq 0, \\
 & s_1 + s_3 = (a_0+a_2)(a_1-a_3)^2 + (a_1+a_3)(a_0-a_2)^2 \geq 0. 
\end{align*}
\endproclaim

\Proof 
These follow from direct calculations. 
\end{proof}

Thus $X_{4,3}^{c0+}$ is a subset of a cube defined by 
$-1/2 \leq s_1/s_0 \leq 1$, $-1/2 \leq s_2/s_0 \leq 1$, 
$-1/2 \leq s_3/s_0 \leq 1$. 
Note that $s_1$, $s_2$ and $s_3$ are not PSD. 
For example $s_1(1/100$, $1/2$, $1/10$, $1) = -229/20000 < 0$. 
The rational map $\Phi_{4,3}^{c0} : \P_+^3 \cdots\to X_{4,3}^{c0+}$ splits as 
\[\Phi_{4,3}^{c0} : \P_+^3 \mapr{\pi} \P_+^3/(\Z/4\Z) 
   \mapr{\Psi_{4,3}^{c0}} X_{4,3}^{c0+}.\]
It is easy to see that 
$\Psi_{4,3}^{c0} : \P_+^3/(\Z/4\Z) \cdots\to X_{4,3}^{c0+}$ is 
a birational map, but is not holomorphic 
at a singular point $\pi(1 \colon 1 \colon 1 \colon 1)$. 
We shall provide more precise structure of $X_{4,3}^{c0+}$ at \Tqcbb. 
The following $\fre_{s,t}(a_0,a_1,a_2,a_3) \in \cH_{4,3}^{c0}$ 
($s$, $t \in \R$) has a possibility to be an extremal element. 
But there exists $(s$, $t)$ such that $\fre_{s,t}$ is not PSD. 

\def\Tqcbbd{Proposition 4.3}%
\proclaim{Proposition 4.3} 
{\sl For $(u \colon v \colon w) 
 \in \P_+^2 - \big\{(1 \colon 0 \colon 1)\big\}$, let 
\begin{align*}
 & g_0^h(u,v,w) := -v\big(u w v^2 - (u+w)(u^2+w^2)v + u w(u-w)^2\big), \\
 & g_1^h(u,v,w) := u v^4 - w(u+2w)v^3 -2u w(u-w)v^2 \\
 & \hskip80pt       - u(2u^3+u^2 w-3w^3)v + w(u^2-w^2)^2, \\
 & g_2^h(u,v,w) := v\big(v^4+(2u^2-3u w+2w^2)v^2 \\
 & \hskip80pt - (u+w)(u^2+w^2)v + (u-w)^2(u^2-u w+w^2)\big), \\
 & g_3^h(u,v,w) := g_1^h(w,v,u), \\
 & \fre_{u,v,w}^h({\bf a}) := \sum_{i=0}^3 g_i^h(u,v,w) s_i({\bf a}). 
\end{align*}
For simplicity, put $g_i(s,t) := g_i^h(s,t,1)$ and 
$\fre_{s,t}({\bf a}) := \fre_{s,t,1}^h({\bf a})$. 
Then the followings hold: } \par
{\parindent=20pt
\Item{\rm (1)} $\fre_{w,v,u}^h - \fre_{u,v,w}^h 
  = (u-w)(v^2-(u+w)^2)((u-w)^2+2(u+w)v+v^2) (s_1-s_3)$. 
\Item{\rm (2)} 
$\fre_{t,1,0}^h = t \fre_{0,t,1}^h - (t^2-1)(t^2+1)^2 s_2$. 
\Item{\rm (3)} {\sl Assume that $s>0$, $t>0$, $t \ne s+1$, $g_0(s,t)>0$  
and $\fre_{s,t} \in \cP_{4,3}^{c0+}$. 
If $f \in \cP_{4,3}^{c0+}$ satisfies $f(0,s,t,1)=0$, 
then there exists $\alpha \geq 0$ such that $f = \alpha \fre_{s,t}$. 
Especially, $\fre_{s,t} \in \cE(\cP_{4,3}^{c0+})$. }
\Item{\rm (4)} {\sl Assume that $s=0$, $t>0$, $t \ne 1$ 
and $\fre_{0,t} \in \cP_{4,3}^{c0+}$. 
If $f \in \cP_{4,3}^{c0+}$ satisfies $f(0,0,t,1)=0$ and 
$\displaystyle \frac{\partial}{\partial b}f(0,0,t,1)=0$, 
then there exists $\alpha \geq 0$ such that $f = \alpha \fre_{0,t}$. 
Especially, $\fre_{0,t} \in \cE(\cP_{4,3}^{c0+})$. }
\Item{\rm (5)} {\sl Assume that $u>0$, $v>0$, 
and $\fre_{u,v,0}^h \in \cP_{4,3}^{c0+}$. 
If $f \in \cP_{4,3}^{c0+}$ satisfies $f(0$, $u$, $v$, $0)=0$ and 
$\displaystyle \frac{\partial}{\partial d}f(0,u,v,0)=0$, 
then there exists $\alpha \geq 0$ such that $f = \alpha \fre_{u,v,0}^h$. 
Especially, $\fre_{u,v,0}^h \in \cE(\cP_{4,3}^{c0+})$. }
\Item{\rm (6)} {\sl If $t = s+1$, then 
\begin{align*}
 \fre_{s,s+1}(a,b,c,d) 
 & = (s+1)(s^2+1)^2 (a-b+c-d)^2(a+b+c+d) & (*) \\
 & = (s+1)(s^2+1)^2 \fre_{0,1}(a,b,c,d). \
\end{align*}
If $f \in \cP_{4,3}^{c0+}$ satisfies $f(0,0,1,1)=0$ and 
$f(0,1,2,1) = 0$, then there exists $\alpha \geq 0$ such 
that $f = \alpha \fre_{0,1}$. 
Especially, $\fre_{s,s+1} \in \cE(\cP_{4,3}^{c0+})$. }
\Item{\rm (7)} {\sl If $g_0(s$, $t) < 0$, then 
$\fre_{s,t} \notin \cE(\cP_{4,3}^{c0+})$ and 
$-\fre_{s,t} \notin \cE(\cP_{4,3}^{c0+})$. }
\Item{\rm (8)} {\sl $\fre_{1,0} = \fre_{1,0,1}^h$ is a zero polynomial. } 

}
\endproclaim

\Proof 
Denote $\displaystyle f_a(a,b,c,d) 
= \frac{\partial}{\partial a}f(a,b,c,d)$ and so on. 

\smallskip

(1), (2) and (8) follow from direct calculation. 

\smallskip

(3) Assume that $f = p_0s_0 + p_1s_1 + p_2s_2 + p_3s_3 
\in \cP_{4,3}^{c0+}$ satisfies $f(0,s,t,1)=0$. 
Then $f_b(0,s,t,1)=0$ and $f_c(0,s,t,1)=0$ hold. 
Let $a_{0,j} = s_j(0,s,t,1)$, $a_{1,j} = (s_j)_b(0,0,t,1)$, 
$a_{2,j} = (s_j)_c(0,s,t,1)$, and $A = (a_{i,j})$. Then 
\[A = \left(
\begin{matrix}
  s^3+t^3-s t+1 & t(s^2-s+t) & s(1+s-t) & t(s t-s+1) \\
  3s^2-t & (2s-1)t & 2s-t+1 & t(t-1) \\
  3t^2-s & s^2-s+2t & -s & 2s t-s+1 
\end{matrix}
\right).\]
Let $B$ be the square matrix add $(1,0,0,0)$ above $A$. 
Then $\det B = (t-s-1) g_0(s,t) \ne 0$. 
Note that $\fre_{s,t} \in \Ker A$. 
Thus, $\Ker A = \R \cdot \fre_{s,t}$. 

\smallskip

(4), (5) Same with (3).  

\smallskip

(6) $(*)$ follows from direct calculation. 
Assume that $f \in \cP_{4,3}^{c0+}$ satisfies 
$f(0$, $0$, $1$, $1)=0$ and $f(0,1,2,1) = 0$. 
Then $f(0,0,1,1)=0$, $f_a(0,0,1,1)=0$ and $f_a(0,1,2,1)=0$. 
then $f_c(0,0,t,1)=0$ holds. 
By the same method as (3), we have the conclusion. 

\smallskip

(7) We may assume $t \ne s+1$. 
 If $\fre_{s,t} \in \cE(\cP_{4,3}^{c0+})$, 
then $g_0(s,t) = \fre_{s,t}(0,0,0,1) \geq 0$. 
On the other hand, 
$\fre_{s,t}(0,0,1,1) = (s+1)(t-s-1)^2((s-1)^2+t^2) > 0$. 
Thus $-\fre_{s,t} \notin \cE(\cP_{4,3}^{c0+})$. 
\end{proof}

\bigskip

The condition that $\fre_{s,t} \in \cE(\cP_{4,3}^{c0+})$ will 
be determined at \Tqccbd. 

\def\Tqcbb{Lemma 4.4}%
\proclaim{Lemma 4.4} 
{\sl Let ${\bf 1} = (1 \colon 1 \colon 1 \colon 1) \in \P_+^3$, 
$\displaystyle Z 
   := A_c^+ - \{{\bf 1}\} - \bigcup_{\tau \in \Z/4\Z} \tau(E_2 \cup E_3)$ and }
\begin{align*}
 & f_{4,3}^{c0}(x_0, x_1, x_2, x_3) \\
 & := (x_1^3 - x_0 x_1 x_3 + x_3^3)^2 
     - x_2 (x_1^3 - x_0 x_1 x_3 + x_3^3)
        (x_0^2 + 3 x_1^2 - 4 x_1 x_3 + 3 x_3^2)  \\
 & \hskip20pt  + x_2^2 \big(x_0^2 (x_1^2 - x_1 x_3 + x_3^2)
    + 2 x_0 x_1 x_3 (x_1 + x_3) \\
 & \hskip80pt
    + x_1^4 - 7 x_1^3 x_3 + 9 x_1^2 x_3^2 - 7 x_1 x_3^3 + x_3^4 \big) \\
 & \hskip20pt  + x_2^3 \big(2 x_0 x_1^2 - x_0 (4 x_1^2 + x_1 x_3 + 2 x_3^2) 
                  + (x_1 + x_3) (x_1^2 - 3 x_1 x_3 + x_3^2)\big) \\
 & \hskip20pt  + x_2^4 (x_1^2 + x_1 x_3 + x_3^2). 
\end{align*}
{\sl Then, the followings hold: } \par
{\parindent=20pt
\Item{(1)} {\sl $\Phi_{4,3}^{c0} : A_c^+ \cdots\to X_{4,3}^{c0+}$ is 
a birational map whose all the exceptional sets 
are $\Phi_{4,3}^{c0}(E_2) = P_2$ 
and $\Phi_{4,3}^{c0}(E_3) = P_3$. 
$\Phi_{4,3}^{c0} \colon Z \lto \Phi_{4,3}^{c0}(Z)$ is an isomorphism. 
$\Bs \Phi_{4,3}^{c0} = \{{\bf 1}\}$ and we can regard 
$\Phi_{4,3}^{c0}({\bf 1})$ as the closed line segments $[P_2P_3]$.}
\Item{(2)} {\sl $\Zar(\partial X_{4,3}^{c0+}) \subset V_{\R}(f_{4,3}^{c0})$, 
$\Phi_{4,3}^{c0}(\overline{B_0}) = \partial X_{4,3}^{c0+}$ 
and $S$ is non-singular.}
\Item{(3)} {\sl $\Delta^0(X_{4,3}^{c0}) = \{P_1, P_2\}$, 
$\Delta^1(X_{4,3}^{c0}) = \{C, (P_1P_2)\}$ 
and $\Delta^2(X_{4,3}^{c0}) = \{S\}$. }
\Item{(4)} {\sl Let $\cL_{(0 \colon s \colon t \colon 1)}^{c0+}$ 
be the local cone of $\cP_{4,3}^{c0+}$ 
at $(0 \colon s \colon t \colon 1)$. 
Take $(0 \colon s \colon t \colon 1) \in B_0$. 
If $\fre_{s,t}$ is PSD, then $\fre_{s,t} \in \cE(\cP_{4,3}^{c0+})$ and 
\[\cL_{(0 \colon s \colon t \colon 1)}^{c0+} = \R_+ \cdot \fre_{s,t}.\]
If $\fre_{s,t}$ is not PSD, 
then $\cL_{(0 \colon s \colon t \colon 1)}^{c0+} = 0$. }

}
\endproclaim

\Proof 
(1), (2) and (3) $\Bs \Phi_{4,3}^{c0} = \{{\bf 1}\}$ is trivial. 
Since $\Phi_{4,3}^{c0}(1,1+b,1+c,1+d)
 = (2(c^2+(b-d)^2) + (b-c+d)^2 : 
         c^2+(b-d)^2 : (b-c+d)^2, c^2+(b-d)^2)
      + \hbox{(higher degree terms)}$, 
we can regard $\Phi_{4,3}^{c0}(1,1,1,1)$ is a line segment $[P_2P_3]$. 

$\Phi_{4,3}^{c0}(E_2) = P_2$ and $\Phi_{4,3}^{c0}(E_3) = P_3$ are 
obtained by the direct calculation. 

$A_c^+$ is a fundamental domain of $\Phi_{4,3}^{c0}$. 
It is easy to see that $\Phi_{4,3}^{c0} \colon \P_+^3 \to \P_+^3$ is 
a generically finite map of degree $4$. 
The Jacobian of $\Phi_{4,3}^{c0}$ is equal to 
\[J_P := - (a-b+c-d)^3 ((a-c)^2+(b-d)^2)^2 (a+b+c+d).\]
Thus $J_P \ne 0$ on $Z$. 
Therefore $\Phi_{4,3}^{c0} \colon Z \lto \Phi_{4,3}^{c0}(Z)$ is 
an isomorphism.  
This also implies $S \subset \partial X_{4,3}^{c0+}$, 
$C \subset \partial X_{4,3}^{c0+}$, $L \subset \partial X_{4,3}^{c0+}$, 
and $\{P_1$, $P_2$, $P_3\} \subset \partial X_{4,3}^{c0+}$. 

We obtain $f_{4,3}^{c0}$ by eliminating $a$, $b$, $c$ from 
$x_i = s_i(a,b,c,0)$ ($i=0$, $1$, $2$, $3$). Using PC, we fave 
\[f_{4,3}^{c0}\big(\Phi_{4,3}^{c0}(a,b,c,d) \big) 
  = a b c d(a-b+c-d)^4 (a+b+c+d)^2 ((a-c)^2 + (b-d)^2)^4 
   \geq 0.\]
Thus $\partial X_{4,3}^{c0} \subset V_{\R}(f_{4,3}^{c0})$, 
and $\Phi_{4,3}^{c0}(\overline{B_0}) = \partial X_{4,3}^{c0+}$. 
Since $J_P \ne 0$ on $B_0$, we have $\Sing(S) = \emptyset$. 

Since $C = \big\{(t^3+1 : t^2 : 0 : t)$ $\big|$ $\hbox{$t>0$}\big\}$, 
$C$ is a cubic curve desined by $x_1^3 + x_2^3 = x_0 x_1 x_3$, $x_2 = 0$ 
and $(x_1+x_3)/x_0 > 0$. 
Note that $C$ has a node at $P_1$ ($t=0$ and $\infty$). 
But $P_3$ ($t=1$) is a non-singular point of $C$. 

Since $L = \big\{(t^3+1 : 0: t(t+1) : 0)$ $\big|$  $\hbox{$t>0$}\big\}$, 
$L$ is a line segment $(P_1P_2]$ desined by $x_1^3 + x_2^3 = x_0 x_1 x_3$, $x_2 = 0$ 
and $0 < x_2/x_0 \leq 1$. 

Thus $\Sing(\partial X_{4,3}^{c0}) = C \cup (P_1P_2) \cup \{P_1, P_2\}$, 
This implies (3). 

\smallskip

(4) follows from \Tqcbbd(4). 
\end{proof}

It is easy to draw a graph of $X_{4,3}^{c0+}$ using Mathematica. 
But it may present incorrect impression. 
It seems that $X_{4,3}^{c0+}$ is a convex set. 
But it is not true. The following observation show us 
that $X_{4,3}^{c0+}$ is not convex near $(1 \colon 0 \colon 0 \colon 0)$. 
Cut $\partial X_{4,3}^{c0+}$ by the plane $V_{\R}(x_1-x_3)$. Note that 
\[f_{4,3}^{c0}(1,x,y,x) = 
 x^2 (2x-3y-1)(2x^3 + x^2y - y^3 - x^2 + 2 y^2 - y).\]
The graph of $V_{\R}(2x^3 + x^2y - y^3 - x^2 + 2 y^2 - y)$ is not 
convex near $(x,y)=(0,0)$. 
Thus $X_{4,3}^{c0+}$ is not convex. 
This also implies that $\fre_{s,t} \notin \cP_{4,3}^{c0+}$ for 
some $(0 \colon s \colon t \colon 1) \in B_0$. 

\bigskip

It is also possible to obtain $\fre_{s,t}$ by the method 
explained in \cite[Remark 1.28]{RefAc}. 

Let $\displaystyle f_i(x_0,x_1,x_2,x_3) := 
 \frac{\partial}{\partial x_i} f_{4,3}^{c0}(x_0,x_1,x_2,x_3)$ and 
\begin{align*}
 & h_i(s,t) := f_i\big(\Phi_{4,3}^{c0}(0,s,t,1)\big), \\
 & g_c(s,t) := s t (t-s-1)^2 (s+t+1) ((s-1)^2+t^2)^2. 
\end{align*}
Then $h_i(s,t) = g_c(s,t) g_i(s,t)$ ($i=0$, $1$, $2$, $3$). 
Thus we have $\displaystyle \fre_{s,t} = \sum_{i=0}^3 g_i(s,t) s_i$. 
We define a rational map $G : \overline{B_0} \cdots\to \P(\cH_{4,3}^{c0})$ by
\[G(0,s,t,1) 
 := \big(g_0(s,t) \colon g_1(s,t) \colon g_2(s,t) \colon g_3(s,t)\big).\]
Note that $(0 \colon 1 \colon 0 \colon 1) \in \Bs G$. 
If $\fre_{s,t} \in \cP_{4,3}^{c0+}$, 
then $G(0,s,t,1) = \fre_{s,t} \in \cF_S$. 
We can extend $G$ to 
$G^S : \partial \P_+^3 \cdots\to \P(\cH_{4,3}^{c0})$ by 
$G^S(x \colon y \colon 1 \colon 0) = G^S(y \colon 1 \colon 0 \colon x) 
= G^S(1 \colon 0 \colon x \colon y) = G^S(0 \colon x \colon y \colon 1) 
 := G(0$, $x$, $y$, $1)$. 

\def\Tqcbf{Lemma 4.5}%
\proclaim{Lemma 4.5} 
{\rm (1)} $\partial \cP_{4,3}^{c0+} = \cF_{P_1} \cup \cF_{P_2}
 \cup \cF_S \cup \cF_C.$ \par
{\parindent=20pt
\Item{(2)} {\sl $B_0$ is a test set of $\cP_{4,3}^{c0+}$. 
In other words, if $f \in \cH_{4,3}^{c0}$ satisfies $f(0,s,t,1) = 0$ 
for all $s \geq 0$, $t \geq 0$, then $f({\bf a}) \geq 0$ for 
all ${\bf a} \in \R_+^4$. }

}
\endproclaim

\Proof 
(1) $\displaystyle \partial \cP_{4,3}^{c0+}
  = \bigcup_{D \in \Delta(X_{4,3}^{c0+})} \cF(D)$ by \Tqbd(1). 
Let $D_3 := \Int(X_{4,3}^{c0+}) \in \Delta^3(X_{4,3}^{c0+})$. 
Sice $\Zar(D_3) = \P_{\R}^3$, $\cF(D_3)$ is not a face component. 

$\Zar\big(\cF((P_1P_2))\big)$ is two dimensional plane 
defined by $p_0 = p_2 = 0$. 
Thus, $\cF((P_1P_2))$ is not a face component. 
Thus we have the conclusion. 

(2) By \Tqcbb(2) and \Tqbby, we have the conclusion. 
\end{proof}

Note that (III) of \Tqcbd \ follows from the above proposition. 

\def\Tqcbg{Lemma 4.6}%
\proclaim{Lemma 4.6} 
{\sl We regard as $\cH_{4,3}^{c0+} = \R^4$ by identifying 
$(p_0$, $p_1$, $p_2$, $p_3) \in \R^4$ with 
$\displaystyle \sum_{i=0}^3 p_i s_i \in \cH_{4,3}^{c0+}$. 
Then, } \par
{\parindent=20pt
\Item{(1)} {\sl $\Zar\big(\cF_{P_1}\big) = V_{\R}(p_0)$. 
Thus $\cF_{P_1} = 
  \big\{f \in \cP_{4,3}^{c0}$ $\big|$ $f$ is at infinity $\big\}$. }
\Item{(2)} $\Zar\big(\cF_{P_2}\big) = V_{\R}(p_0+p_2)$. 
\Item{(3)} {\sl $\disc_S\big(g_0(s,t), \, g_1(s,t), \, g_2(s,t), 
     \, g_3(s,t)\big) = 0$ for all $s$, $t \in \R$. }
\Item{(4)} {\sl $\disc_S\big(g_0(s,t), \, g_3(s,t), \, g_2(s,t), 
     \, g_1(s,t)\big) = 0$ for all $s$, $t \in \R$. }
\Item{(5)} $\Zar\big(\cF_C\big) = V_{\R}(\disc_C)$. 
\Item{(6)} $\Zar\big(\cF_S\big) = V_{\R}(\disc_S)$. 

}
\endproclaim

\Proof 
(1) and (2) are trivial. 

(3) and (4) follow from direct calculation. 

(5) follows from study of $\cP_{3,3}^{c+}$. See \cite[\S 3]{RefAc}. 

(6) follows from (3). 
\end{proof}

Now, we shall observe $\cF_{P_2}$. 
Remember that $\fre_{1,0} = 0$. 
In other word, $g_i(1$, $0) = 0$ ($i=0$, $1$, $2$, $3$). 
This $\fre_{1,0}$ corresponts to $\cF_{P_2}$. 
Put $\displaystyle g_i^{P_2}(c) := \lim_{h \to 0} \frac{g_i(c h+1, h)}{4h^2}$. 
Then $g_0^{P_2}(c)=1$, $g_1^{P_2}(c)=c(c-2)$, $g_2^{P_2}(c)=-1$, 
$g_3^{P_2}(c)=c(c+2)$. 

\def\Tqcbh{Lemma 4.7}%
\proclaim{Lemma 4.7} 
{\sl For $c \in \R$, let 
\[\fre_c^{P_2} := s_0 + c(c-2) s_1 - s_2 + c(c+2) s_3,\]
and $\fre_{\infty}^{P_2} := s_1 + s_3$. Then the followings hold: } \par
{\parindent=20pt
\Item{\rm (1)} {\sl 
$\fre_c^{P_2} \in \cF_{P_2} \cap \cF_S$ and 
$\cF_{P_2} \cap \cF_S \cap \cF_{P_1} = \R_+ \cdot \fre_{\infty}^{P_2}$. }
\Item{\rm (2)} $\partial \cF_{P_2} \subset \cF_S$. 
\Item{\rm (3)} $\cF_{P_2} \cap \cF_S 
  = V\big(p_0+p_2, \, 8p_0(p_1+p_3) - (p_1-p_3)^2\big)$. \par
\Item{\rm (4)} {\sl $\fre_c^{P_2} \in \cE(\cP_{4,3}^{c0+})$ 
 for all $c \in \P_{\R}^1$. }

}
\endproclaim

\Proof 
(0) We shall show that $\fre_c^{P_2} \in \cP_{4,3}^{c0+}$ for all 
$c \in \P_{\R}^1$. 

Let $c_2(u,v) := (u-1)^2+v(u+1)$ and $c_1(u,v) := 2(u-1)v(v-u-1)$. 
$c_2(u,v) \geq 0$ for $u \geq 0$, $v \geq 0$. 
Then, 
\[\fre_c^{P_2}(0,u,v,1) 
 = v c_2(u,v) \left(c + \frac{c_1(u,v)}{2v c_2(u,v)}\right)^2 
      + \frac{(u+1)((u-1)^2+v^2)^2}{c_2(u,v)} \geq 0.\]
Thus, $\fre_c^{P_2}$ is PSD for $c \in \R$. 
$\fre_{\infty}^{P_2} = s_1+s_3$ is PSD by \Tqcbbb. 

\smallskip

(1) Since $\displaystyle 
  \fre_c^{P_2} = \lim_{h \to 0} \fre_{ch+1,h}/4h^2$, 
we have $\fre_c^{P_2} \in \cF_S$. 
Since $\fre_c^{P_2}(1,0,1,0) = 0$, we have $\fre_c^{P_2} \in \cF_S$. 
It is easy to see that $\fre_{\infty}^{P_2} \in \cF_{P_1}$ and 
$\dim (\cF_{P_2} \cap \cF_C \cap \cF_{P_1}) = 1$. 
Thus $\cF_{P_2} \cap \cF_S \cap \cF_{P_1} = \R_+ \cdot \fre_{\infty}^{P_2}$. 

\smallskip

(2) We shall determine $(\cF_{P_2} \cap \cF_C) - \cF_{P_1}$. Note that 
\[\disc_S(p_0,p_1,-p_0,p_3)
 = 2p_0 \big((p_0-p_1)^2+(p_0-p_3)^2\big)\big(8p_0(p_1+p_3) 
     - (p_1-p_3)^2\big)^3.\] 
Thus let
\[\breve{V}_C := 
   \big\{ (1,p_1,-1,p_3) \in \breve{\cH}_{4,3}^{c0} \; \big| \; 
          \hbox{$8(p_1+p_3) - (p_1-p_3)^2 = 0$}\big\}.\]
Then $\breve{V}_C 
 = \big\{(1$, $c(c-2)$, $-1$, $c(c+2))$ $\big|$ $c \in \R\big\}$. 
Each point in $\breve{V}_C$ corresponts to $\fre_c^{P_2}$. 
Since $\R_+ \cdot \breve{V}_C \cup \R_+ \cdot \fre_{\infty}^{P_2}$ 
is a conic closed convex cone, it must agree with $\cF_C$, 
and $\partial \cF_C$ is generated by $\fre_c^{P_2}$ 
($c \in \P_{\R}^1$). 

\smallskip

(3) follows from (1) and (2). 

\smallskip

(4) Put $D_{P_2} := \big\{(p_0 \colon p_1 \colon p_2 \colon p_3) 
 \in \P_{\R}^3$ 
$\big|$ $p_0+p_2 \geq 0$, $8p_0(p_1+p_3) \geqq (p_1-p_3)^2\big\}$. 
Then $\P(\cF_{P_2}) = D_{P_2}$, 
and $\fre_c^{P_2} \in \partial \cF_{P_2}$. 
Any point of $\partial D_{P_2}$ is an extremal point of $D_{P_2}$. 
\end{proof}

\bigskip

To characterize $\fre_c^{P_2}$, we need an infinitesimal local cone. 
Let $\pi \colon X \to A = \P_+^2$ be the blowing up 
at $(1 \colon 0 \colon 1)$, and put $\fre_c^b(x,y,z) 
:= \fre_c^{P_2}(x z$, $y z+1$, $z$, $1)/z^2$. 
Then $\fre_c^b(x,y,0) = 2(c x+y-t)^2$. 
This zero locus $V_X(c x+y-t$, $z)$ characterizes $\fre_c^{P_2}$. 

\smallskip

Next we shall study $\cF_S \cap \cF_C$. 
Remember that $\disc_C$ is the edge discriminant of $X_{3,3}^{c+}$ 
and $X_{3,3}^{c0+}$. Let 
\[\cD_C := \big\{(1,x,y,z) \in \breve{\cH}_{4,3}^{c0+} \; \big| \;
 \hbox{$y \geq -1$ and $d_C(x,z) \geq 0$} \big\}.\]
Then $\cD_C$ is a closed convex set such 
that $\breve{\cP}_{4,3}^{c0+} \subset \cD_C$, 
and $\big(\partial \breve{\cP}_{4,3}^{c0+}\big) \cap \Int(\cD_C)
  \subset V_{\R}(\disc_S)$ by \Tqcbg. 
We need the following polynomial to describe the cusp loci 
of $V_{\R}(\disc_S)$. 
\begin{align*}
 & f_{Q_0}(x,y) := 4(x+1)^2 + (y-3)^2, \\
 & f_{L_S}(x,y) := 2x+y-1, \\
 & f_{C^s}(x,y) := y^2+4x(y+1)-2y+13, \\
 & f_S^{cusp}(x,y) 
   = 260403739669 + 153581431744 x + 102255553008 x^2 + 5758906656 x^3 \\
 & \hskip20pt + 2375407488 x^4 - 2980119168 x^5 + 472233216 x^6 
     - 115722240 x^7 \\
 & \hskip20pt 
     + 17307648 x^8 - 438272 x^9 + 4096 x^{10} + 89440948796 y 
     + 32061417248 x y \\
 & \hskip20pt 
     + 8138124864 x^2 y - 17528885472 x^3 y - 2067065472 x^4 y 
     - 828572544 x^5 y \\
 & \hskip20pt 
     + 1188607488 x^6 y - 112318464 x^7 y - 15593472 x^8 y - 126976 x^9 y 
     + 8192 x^{10} y \\
 & \hskip20pt 
     - 223071977286 y^2 - 16231383328 x y^2 - 12833341936 x^2 y^2 
     + 40377065344 x^3 y^2 \\
 & \hskip20pt 
     + 5505244544 x^4 y^2 + 4819181440 x^5 y^2 - 264563968 x^6 y^2 
     + 218927104 x^7 y^2 \\
 & \hskip20pt 
     + 9482240 x^8 y^2 + 176128 x^9 y^2 + 4096 x^{10} y^2 + 30713189004 y^3 
     + 8960225536 x y^3 \\
 & \hskip20pt 
     + 17703049984 x^2 y^3 - 2170474624 x^3 y^3 - 7085133440 x^4 y^3
     - 4728214912 x^5 y^3 \\
 & \hskip20pt 
     - 1856392192 x^6 y^3 - 112496640 x^7 y^3 - 3928064 x^8 y^3 
     - 135168 x^9 y^3 \\
 & \hskip20pt 
     + 61229381323 y^4 - 32671427200 x y^4 - 16135419808 x^2 y^4 
     - 19363454784 x^3 y^4 \\
 & \hskip20pt 
     + 2347438208 x^4 y^4  + 668450944 x^5 y^4 + 1133005568 x^6 y^4 
     + 47364096 x^7 y^4 \\
 & \hskip20pt 
     + 1464320 x^8 y^4 - 40004520712 y^5 + 14114790976 x y^5 
     - 921252992 x^2 y^5 \\
 & \hskip20pt 
     + 9081775296 x^3 y^5 + 71177344 x^4 y^5 + 679918976 x^5 y^5 
     - 112298496 x^6 y^5 \\
 & \hskip20pt 
     - 6821888 x^7 y^5 + 10688483692 y^6 
     - 1398548800 x y^6 + 3457102112 x^2 y^6 \\
 & \hskip20pt 
     - 1135819904 x^3 y^6  + 55287936 x^4 y^6 - 134577536 x^5 y^6 
     - 18625280 x^6 y^6 \\
 & \hskip20pt 
     - 870429832 y^7  + 226903552 x y^7 - 733186304 x^2 y^7 
     - 48610432 x^3 y^7 \\
 & \hskip20pt 
     - 35363712 x^4 y^7 - 12108928 x^5 y^7 - 108565637 y^8 
     - 133149760 x y^8 \\
 & \hskip20pt + 1725104 x^2 y^8 + 6646560 x^3 y^8 - 2811392 x^4 y^8 
   + 4147404 y^9 + 9240992 x y^9 \\
 & \hskip20pt + 5649472 x^2 y^9 - 26336 x^3 y^9 + 2233722 y^{10} 
   + 1416544 x y^{10} + 84944 x^2 y^{10} \\
 & \hskip20pt + 121340 y^{11} + 16896 x y^{11} + 517 y^{12}. 
\end{align*}
Note that 
\[\disc_S(1,x,y,x) = f_{L_S}(x,y)^2 f_{C^s}(x,y)^2 
    (16x^3 - x^2y + 18x y - x^2 - y^2 + 18 x + 25 y + 26).\]

\def\Tqcbi{Lemma 4.8}%
\proclaim{Lemma 4.8} 
{\sl Consider on $\breve{\cH}_{4,3}^{c0}:(1,x,y,z) \cong \R^3$. Then 
\[\Sing\big(V_{\R}(\disc_S(1,x,y,z))\big) \cap \cP_{4,3}^{c0+}
  = \{Q_0\} \cup L^s \cup C_1^{cusp} \cup C_2^{cusp} 
     \cup C_3^{cusp} \cup C_4^{cusp} 
  \subset \breve{\cH}_{4,3}^{c0},\]
where $Q_0$, $L^s$ and $C_i^{cusp}$ are as follows:} \par
{\parindent=20pt
\Item{(1)} {\sl $Q_0$ is a point defined 
by $(x,y,z)=(-1,3,-1)$. 
This lies on $\partial \breve{\cP}_{4,3}^{c0+}$.}
\Item{(2)} {\sl $L^s$ is the half line defined by $x=z=(y-1)/2$, 
$y \geq -1$. 
Note that 
$L^s \cap \partial \cP_{4,3}^{c0+} = \{Q_0\}$. }
\Item{(3)} {\sl $C_s$ is the hyperbolic curve $C^s$ 
defined by $x=z$ and $f_{C^s}(x,y) = 0$ in $\breve{\cH}_{4,3}^{c0}$. 
Note that $C^s$ is a singular locus of of $V_{\R}(\disc_S(1,x,y,z))$, 
and $C^s \cap \partial \breve{\cP}_{4,3}^{c0+} = \{Q_0\}$. } 
\Item{(4)} {\sl Let $x = \alpha_i(y)$ be all the four real roots 
of $f_S^{cusp}(x,y) = 0$ when we regard $y$ to be a constant 
where $y \geq 3$ and $\alpha_1(y) \leq \alpha_2(y) 
 \leq  \alpha_3(y) \leq \alpha_4(y)$. 
Note that $\alpha_1(3) = \alpha_2(3) = \alpha_3(3) = \alpha_4(3) = 1$.
Then, the following four branches are cusps of $S$.} 
\begin{align*}
 & C_1^{cusp} 
   = \big\{ (1 \colon \alpha_1(y) \colon y \colon \alpha_4(y)) 
    \in \breve{\cH}_{4,3}^{c0} \; \big| \; \hbox{$y > 3$} \big\}, \\
 & C_2^{cusp} 
   = \big\{ (1 \colon \alpha_2(y) \colon y \colon \alpha_3(y)) 
    \in \breve{\cH}_{4,3}^{c0} \; \big| \; \hbox{$y > 3$} \big\}, \\
 & C_3^{cusp} 
   = \big\{ (1 \colon \alpha_3(y) \colon y \colon \alpha_2(y)) 
    \in \breve{\cH}_{4,3}^{c0} \; \big| \; \hbox{$y > 3$} \big\}, \\
 & C_4^{cusp} 
   = \big\{ (1\colon \alpha_4(y) \colon y \colon \alpha_1(y)) 
    \in \breve{\cH}_{4,3}^{c0} \; \big| \; \hbox{$y > 3$} \big\}. 
\end{align*}

}
\endproclaim

\Proof 
Let $f(x$, $y$, $z) := \disc_S(1$, $x$, $y$, $z)$ and 
$\displaystyle f_x := \frac{\partial f}{\partial x}$ and so on. 
$\Sing\big(V_{\R}(\disc_S(1$, $x$, $y$, $z))\big)$ can be obtained by solving 
the system of equations $f(x$, $y$, $z) = f_x(x$, $y$, $z) = 
 f_y(x$, $y$, $z) = f_z(x$, $y$, $z) = 0$. 
But it is next to impossible to proceed this calculation. 
Instead of it, we eliminate $z$ from $f_x(x$, $y$, $z) = 0$, 
$f_y(x$, $y$, $z) = 0$, and $f_z(x$, $y$, $z) = 0$. 
During this elimination process, we obtain some factors 
which include $\Sing\big(V_{\R}(\disc_S(1,x,y,z))\big)$. 
$2x+y-1$, $4(x+1)^2 + (y-3)^2$, $y^2+4x(y+1)-2y+13$ and $f_S^{cusp}(x$, $y)$ 
are such factors. Note that 
\begin{align*}
  \displaystyle \disc_S(1,x,y,x) 
 & = (2x+y-1)^2 (y^2+4x y + 4 x - 2 y + 13)^2 \\
 & \hskip30pt  (16x^3 - x^2y + 18x y - x^2 - y^2 + 18 x + 25 y + 26). 
\end{align*}
\end{proof}


$\Sing\big(V_{\R}(\disc_S(1,x,y,z))\big)$ may has other loci. 
But we will see that 
\[\Sing\big(V_{\R}(\disc_S(1,x,y,z))\big) \cap \partial \breve{\cP}_{4,3}^{c0+}
  \subset \{Q_0\} \cup C_1^{cusp} \cup C_2^{cusp} 
     \cup C_3^{cusp} \cup C_4^{cusp},\]
during discussion from now. 

\bigskip

{\it Proof of \Tqcbd.} 
We take the section of $\breve{\cP}_{4,3}^{c0+}$ by the hyperplane 
\[H_r := \big\{(1,x,y,z) \in \breve{\cH}_{4,3}^{c0} \; \big| \; y=r \big\}.\]
We regard $H_r$ as $(x$, $z)$-plane. Put 
\begin{align*}
 & D_r := H_r \cap \breve{\cP}_{4,3}^{c0+}
        = \big\{(x,z) \in H_r \; \big| \; 
            \hbox{$(1,x,r,z) \in \breve{\cP}_{4,3}^{c0+}$} \big\}, \\
 & D_C := \cD_C \cap H_r = \big\{(x,z) \in H_r \; \big| \; 
    \hbox{$d_C(x,z) \geq 0$} \big\}, \\
 & V_C := \partial D_C = \big\{(x,z) \in H_r \; \big| \;
      \hbox{$d_C(x,z) = 0$} \big\}, \\
 & V_S^r := \big\{(x,z) \in H_r \; \big| \;
      \hbox{$\disc_S(1,x,r,z) = 0$} \big\} 
        - (C^s \cup L^s) \cap H_r. 
\end{align*}

(O-1) If $r<-1$, then $D_r = \emptyset$, by \Tqcbg(2). 

(O-2) If $r =-1$, then the condition of (1) of \Tqcbd \ determines the 
set $\breve{\cP}_{4,3}^{c0+} \cap H_{-1}$, because of \Tqcbh. 

\bigskip
\begin{center}
\includegraphics[width=120mm,clip]{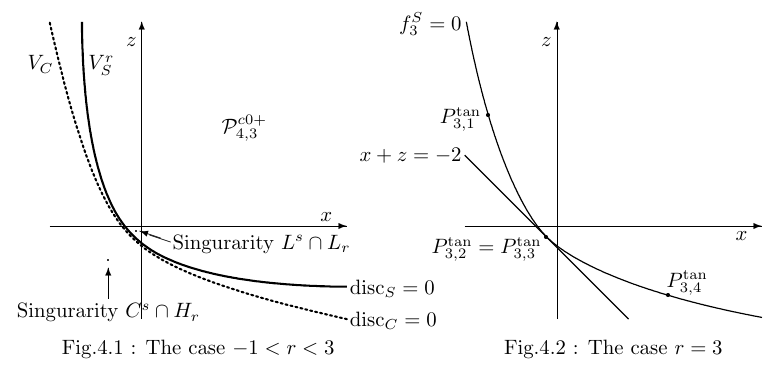}
\end{center}
\bigskip

(I) When $-1<r<3$, $V_S^r$ is as Fig.4.1. 
Two points $C^s \cap H_r$ and $L^s \cap H_r$ are 
all the isolated singularities of $V_{\R}(\disc_S) \cap H_3$. 
$V_S^r$ is a smooth curve in $D_C$ and 
enclose a convex set $\breve{\cP}_{4,3}^{c0+} \cap H_r$. 
Thus, 
\[D_r = \big\{(x,z) \in \R^2 \; \big| \; 
  \hbox{$\disc_S(1,x,r,z) \geq 0$ and $d_C(x,z) \geq 0$}\big\}.\]
Thus, the conditions of (2) of \Tqcbd \ 
determines $\breve{\cP}_{4,3}^{c0+} \cap H_r$. 

\medskip

(II) Consider the case $r=3$. 
Let 
\begin{align*}
 & f_3^S(x,z) := 
   x^6 - 4 x^5 z + 7 x^4 z^2 - 8 x^3 z^3 + 7 x^2 z^4 - 4 x z^5 + z^6 \\
 & \hskip30pt 
   - 174 x^5 - 342 x^4 z - 508 x^3 z^2 - 508 x^2 z^3 - 342 x z^4 - 174 z^5 \\
 & \hskip30pt - 414 x^4 - 712 x^3 z - 1332 x^2 z^2 - 712 x z^3 - 414 z^4 \\
 & \hskip30pt - 800 x^3 - 4320 x^2 z - 4320 x z^2 - 800 z^3  \\
 & \hskip30pt - 6592 x^2 - 16512 x z - 6592 z^2 - 16384 x - 16384 z  - 11776. 
\end{align*}
Then $\disc_S(1,x,3,z) = -2(x+z+2)^2 f_3^S(x,z)$. 
As Fig 4.2, $V_{\R}(f_3^S)$ tangents $V_C$ at three points $P_{3,1}^{\tan}$, 
$P_{3,4}^{\tan}$ and $P_{3,2}^{\tan}=P_{3,3}^{\tan}$ 
(these symbols will be explained in (III)). 
Moreover $V_{\R}(f_3^S) \subset D_C$. 
Thus, 
\[D_3 = \big\{(x,z) \in \R^2 \; \big| \; 
  \hbox{$\disc_S(1,x,3,z) \geq 0$ and $d_C(x,z) \geq 0$}\big\},\] 
and the conditions of (2) of \Tqcbd \ 
determines $\breve{\cP}_{4,3}^{c0+} \cap H_3$. 

Note that $V_{\R}(f_3^S) \cap V_{\R}(z+11.851831\cdots) = \emptyset$, 
and $V_{\R}(f_3^S) \cap V_{\R}(z-z_0)$ consists of two points 
for $z_0 > -11.851831\cdots$. 

\smallskip

(III) Consider the case $r>3$. 
Then, $V_S^r$ has just four cusps 
$P_{r,i}^{cusp} := C_i^{cusp} \cap H_r
  = \big(\alpha_i(r)$, $\alpha_{5-i}(r)\big)$ ($i=1$, $2$, $3$, $4)$. 
Since $V_S^r$ is symmetric with respect to the line $V_{\R}(x-z)$, 
it is enough to consider the part $z \geq x$. 

We observe $\cF_S \cap \cF_C \cap \breve{\cH}_{4,3}^{c0}$. 
Let 
\begin{align*}
 & L_x := \big\{ (0 \colon 0 \colon w \colon 1) \in \P_+^3 \; \big| \; 
      \hbox{$w \in [0,\infty]$} \big\}, \\
 & L_y := \big\{ (0 \colon w \colon 0 \colon 1) \in \P_+^3 \; \big| \; 
      \hbox{$w \in [0,\infty]$} \big\}, \\
 & L_z := \big\{ (0 \colon x \colon y \colon 0) \in \P_+^3 \; \big| \; 
      \hbox{$(x \colon y) \in \P_+^1$} \big\}. 
\end{align*}
Note that $\partial \Cls_{\P_+^3}(\overline{B_0}) = L_x \cup L_y \cup L_z$. 

We define a rational map $G^S : \P_+^3 \cdots\to \P(\cH_{4,3}^{c0})$ 
just before \Tqcbf. 
Since $g_0(w,0) = 0$, 
$G^S(L_y) \cap \breve{\cH}_{4,3}^{c0} = \emptyset$. 
Since $G^S(0 \colon x \colon y \colon 0) 
 = G^S(0 \colon 0 \colon x/y \colon 1)$, 
we have $G^S(L_z) = G^S(L_x)$. 
Since 
\[\disc_C\big(g_0(0,w), \, g_1(0,w), \, g_3(0,w)\big) = 0,\]
we have $G^S(L_x) \subset V_{\R}(\disc_C) \cap V_{\R}(\disc_S)$. 
Put $C_x^{\tan} := G^S(L_x)$. 

Similarly, we define a rational map 
$G' : \P_+^3 \cdots\to \P(\cH_{4,3}^{c0})$ by
\[G'(0,x,y,1) := \big(g_0(x,y):g_3(x,y):g_2(x,y):g_1(x,y)\big).\]
Let $C_z^{\tan} := G'(L_x)$. 
Then $C_x^{\tan} \cup C_z^{\tan} 
 \subset V_{\R}(\disc_C) \cap V_{\R}(\disc_S)$. 

Put $H_{\geq 3} := \big\{(1,x,r,z) 
   \in \breve{\cH}_{4,3}^{c0}$ $\big|$ $r \geq 3 \big\}$. 
We regard $H_{\geq 3} \subset 
 \breve{\cH}_{4,3}^{c0} \subset \P(\cH_{4,3}^{c0})$. 
We shall determine $C_x^{\tan} \cap H_{\geq 3}$. 
Let $\delta := 0.2955977425\cdots$ be the real 
root of $t^3+t^2+3t-1 = 0$. 
Then, all the real roots of 
$g_2(0$, $t)/g_0(0$, $t) = 3$ are $t=1$, $\delta$. 
We put \par
\begin{align*}
 & C_1^{\tan} := \left\{ G'(0,0,w,1) 
     \in \P(\cH_{4,3}^{c0}) \; \big| \; 
       \hbox{$0 < w \leq \delta$} \right\} \subset C_z^{\tan}, \\
 & C_2^{\tan} : = \left\{ G^S(0,0,w,1) 
     \in \P(\cH_{4,3}^{c0}) \; \big| \; 
       \hbox{$w \geq 1$} \right\} \subset C_x^{\tan}, \\
 & C_3^{\tan} := \left\{ G'(0,0,w,1) 
     \in \P(\cH_{4,3}^{c0}) \; \big| \; 
       \hbox{$w \geq 1$} \right\} \subset C_z^{\tan}, \\
 & C_4^{\tan} : = \left\{ G^S(0,0,w,1) 
    \in \P(\cH_{4,3}^{c0}) \; \big| \; 
       \hbox{$0 < w \leq \delta$} \right\} \subset C_x^{\tan}, \\
 & P_{r,i}^{\tan} 
     := C_i^{\tan} \cap H_r \in \breve{\cP}_{4,3}^{c0+} 
    \quad \hbox{($i=1$, $2$, $3$, $4$).} 
\end{align*}
Then $C_x^{\tan} \cap H_{\geq 3} = C_1^{\tan} \cup C_3^{\tan}$ and 
$C_z^{\tan} \cap H_{\geq 3} = C_2^{\tan} \cup C_4^{\tan}$. 
Note that 
$\cF_S \cap \cF_C \cap 
 \big\{ G^S(0,0,w,1) \in \P(\cH_{4,3}^{c0})$ 
  $\big|$ $\delta < w < 1$ $\big\} = \emptyset$. 

\proclaim{Lemma 4.9} 
$C_1^{\tan} \cup C_2^{\tan} \cup C_3^{\tan} \cup C_4^{\tan} 
  \subset \Zar(\cF_S \cap \cF_C) \cap H_{\geq 3}$. 
\endproclaim

\Proof 
Clear. 
\end{proof}

\bigskip

Put $C^{cusp} := \Cls_{\breve{\cH}_{4,3}^{c0}}(C_1^{cusp} \cup C_2^{cusp} 
   \cup C_3^{cusp} \cup C_4^{cusp})$. 
Let's determine $C_x^{\tan} \cap C^{cusp}$. 
Since $C_x^{\tan} = G^S(L_x) 
\subset V_{\R}(\disc_C) \cap V_{\R}(\disc_S)$, and 
\[G^S(0,0,w,1) 
 = \left(1 : \frac{1-2w^3}{w^2} : \frac{(w^2+1)^2 - w}{w} 
      : \frac{w^3-2}{w} \right),\]
we put $G_x^S(w) := (1-2w^3)/w^2$, 
$G_y^S(w) := ((w^2+1)^2 - w)/w$ and $G_z^S(w) :=(w^3-2)/w = G_x^S(1/w)$. 

\def\Tqcbj{Lemma 4.10}%
\proclaim{Lemma 4.10} 
{\sl $\eta(x,y) =  61 + 62 x + 56 y + 32 x^2 + 30 x y - 6 y^2
    + 9 x^3 + 4 x^2 y - 6 x y^2 - 16 y^3 
    + x^4 - 4 x^2 y^2 - 6 x y^3 + y^4 - x^3 y^2$ has 
the following properties:} \par
{\parindent=20pt
\Item{(1)} {\sl If $(1 \colon x \colon y \colon z) 
   \in C_x^{\tan} \cup C_z^{\tan}$, then $\eta(x+z$, $y)=0$. }
\Item{(2)} {\sl Let $r>3$. 
On a plane $H_r$, the zero locus $\eta(x+z,r)=0$ is 
the union of two lines. 
One is the line $P_{r,1}^{\tan}P_{r,4}^{\tan}$, 
and the other is the line $P_{r,2}^{\tan}P_{r,3}^{\tan}$. 
$\eta(x+z,r)<0$ between these two lines, and $\eta(x+z,r)>0$ outside. }

}
\endproclaim

\Proof 
(1) follows from  
$\eta\big(G_x^S(w)+G_z^S(w)$, $G_y^S(w)\big) = 0$. 

(2) $\eta(x,r) = 0$ has just two real roots for $r > 3$, 
and $\eta(y-3,y) < 0 $ for $y<3$. 
\end{proof}

\bigskip

Note that 
\[f_S^{cusp}\big(G_x^S(w), \, G_y^S(w)\big) 
 = \frac{(w-1)^4(w^2+1)^4(w^4-6w^2-8w+1)^2 f_{38}(w)}{w^{22}},\]
here $f_{38}(w)$ is a polynomial of degree $38$ 
whose real roots are two negative numbers 
$w = -8.590880\cdots$, $-2.4445756\cdots$. 
Let $\tau_1 := 0.1150\cdots$ and $\tau_2 := 2.9343\cdots$ be 
the real roots of $w^4 - 6 w^2 - 8 w + 1 = 0$, and 
\[r_1 := \frac{g_2(0,\tau_1)}{g_0(0,\tau_1)} = 7.9207039574\cdots, \quad
  r_2 := \frac{g_2(0,\tau_2)}{g_0(0,\tau_2)} = 30.474537321\cdots.\]
be the real roots of $r^4-28r^3-90r^2-92r+16353=0$. 
Then, all positive the roots of $f_S^{cusp}\big(G_x^S(w)$, $G_y^S(w)\big) = 0$ 
are $w = 1$, $\tau_1$, $\tau_2$. 
In the case $w=1$, $G^S(0,0,w,1) 
 = (1 \colon -1 \colon 3 \colon -1) = Q_0$. 
Thus, $C_x^{\tan} \cap C^{cusp}$ consists of three points $Q_0$, 
$P_{r_1,1}^{\tan}=P_{r_1,1}^{cusp}=G^S(0,0,\tau_1,1)$, 
and 
$P_{r_2,2}^{\tan}=P_{r_2,2}^{cusp} = G^S(0,0,\tau_2,1)$. 
Similarly, $C_z^{\tan} \cap C^{cusp}$ consists of three points $Q_0$, 
$P_{r_1,4}^{\tan}=P_{r_1,4}^{cusp}=G'(0,0,\tau_1,1)$, 
and $P_{r_2,3}^{\tan}=P_{r_2,3}^{cusp} = G'(0,0,\tau_2,1)$. 

\def\Tqcca{Lemma 4.11}%
\proclaim{Lemma 4.11} 
{\sl In $\breve{\cH}_{4,3}^{c0} \cong \R^3:(x,y,z)$, 
$\kappa_1 (x+z) + \kappa_2 y = 1$ defines the plane which passes through 
$P_{r_1,1}^{\tan}$, $P_{r_2,2}^{\tan}$, $P_{r_2,3}^{\tan}$ 
and $P_{r_1,4}^{\tan}$. }
\endproclaim

\Proof 
Note that 
$P_{r_1,1}^{\tan} = P_{r_1,1}^{cusp} 
= (\alpha_1(r_1)$, $r_1$, $\alpha_4(r_1))$ and so on. 
\begin{align*}
 & \alpha_2(r_1) + \alpha_3(r_1) 
    = \frac{g_1(0,\tau_1)+g_3(0,\tau_1)}{g_0(0,\tau_1)}
    = G_x^S(\tau_1) + G_z^S(\tau_1), \\
 & \alpha_1(r_2) + \alpha_4(r_2) 
    = \frac{g_1(0,\tau_2)+g_3(0,\tau_2)}{g_0(0,\tau_2)}
    = G_x^S(\tau_2) + G_z^S(\tau_2). 
\end{align*}
Solve $\kappa_1 (G_x^S(w)+G_z^S(w)) + \kappa_2 G_y^S(w) = 1$ 
for $w=\tau_1$ and $\tau_2$. 
Then, we obtain 
\begin{align*}
 \kappa_1 
  & := \frac{s^2 t^2 - t^3 + 2 t^2 - t}{s^4 - 2 s^3 t - 2 s t^3 + t^4
         - 4 s^2 t + 5 s t^2 - 2 t^3 + 2 s^2 - 2 s t - s + 1} \\
  & = 0.0129074031\cdots, \\
 \kappa_2 
  & : = \frac{- s t^2 + 2 t^2 + s - 2 t}{s^4 - 2 s^3 t - 2 s t^3 + t^4
         - 4 s^2 t + 5 s t^2 - 2 t^3 + 2 s^2 - 2 s t - s + 1} \\
  & = 0.0318925844\cdots, 
\end{align*}
where $s = \tau_1 + \tau_2$, $t = \tau_1 \tau_2$. 
Let $\gamma$, $\delta$ be all the imaginal roots $w^4 - 6 w^2 - 8 w + 1 = 0$, 
and put $s_2 := \gamma+\delta$, $t_2 := \gamma \delta$. 
Then $s+s_2=0$, $t t_2=1$, $t+t_2+s+s_2 = -6$, $t s_2+s t_2=8$. 
When we eliminate $s$, $t$, $s_1$, $t_1$ from these relations, we have 
\begin{align*}
 & 817808203 \kappa_1^6 - 546807084 \kappa_1^5 + 129155640 \kappa_1^4 
    - 13342016 \kappa_1^3 \\
 & \hskip50pt + 556080\kappa_1^2 - 10176 \kappa_1 + 64 = 0, \\
 & 43042537 \kappa_2^6 - 4514514 \kappa_2^5 - 188769 \kappa_2^4 
    - 38684 \kappa_2^3 + 4119 \kappa_2^2 - 114 \kappa_2 + 1 = 0. 
\end{align*}
\end{proof}

\bigskip

Now, we shall complete the proof of \Tqcbd. 
To prove (3), (4), (5) of \Tqcbd, we put 
\begin{align*}
 & D_r^{(3)} := \big\{ (x,z) \in H_r \; \big| \; 
   \hbox{$\kappa_1 (x+z) + \kappa_2 r \geq 1$, 
         $\disc_S(1,x,r,z) \geq 0$, $d_C(x,z) \geq 0$} \big\}, \\
 & D_r^{(4)} := 
  \left\{ (x,z) \in H_r \; \left| \; 
   \vcenter{\hbox{$\kappa_1 (x+z) + \kappa_2 r < 1$, $\eta(x+z$, $r) > 0$,} 
            \hbox{$\disc_S(1,x,r,z) \geq 0$, $d_C(x,z) \geq 0$}}
         \right.\right\}, \\
 & D_r^{(5)} := \big\{ (x,z) \in H_r \; \big| \; 
    \hbox{$\kappa_1 (x+z) + \kappa_2 r < 1$, 
    $\eta(x+z$, $r) \leq 0$, $d_C(x$, $z) \geq 0$} \big\}. 
\end{align*}

\begin{center}
\includegraphics[width=120mm,clip]{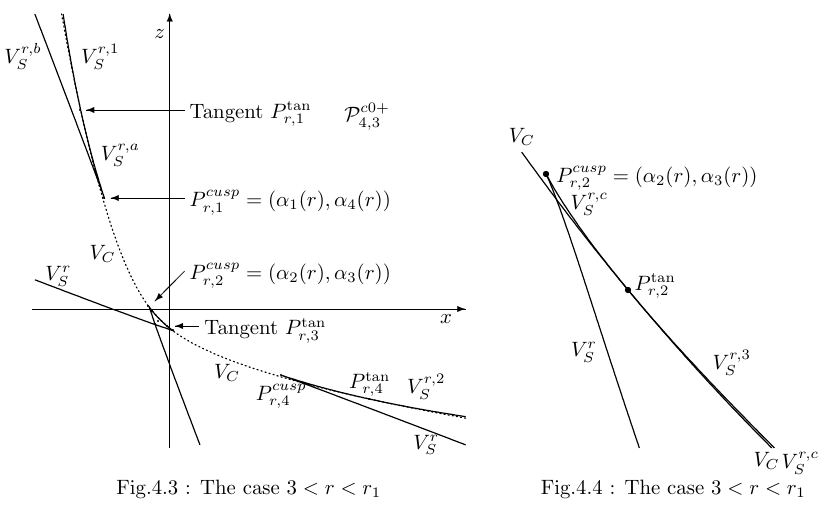}
\end{center}
\bigskip

As Fig. 4.3, we divide the part $z > \alpha_4(y)$ of $V_S^r$ 
at $P_{r,1}^{cusp}$, and denote the right part by $V_S^{r,a}$ 
and the left part by $V_S^{r,b}$. 
We mean $V_S^{r,a} \cap V_S^{r,b} = \{P_{r,1}^{cusp}\}$. 
Similarly, let $V_S^{r,c}$ be the smooth interval 
between $P_{r,2}^{cusp}$ and $P_{r,3}^{cusp}$ of $V_S^r$. 
We mean $P_{r,2}^{cusp}$, $P_{r,3}^{cusp} \in V_S^r$. 

\smallskip

(III-1) If $3<r<r_1$, then 
$P_{r,1}^{cusp} = \big(\alpha_1(r)$, $\alpha_4(r)\big) \in \Int(D_C)$, 
and $V_S^{r,a}$ tangents to $V_C$ at $P_{r,1}^{\tan}$, as Fig. 4.3. 
This implies that $P_{r,1}^{\tan} 
 \in (\partial \cF_C) \cap (\partial \cF_S)$. 
We divide the curve segment $V_S^{r,a}$ at the point $P_{r,1}^{\tan}$, 
and denote the upper part by 
\[V_S^{r,1} := \big\{(x,z) \in H_r \; \big| \; 
   \hbox{$\disc_S(x,r,z)=0$, $d_C(x,z) \geq 0$, 
   $z \geq z(P_{r,1}^{\tan})$} \big\},\] 
where $z(P)$ is the $z$-coordinate of the point $P \in H_r$. 
Then $V_S^{r,1} = \cF_S \cap V_S^{r,a}$. 
Every $P \in V_S^{r,a} - V_S^{r,1}$ is obtained as 
$P = G(0 \colon s \colon t \colon 1)$ for a certain 
$(s$, $t) \in \C^2 - \overline{B_0}$. 

Let $V_S^{r,2}$ be the symmetric set of $V_S^{r,1}$ with 
respect to the line $x=z$ on $H_r$. 

Similarly, $\big(\alpha_2(r)$, $\alpha_3(r)\big) \in \Int(D_C)$, 
and $V_S^{r,c}$ tangents to $V_C$ at $P_{r,2}^{\tan}$, as Fig. 4.4. Let 
\[V_S^{r,3} := \big\{ (x,z) \in H_r \; \big| \; 
   \hbox{$\disc_S(x,r,z)=0$, $d_C(x,z) \geq 0$, 
    $z(P_{r,3}^{\tan}) \leq z \leq z(P_{r,2}^{\tan})$} \big\}\]
be the interval of $V_S^{r,c}$ between 
$P_{r,2}^{\tan}$ and $P_{r,3}^{\tan}$. 
Then $V_S^{r,2} = \cF_S \cap V_S^{r,c}$. 
By \Tqcbj, 
\[V_S^{r,1} \cup V_S^{r,2} \cup V_S^{r,3} 
 = \left\{ (1,x,r,z) \in \partial \breve{\cP}_{4,3}^{c0+} \; \left| \; 
   \vcenter{\hbox{$\disc_S(1,x,r,z)=0$, $d_C(x,z) \geq 0$, }
            \hbox{$\eta(x+z,r) \geq 0$}} \right.\right\}.\]
So, $D_r = D_r^{(3)} \cup D_r^{(4)} \cup D_r^{(5)}$. 

\medskip

(III-2) If $r=r_1$, then 
$P_{r_1,1}^{\tan} = \big(\alpha_1(r_1)$, $\alpha_4(r_1)\big)$, 
$P_{r_1,4}^{\tan} = \big(\alpha_4(r_1)$, $\alpha_1(r_1)\big) 
  \in (\partial \cF_C) \cap (\partial \cF_S)$. 
The line defined by $\kappa_1(x+z) + \kappa_2 r_1 = 1$ agrees with 
the line $P_{r_1,1}^{\tan}P_{r_1,4}^{\tan}$. 
Others are similar as (III)-1. 

\medskip

(III-3) Consider the case $r_1 < r < r_2$. 
About $V_S^{r,3}$ the situation is same as (III-1). 

The situation of $V_S^{r,1}$ and $V_S^{r,2}$ changes. 
If $r>r_1$, then $\big(\alpha_1(r)$, $\alpha_4(r)\big) \notin D_C$. 
and $P_{r,1}^{\tan} \notin D_C$ as Fig. 4.5. 
In this case, $V_C$ and $V_S^{r,a}$ intersect at a point $Q_r^a$ 
transversally. So, $\cF_S \cap V_S^{r,a}$ agrees with the 
following $V_S^{r,4}$ in this case. 
\[V_S^{r,4} := \big\{(x,z) \in H_r \; \big| \; 
   \hbox{$\disc_S(x,r,z)=0$, $d_C(x,z) \geq 0$, 
   $z \geq z(Q_r^a)$} \big\},\] 
be the interval of $V_S^{r,a}$ upper than $Q_r^a$. 
Let $V_S^{r,5}$ be the symmetric set of $V_S^{r,4}$ with 
respect to $V_{\R}(x-z)$. 
Then, 
\begin{align*}
 & V_S^{r,4} \cup V_S^{r,5} 
 = \left\{ (1,x,r,z) \in \partial \breve{\cP}_{4,3}^{c0+} \; \left| \; 
   \vcenter{\hbox{$\disc_S(1,x,r,z)=0$, $d_C(x,z) \geq 0$, }
            \hbox{$\kappa_1(x+z)+\kappa_2 r \geq 1$}}
      \right.\right\}, \\
 & V_S^{r,3} 
 = \left\{ (1,x,r,z) \in \partial \breve{\cP}_{4,3}^{c0+} \; \left| \; 
   \vcenter{\hbox{$\disc_S(1,x,r,z)=0$, $d_C(x,z) \geq 0$, }
            \hbox{$\eta(x+z,r) \geq 0$, $\kappa_1(x+z)+\kappa_2 r < 1$}}
     \right.\right\}. 
\end{align*}
So, $D_r = D_r^{(3)} \cup D_r^{(4)} \cup D_r^{(5)}$. 

\bigskip
\begin{center}
\includegraphics[width=120mm,clip]{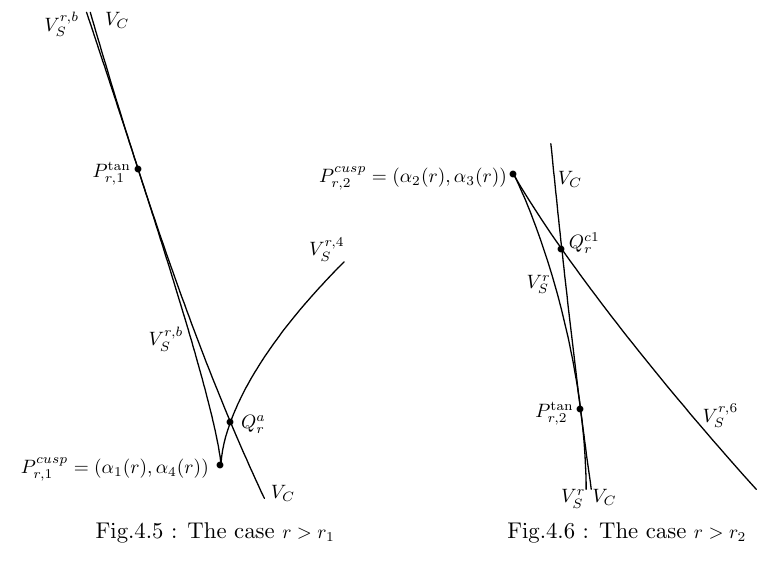}
\end{center}
\bigskip

(III-4) If $r=r_2$, then 
$P_{r_2,2}^{\tan} = \big(\alpha_2(r_2)$, $\alpha_3(r_2)\big)$,  
$P_{r_2,3}^{\tan} = \big(\alpha_2(r_3)$, $\alpha_3(r_2)\big) 
  \in (\partial \cF_C) \cap (\partial \cF_S)$. 
Others are similar as (III-3). 

\smallskip

(III-5) 
If $r>r_2$, then $\big(\alpha_2(r)$, $\alpha_3(r)\big) \notin D_C$, 
and $P_{r,2}^{\tan}$, $P_{r,3}^{\tan} \notin D_C$ as Fig.4.6. 
In this case, $V_C$ and $V_S^{r,c}$ intersect at two points $Q_r^{c1}$, 
$Q_r^{c2}$ transversally. 
So, $\cF_S \cap V_S^{r,c}$ agrees with the 
following new $V_S^{r,6}$ in this case. 
\[V_S^{r,6} := \big\{ (x,z) \in H_r \; \big| \; 
   \hbox{$\disc_S(x,r,z)=0$, $d_C(x,z) \geq 0$, 
    $z(Q_r^{c2}) \leq z \leq z(Q_r^{c3})$} \big\}\]
be the interval of $V_S^{r,c}$ between $Q_r^{c1}$ and $Q_r^{c2}$. 
Then 
\[V_S^{r,4} \cup V_S^{r,5} \cup V_S^{r,6} 
 = \big\{ (1,x,r,z) \in \partial \breve{\cP}_{4,3}^{c0+} \; \big| \; 
   \hbox{$\disc_S(1,x,r,z)=0$, $d_C(x,z) \geq 0$}\big\}.\]
If $r>r_2$, then $\kappa_1 (x+z) + \kappa_2 r \geq 1$ holds for any $(x$, $z) \in D_C$. 
Thus $D_r = D_r^{(3)}$ in this case. 

\medskip

By (III-1)---(III-5) and \Tqcca, we conclude that the conditions of 
(3), (4), (5) of \Tqcbd \ determine $\breve{\cP}_{4,3}^{c0+}$ 
when $r > 3$. \QED

Next we observe $\partial \cF_{P_1}$. 
Note that $\fre_{s,t} \in \cF_{P_1}$ when $g_0(s,t) = 0$. 

\def\Tqccbdm{Proposition 4.12}%
\proclaim{Proposition 4.12} 
{\sl Let 
\begin{align*}
 & h_{\xi}(t) := t^4-3t^3-27t^2-64t+2, \\
 & h_{\mu}(t) := t^4+t^3-2t^2-3t+1, \\
 & h_{\nu,a}(t) := t^4-7t^3+13t^2-20t+2, \\
 & h_{\nu,b}(t) := t^4-4t^3+3t^2-6t+2.
\end{align*}
Take the real roots of these polynomials as follows: 
\begin{align*}
 & V_{\R}(h_{\xi}) 
    = \big\{\xi_1 := 0.0308472031\cdots, \; 
            \xi_2 := 7.631998798\cdots\big\}, \\
 & V_{\R}(h_{\mu}) 
    = \big\{\mu_1 := 0.2882309962\cdots, \;
            \mu_4 := 1.4587325322\cdots\big\}, \\
 & V_{\R}(h_{\nu,a}) 
    = \big\{\nu_1 := 0.1070225045\cdots, \; 
            \nu_4 := 5.2319384324\cdots\big\}, \\
 & V_{\R}(h_{\nu,b}) 
    = \big\{\nu_2 := 0.3713081034\cdots, \; 
            \nu_3 := 3.586633132\cdots\big\}. 
\end{align*}
Moreover, put $\mu_2 := 1/\mu_4$ and $\mu_3 := 1/\mu_1$. 
Then the following hold: }
{\parindent=20pt
\Item{\rm (1)} {\sl $s_1+s_3+c s_2 \in \cP_{4,3}^{c0+}$, 
if and only if $0 \leq c \leq 16$. 
Moreover $s_1+s_3+16 s_2 = (1/64) \fre_{1,4}$ 
and $s_1+s_2 = \fre_{\infty}^{P_2}$. } \par
\Item{\rm (2)} {\sl $s_1+c s_2$ and $s_3+c s_2$ are PSD, 
if and only if $\xi_1 \leq c \leq \xi_2$. }
\Item{\rm (3)} {\sl There exists $\alpha_i > 0$ {\rm ($i=1$, $2$, $3$, $4$)} 
such that }
\begin{align*}
 & \fre_{\mu_1,\nu_1} = \alpha_1(s_1 + \xi_1 s_2), \quad
   \fre_{\mu_2,\nu_2} = \alpha_2(s_3 + \xi_1 s_2), \\
 & \fre_{\mu_3,\nu_3} = \alpha_3(s_3 + \xi_2 s_2), \quad
   \fre_{\mu_4,\nu_4} = \alpha_4(s_1 + \xi_2 s_2). 
\end{align*}
\Item{\rm (4)} {\sl $\cF_{P_1}$ is given as the following. 
Normalize $f \in \cF_{P_1}$ as $f = x s_1 + y s_2 + (1-x) s_3$, 
and correspond this $f$ to the point $(x$, $y) \in \R^2$. 
Let 
\begin{align*}
 & D(P_1) := \big\{ (x,y) \in \R^2 \; \big| \; 
   \hbox{$x s_1 + y s_2 + (1-x) s_3 \in \cF_{P_1}$}\big\}, \\
 & V_S^u := \big\{ (x,y) \in \R^2 \; \big| \; 
   \hbox{$0 \leq x \leq 1$, $4 \leq y \leq 8$, 
          $\disc_S(0,x,y,1-x) = 0$}\big\}, \\
 & V_S^l := \big\{ (x,y) \in \R^2 \; \big| \; 
   \hbox{$0 \leq x \leq 1$, $0 < y \leq 4$, 
          $\disc_S(0,x,y,1-x) = 0$}\big\} \cup \{(1/2,0)\}. 
\end{align*}
Then, $D(P_1)$ is a convex domain enclosed by $V_S^u$, $V_S^l$ and 
lines $x=0$, $x=1$. 
We can identify $D(P_1)$ with $\P(\cF_{P_1}) \subset \P(\cH_{4,3}^{c0})$.}

}
\endproclaim

\Proof 
(1) Let $f_t := s_1+s_3+t s_2$, $w_f(u) := u+1/u$, 
$\displaystyle v_f(t,u) := \frac{u}{2(u+1)}\big(t+2-w_f(u)\big)$, and 
$r_f(t,u) := -w_f(u)^2+2(3t+2)w_f(u)-(t-2)^2$. Then 
\[f_t(0,u,v,1) = (u+1)(v-v_f(t,u))^2 + \frac{u^2 r_f(t,u)}{4(u+1)}.\]
Note that $w_f(u) \geq 2$. 

Consider the case $w_f(u) > t+2$. 
Then $v_f(t,u)<0$ and $f_t(0,u,v,1)$ is monotonically increasing 
with respect to $v$ in $v \geq 0$. Thus 
$f_t(0,u,v,1) \geq f_t(0,u,0,1) = t u(u+1) \geq 0$. 

Consider the case $2 \leq w_f(u) \leq t+2$. 
Then $r_f(t,u) \geq r_f(t,1) = t(16-t)$. 
Thus, $f_t(0,u,v,1) \geq 0$ if $0 \leq t \leq 16$. 
If $t<0$ or $t>16$, then $v_f(t,2) > 0$ and 
$f_t(0,1,v_f(t,1),1) = r_f(t,1)/8 = t(16-t)/8 < 0$. 

Thus $f_t \in \cP_{4,3}^{c0+}$ if and only if $0 \leq t \leq 16$. 

Since $f_{16}(0,1,v,1) = 2(v-4)^2$ and $g_1(1,4)=64$, 
we have $f_{16} = \fre_{1,4}/64$. 

\smallskip

(2) Let $g_t := s_1+t s_2$, $v_g(t,u) := u(t+1-u)/2$, and 
$r_g(t,u) := -u^3+(2t+2)u^2-(t-1)^2u+4t$. Then 
\[g_t(0,u,v,1) = (v-v_g(t,u))^2 + \frac{u}{4}r_g(t,u).\]
If $u>t+1$, then $g_t(0,u,v,1) \geq g_t(0,u,0,1) = t u(u+1) \geq 0$. 

Assume that $0 \leq u \leq t+1$. 
Observe the cubic function $r_g(t,u)$. 
The roots of $(\partial/\partial u)r_g(t,u) = 0$ are 
$u_{\pm} := (2(t+1) \pm \sqrt{t^4+14t+1})/3$. 
Note that $0 \leq u_- < t+1 < u_+$. 
Thus $\min g_t(0,u,v,1) = \min r_g(t,u) = r_g(t,u_-)$. 
If $g_y \in \cE(\cP_{4,3}^{c0+})$ then $r_g(t,u_-) = 0$ and 
the cubic equation $r_g(t,u)=0$ has a double root at $u=u_-$. 
Then $\Disc_3(-1$, $2t+2$, $-(t-1)^2$, $4) = 0$. 
Note that $\Disc_3(-1$, $2t+2$, $-(t-1)^2$, $4t) = 16t \cdot h_{\xi}(t)$. 
Thus $t = \xi_1$ or $\xi_2$ if $g_t \in \cE(\cP_{4,3}^{c0+})$. 

We can also see that $g_t$ is PSD if and only if $\xi_1 \leq t \leq \xi_2$. 
Since $\disc_S(0$, $1$, $t$, $0) = -t^2h_{\xi}(t)$, 
$g_t \in \cF_S$ if and only if $t = \xi_1$ and $\xi_2$. 

Since $\disc_S(0,x,y,z)=\disc_S(0,z,y,x)$ and 
\[s_3(a,b,c,d)+t s_2(a,b,c,d) = s_1(b,a,d,c)+t s_2(b,a,d,c),\]
$s_3+t s_2$ is PSD if and only if $\xi_1 \leq t \leq \xi_2$. 

\smallskip

(3) Assume that $t=\xi_1$ or $\xi_2$. 
Then $r_g(\xi_i,u_0)=0$ for $\exists u_0 \in \R$ and $h_{\xi}(\xi_i)=0$. 
Eliminate $t$ from $r_g(t,u) = 0$ and $h_{\xi}(t) = 0$, we obtain 
\[h_{\mu}(u)^2 (u^4-16u^3+48u^2-384u+512) = 0.\]
$u = u_0$ must be a multiple root of the above equation. 
Thus $h_{\mu}(u_0) = 0$, and $u_0 = \mu_1$ or $\mu_4$. 
Let $\nu_1 := v_g(\xi_1, \mu_1)$ and $\nu_4 := v_g(\xi_2, \mu_4)$. 
Then $g_{\xi_1}(0,\mu_1,\nu_1,1) = 0$ and $g_{\xi_2}(0,\mu_4,\nu_4,1) = 0$. 
Thus $g_{\xi_1} \in \R_+ \cdot \fre_{\mu_1,\nu_1}$ 
and $g_{\xi_2} \in \R_+ \cdot \fre_{\mu_4,\nu_4}$ by \Tqcbbd (5). 
Eliminate $t$ and $u$ from $v = v_g(t, u)$, $h_{\xi}(t) = 0$ 
and $h_{\mu}(u) = 0$, we obtain $h_{\nu,a}(\nu_1) = h_{\nu,a}(\nu_4) = 0$. 

Let $h_t := s_1+t s_2$. Then $h_t(0,u,v,1) = u^3 g_t(0,1/u,v/u,1)$. 
$\mu_2=1/\mu_4$ and $\mu_3=1/\mu_2$ are roots of $u^4 h_{\mu}(1/u)$. 
Let $\nu_3 := \nu_1/\mu_1$ and $\nu_2 := \nu_4/\mu_4$. 
Then $h_{\xi_2}(0,\mu_2,\nu_2,1) = 0$ and $h_{\xi_1}(0,\mu_4,\nu_4,1) = 0$. 
Thus $h_{\xi_2} \in \R_+ \cdot \fre_{\mu_2,\nu_2}$ 
and $h_{\xi_2} \in \R_+ \cdot \fre_{\mu_3,\nu_3}$. 
Eliminate $u$ from $v = u v_g(t, 1/u)$, $h_{\xi}(t) = 0$ 
and $h_{\mu}(1/u) = 0$, we obtain $h_{\nu,b}(\nu_2) = h_{\nu,b}(\nu_3) = 0$. 

\smallskip

(4) For $f = p_0 s_0 + p_1 s_1 + p_2 s_2 + p_3 s_3 \in \cH_{4,3}^{c0}$, 
$\disc(P_1) = p_0$ and $\disc(P_2) = p_0+p_2$. 
By \Tqcbf, $\partial \cF_{P_1} \subset 
\cF_{P_2} \cup \cF_S \cup \cF_C.$
$\disc(P_2)=0$ corresponds to $y=0$. 
Thus, $D(P_1)$ must be included in the upper half space $y \geq 0$. 
Since $\disc_C(0$, $x$, $(1-x)) = -x^2(1-x)^2$ and $(1/2$, $1) \in D(P_1)$, 
$D(P_1)$ is included in the stlipe $0 \leq x \leq 1$. 
$V_S^u$ and $V_S^l$ are curves as Fig 4.7. 
Thus, we have the conclusion. 
\end{proof}

\bigskip
\begin{center}
\includegraphics[width=120mm,clip]{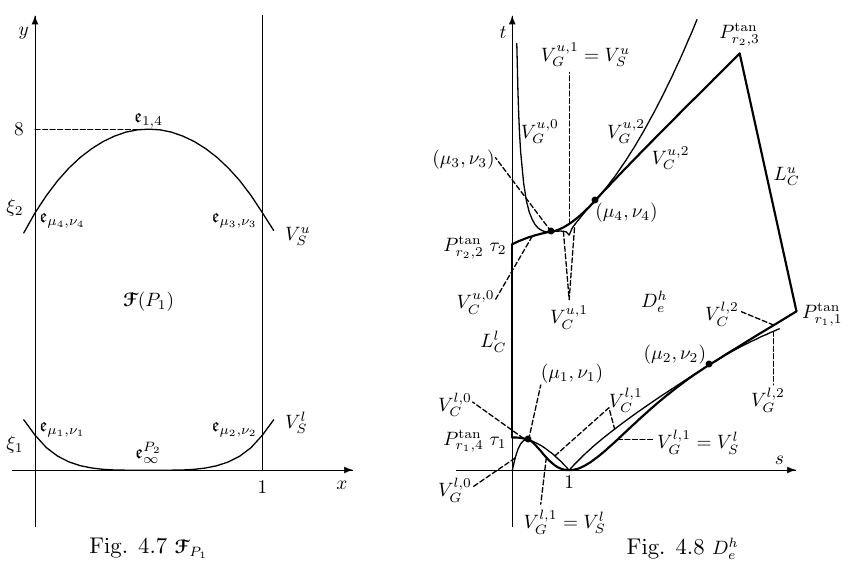}
\end{center}
\bigskip

Let 
\begin{align*}
 & D_e^h 
   := \big\{ (u \colon v \colon w) \in \P_+^2 \; \big| \; 
     \hbox{$\fre_{u,v,w}^h \in \cE(\cP_{4,3}^{c0+})$}\big\} 
   = \big\{ (u \colon v \colon w) \in \P_+^2 \; \big| \; 
     \hbox{$\fre_{u,v,w}^h \in \cP_{4,3}^{c0+}$}\big\}, \\
 & d_e^{Ch}(u,v,w) := 
  \frac{\disc_C\left(g_0^h(u,v,v), \, g_1^h(u,v,w), \, g_3^h(u,v,w)\right)}{
    u^2 w^2 (u+w-v)^2 ((u-w)^2+v^2)^2}, \\
 & d_e^C(s,t) := d_e^{Ch}(s,t,1). 
\end{align*}
$d_e^{Ch}(u,v,w)$ is a homogeneous polynomial of degree 10. 
Let $L_w := V_+(w) \subset \P_+^2$ be the line segment at infinity. 
For $(u \colon v \colon w) \in \P_+^2 - L_w$, 
let $s := u/w$, $t := v/w$ and regard $\P_+^2 - L_w$ to be the 
the first quadrant of the $(s$, $t)$-plane $\R_+^2$. 
The point $(s,t)=(1,0) \notin D_e^h$ because $\fre_{1,0}=0$. 
For completion of $D_G^h$, it is better to put $\fre_{\infty}^{P_2}
 = s_1+s_3$ at $(s$, $t)=(1$, $0)$. 
In the quadrant $s \geq 0$ and $t \geq 0$, 
the curve $V_C := V_{\R}(d_e^C(s,t))$ has two connected 
components $V_C^l$ and $V_C^u$. 
Similarly, $V_G := V_{\R}(g_0(s,t))$ has two connected 
components $V_G^l$ and $V_G^u$. 
$V_C^l$ and $V_G^l$ are included in $t<s+1$, 
and $V_C^u$, $V_G^u$ are included in $t>s+1$. 

$V_C^l \cap V_G^l = \big\{(\mu_1,\nu_1)$, $(\mu_2,\nu_2)\big\}$, and 
$V_C^u \cap V_G^u = \big\{(\mu_3,\nu_3)$, $(\mu_4,\nu_4)\big\}$. 
Divide $V_C^l$ and $V_G^l$ by the points $(\mu_1,\nu_1)$ and $(\mu_2,\nu_2)$, 
and define $V_C^{l,i}$ and $V_G^{l,i}$ ($i=0$, $1$, $2$) as Fig. 4.8. 
Similarly, we divide $V_C^u$ and $V_G^u$ by the points 
$(\mu_3,\nu_3)$ and $(\mu_4,\nu_4)$,  
and we define $V_C^{u,i}$ and $V_G^{u,i}$ ($i=0$, $1$, $2$) as Fig. 4.8. 
The segment $V_G^{l,1}$ corresponds to $V_S^l$, 
and $V_G^{u,1}$ corresponds to $V_S^u$. 

\def\Tqccbd{Theorem 4.13}%
\proclaim{Theorem 4.13} 
\[D_e^h = \left\{\ (u \colon v \colon w) \in \P_+^2 \  \left| \  
  \vcenter{\hbox{$g_0^h(u,v,w) \geq 0$, $v>0$ and }
           \hbox{one of the following (1) or (2) holds.}}
\right.\right\}.\]
{\parindent=20pt
\Item{\rm (1)} $d_e^{Ch}(u,v,w) \geq 0$. 
\Item{\rm (2)} {\sl $g_1^h(u,v,w) \geq 0$ and $g_3^h(u,v,w) \geq 0$.}

}
\endproclaim

\Proof 
We already proved that if $\fre_{u,v,w}^h$ is PSD 
then $\fre_{u,v,w}^h \in \cE(\cP_{4,3}^{c0+})$. 
By \Tqcbbd, $g_0^h(u,v,w) \geq 0$ is required. 

\smallskip

(i) Consider the case  $g_0^h(u,v,w) > 0$. 

Let $p_i = g_i^h(u,v,w)/g_0^h(u,v,w)$ ($i=1$, $2$, $3$). 
$\fre_{u,v,w}^h$ is PSD, if and only if 
$(p_1$, $p_2$, $p_3)$ satisfy the condition of \Tqcbd. 
If $\fre_{u,v,w}^h$ is PSD, 
then $\fre_{u,v,w}^h \in \cF_S$ and $\disc_S(1,p_1,p_2,p_3) = 0$. 
Conditions about $\eta(p_1+p_3, p_2)$ and 
$\kappa_1 (p_1 + p_3) + \kappa_2 p_2 - 1$ do 
not have special sence in this case. 
Thus, $\fre_{u,v,w}^h$ is PSD, 
if and only if $d_C(p_1,p_3) \geq 0$. 
That is, $\disc_C(1,p_1,p_3) \geq 0$ or `$p_1 \geq 0$ and $p_3 \geq 0$'. 
$\disc_C(1,p_1,p_3) \geq 0$ is equivalent to $u \geq 0$, $w \geq 0$ 
and $d_e^{Ch}(u,v,w) \geq 0$. 
Thus, we have the conclusion. 

\smallskip

(ii) Consider the case $g_0^h(u,v,w) = 0$. 

In this case, $V_S^l$ and $V_S^u$ of \Tqccbdm \ appears in $\partial D_e^h$. 
By \Tqccbdm, $V_G^l \cup V_G^u$ is determined by $g_0(s,t)=0$, 
$g_1(s,t) \geq 0$ and $g_3(s,t) \geq 0$. 
\end{proof}

\bigskip

By the avobe theorem, $\fre_{0,t}$ is PSD, 
if and only if $\tau_1 \leq t \leq \tau_2$. 
Similarly, 
$\fre_{t,1,0}^h = t \fre_{0,t,1}^h - (t^2-1)(t^2+1)^2 s_2$ is PSD, 
if and only if $1/\tau_2 \leq t \leq 1/\tau_1$. 

We shall observe $\partial D_e^h$ precisely. 
$\cF_S \cap \cF_{P_1}$ and $\cF_S \cap \cF_{P_2}$ 
are determined already. 
We observe the part of $\partial D_e^h$ corresponding 
to $\cF_S \cap \cF_C$. 

Let $L_C^l$ be the line segment defined 
by $s=0$ and $\tau_1 \leq t \leq \tau_2$, 
and put $V_{SC}^1 := V_C^{l,0} \cup L_C^l \cup V_C^{u,0}$. 
Since $V_{SC}^1 \subset V(\disc_C) \cap \partial D_e^h$, 
if $(s,t) \in V_{SC}^1$, there exists $\rho \in \P_{\R}^1$ such 
that $\fre_{s,t}(0,0,\rho,1) = 0$. 
We denote this $\rho$ by $\rho(s,t)=\rho(s \colon t \colon 1)$. 
Note that $\rho(0,\tau_1)=\tau_1$, $\rho(\mu_1$, $\nu_1)=0$. 
If $(s,t) \in V_C^{l,0}$, $\rho(s,t)$ is monotonically decreasing 
from $\tau_1$ to $0$ with respect to $s$. 
Similarly, $\rho(0,\tau_2)=\tau_2$, $\rho(\mu_3$, $\nu_3)=+\infty$, 
and of $(s,t) \in V_C^{u,0}$, $\rho(s,t)$ is monotonically increasing 
from $\tau_2$ to $+\infty$ with respect to $s$. 
If $(s,t) \in L_C^l$, then $\rho(s,t)=t$. 
So, each $\rho \in [0$, $+\infty]$, 
there exists unique $(s$, $t) \in V_{SC}^1$ such that $\rho(s,t) = \rho$. 
That is, $\fre_{s,t}(0,0,u,1) = 0$. 
Note that $(s,t)=(0,\tau_1)$ corresponds 
to $P_{r_1,4}^{\tan} = P_{r_1,4}^{cusp}$, 
and $(s,t)=(0,\tau_2)$ corresponds 
to $P_{r_2,2}^{\tan} = P_{r_2,2}^{cusp}$. 

When $w=0$, let $L_C^u$ be the interval of $L_w = V_+(w)$ between 
$(1 \colon \tau_1 \colon 0)$ and $(1 \colon \tau_2 \colon 0)$. 
Note that $V_C^{l,2} \cap L_w = (1 \colon \tau_1 \colon 0)$ and 
$V_C^{u,2} \cap L_w = (1 \colon \tau_2 \colon 0)$. 
Put $V_{SC}^2 := V_C^{u,2} \cup L_C^u \cup V_C^{l,2}$. 
Note that $\rho(\mu_4,\nu_4)=0$, 
$\rho(1 \colon t \colon 0) = 1/t$, and $\rho(\mu_3,\nu_3)=+\infty$. 
So, each $\rho \in [0$, $+\infty]$, 
there exists unique $(u \colon v \colon w) \in V_{SC}^2$ such 
that $\rho(u \colon v \colon w) = \rho$. 

$L_C^l$ corresponts to $C_x^{\tan}$, 
and $L_C^r$ corresponts to $C_z^{\tan}$. 
$P_{r,1}^{\tan}$ moves on the interval of $L_C^u$ defined by 
$1/\tau_2 \leq v/u \leq 1$. $Q_r^a$ moves on $V_C^{l,2}$. 
$P_{r,2}^{\tan}$ moves on the interval of $L_C^l$ defined by 
$1 \leq t \leq \tau_2$. $Q_r^{c1}$ moves on $V_C^{u,0}$. 

If $(s,t) \in V_C^{l,0} \cup V_C^{l,2} \cup V_C^{u,0} \cup V_C^{2,2}$ 
and $\rho = \rho(s,t)$, then $s$ and $\rho$ satisfy the following relation: 
\begin{align*}
  & (\rho^3+1)^2(\rho^4-8\rho^3-6\rho^2+1)s^4 \\
  & + (3\rho+1)(-\rho^9-3\rho^8-2\rho^7-6\rho^6-14\rho^5
                +6\rho^4-2\rho^3-6\rho^2-5\rho+1) s^3 \\
  & - 2(\rho^{10}+12\rho^8+26\rho^7-\rho^6+4\rho^5
                -\rho^4+26\rho^3+12\rho^2+1) s^2 \\
  & + (\rho+3)(\rho^9-5\rho^8-6\rho^7-2\rho^6+6\rho^5
                -14\rho^4-6\rho^3-2\rho^2-3\rho-1) s \\
  & + (\rho^3+1)^2(\rho^4-6\rho^2-8\rho+1) = 0. 
\end{align*}

Especially, we have the following: 

\proclaim{Proposition 4.14} 
{\sl For $t \in [0$, $+\infty]$. 
let $\cL_t^C \subset \cF_C$ be 
the local cone of $\cP_{4,3}^{c0+}$ 
at $(0 \colon 0 \colon t \colon 1) \in \P_+^3$. 
Take $(u_i \colon v_i \colon w_i) \in V_{SC}^i$ such that 
$\rho(u_i \colon v_i \colon w_i) = 1$ {\rm ($i=1$, $2$)}. Then }
\[\cL_t^C = \R_+ \cdot \fre_{u_1,v_1,w_1}^h 
  + \R_+ \cdot \fre_{u_2,v_2,w_2}^h.\]
\endproclaim

\def\Tqccb{Theorem 4.15}%
\proclaim{Theorem 4.15} 
{\sl All the elements of $\cE(\cP_{4,3}^{c0+})$ is 
the positive multiple of $\fre_{u,v,w}^h$ 
{\rm ($(u \colon v \colon w) \allowbreak \in D_e^h$)} or 
$\fre_t^{P_2}$ {\rm ($t \in \P_{\R}^1$)}. }
\endproclaim

{\it Proof of \Tqai.} 
Let $e_1$,$\ldots$, $e_{20}$ be all the monomials in $\cH_{4,3}$. 
Assume that $(s \colon t \colon 1) \in D_e^h$, $s>0$, $t>0$ and $t \ne s+1$. 
Put $u := \sqrt{s}$, $v := \sqrt{t}$ and 
$E_{s,t}(a,b,c,d) := \fre_{s,t}(a^2$, $b^2$, $c^2$, $d^2)$. 
$V_{\R}(E_{s,t})$ contains at least 27 isolated points. 
Among $V_{\R}(E_{s,t})$, we choose the following 20 points: 
${\bf a}_1 = (1 \colon 1 \colon 1 \colon 1)$, 
${\bf a}_2 = (-1 \colon 1 \colon 1 \colon 1)$, 
${\bf a}_3 = (1 \colon -1 \colon 1 \colon 1)$, 
${\bf a}_4 = (1 \colon 1 \colon -1 \colon 1)$, 
${\bf a}_5 = (1 \colon 1 \colon 1 \colon -1)$, 
${\bf a}_6 = (1 \colon 1 \colon -1 \colon -1)$, 
${\bf a}_7 = (1 \colon -1 \colon 1 \colon -1)$, 
${\bf a}_8 = (1 \colon -1 \colon -1 \colon 1)$, 
${\bf a}_9 = (0 \colon u \colon v \colon 1)$, 
${\bf a}_{10} = (1 \colon 0 \colon u \colon v)$, 
${\bf a}_{11} = (v \colon 1 \colon 0 \colon u)$, 
${\bf a}_{12} = (u \colon v \colon 1 \colon 0)$, 
${\bf a}_{13} = (0 \colon u \colon v \colon -1)$, 
${\bf a}_{14} = (-1 \colon 0 \colon u \colon v)$, 
${\bf a}_{15} = (v \colon -1 \colon 0 \colon u)$, 
${\bf a}_{16} = (u \colon v \colon -1 \colon 0)$, 
${\bf a}_{17} = (0 \colon u \colon -v \colon 1)$, 
${\bf a}_{18} = (1 \colon 0 \colon u \colon -v)$, 
${\bf a}_{19} = (-v \colon 1 \colon 0 \colon u)$, 
${\bf a}_{20} = (u \colon -v \colon 1 \colon 0)$. 
Let $a_{i,j} := e_j({\bf a}_i)$ and $A := (a_{i,j})$. Then 
\[\det A = \pm 1048576 s^2 t^2 (t-s-1)^4 ((s-1)^2+t^2)^4 \ne 0.\]
Thus, there exists no $g \in \cH_{4,3} - \{0\}$ such that 
$g({\bf a}_i) = 0$ for all $1 \leq i \leq 20$. 
Thus $E_{s,t} \notin \Sigma_{4,6}$.  \QED

\bigskip

It seems that if $(s,t) \in V_C^{l,0} \cup V_G^{l,1} \cup V_C^{l,2} 
 \cup V_C^{u,0} \cup V_G^{u,1} \cup V_C^{u,2} - (L_C^l \cup L_C^u)$, 
then $\fre_{s,t} \in \cE(\cP_{4,3}^+)$. 
If $(s,t) \in \Int(D_e^h) \cup L_C^l \cup L_C^u$, 
then $\fre_{s,t} \notin \cE(\cP_{4,3}^+)$. 
This suggests that $\cE(\cP_{4,3}^+)$ is not so simple. 

If $(s,t) \in V_C^{l,0} \cup V_C^{l,2} \cup V_C^{u,0} \cup V_C^{u,2}
  - (L_C^l \cup L_C^u) 
  - \big\{(\mu_i,\nu_i)$ $\big|$ $i=1$, $2$, $3$, $4\big\}$, 
then $\fre_{s,t}(a^2,b^2,c^2,d^2)$ has 35 isolated zeros, 
because $\fre_{s,t}(0,0,r,1) = 0$ by $r=\rho(s,t)>0$, $t \ne 1$. 
So, $\fre_{s,t}(a^2,b^2,c^2,d^2)$ will be an extremal 
element of $\cP_{4,6}$ which is irreducible. 

\bigbreak
{\bf 4.2. Structure of $\cP_{4,3}^{c+}$}%
\par\penalty1000\vskip0.4em plus0.1em minus0.1em
We have not complete any of (I1), (I2), (I3) for $\cP_{4,3}^{c+}$. 
But, we shall give (I4) and some information about $X_{4,3}^{c+}$. 

We choose $s_0 := S_3 - S_{1,1,1}$, $s_1 := S_{2,1,0} - S_{1,1,1}$, 
$s_2 := S_{2,0,1} - S_{1,1,1}$, $s_3 := S_{1,2,0} - S_{1,1,1}$, 
$s_4 := S_{1,1,1}$ as a basis of $\displaystyle \cH_{4,3}^c$, 
and define $\Phi_{4,3}^c : 
  \P_+^3 \lto \P_+^4$ by $\Phi_{4,3}^c(a) 
 = \big(s_0(a):s_1(a):s_2(a):s_3(a):s_4(a)\big)$. 
Put $X_{4,3}^{c+} := \Phi_{4,3}^c(\P_+^3)$. 
$\Psi_{4,3}^{c0} : \P_+^3/(\Z/4\Z) \cdots\to X_{4,3}^{c0+}$ split as 
\[\Psi_{4,3}^{c0} : \P_+^3/(\Z/4\Z) \mapr{\Psi_{4,3}^c} X_{4,3}^{c+} 
   \mapr{\rm pr} X_{4,3}^{c0+}.\]

\def\Tqccd{Proposition 4.16}%
\proclaim{Proposition 4.16} 
{\sl Let 
\begin{align*}
 & f_{4,3}^c(x_0,x_1,x_2,x_3,x_4) \\
 & := x_1^3 - x_0 x_1 x_3 + x_3^3 
      + x_1^2 x_2 + x_1 x_2^2 + x_2^2 x_3 + x_2 x_3^2 
      - x_0 x_1 x_2 - x_0 x_2 x_3 - x_1 x_2 x_3 \\
 & \hskip12pt 
    + x_4\Big(x_0^2 + 5 x_1^2 + x_2^2 + 5 x_3^2 
     - 2 x_0 x_1 - 2 x_0 x_2 - 2 x_0 x_3 
      + 2 x_1 x_2 - 6 x_1 x_3 + 2 x_2 x_3 \Big) 
\end{align*}
Then $\displaystyle \Zar(X_{4,3}^{c+})
  = \big\{ {\bf x} \in \P_{\R}^4$ $\big|$ 
$f_{4,3}^c({\bf x})=0$, $f_{4,3}^{c0}({\bf x}) \geq 0 \big\}$ 
with the coordinate system $x_i = s_i(a_0 \colon \cdots \colon a_3)$ 
($i=0$,$\ldots$, $4$). 
This cubic hypersurface $V_{\R}(f_{4,3}^c)$ has an isolated singularity 
at $\Phi_{4,3}^c(1 \colon 1 \colon 1 \colon 1) 
  = (0 \colon 0 \colon 0 \colon 0 \colon 1)$. }
\endproclaim

\Proof 
Using PC, we have $f_{4,3}^c(s_0$, $s_1$, $s_2$, $s_3$, $s_4) = 0$. 
Define ${\rm pr} : X_{4,3}^{c+} \cdots\to X_{4,3}^{c0+}$ by 
${\rm pr}(x_0 \colon \cdots \colon x_4) = (x_0 \colon \cdots \colon x_3)$. 
This is a birational map. 
By \Tqcbb, we have the conclusion. 
\end{proof}

\proclaim{Proposition 4.17} 
{\sl $X_{4,3}^{c+}$ does not have the main component. }
\endproclaim

\Proof 
Assume that $X_{4,3}^{c+}$ has the main component. 
Note that $\cE(\cP_{4,3}^{c0+}) 
   = \cE(\cP_{4,3}^{c+}) \cap \cH_{4,3}^{c0+}$. 
Let $f$ be an element of the main component 
such that $f \in \cE(\cP_{4,3}^{c+}) 
   - \cE(\cP_{4,3}^{c0+})$. 
Then, there exists 
${\bf a} = (a \colon b \colon c \colon 1) \in \Int(\P_+^3)$ such 
that $f({\bf a}) = 0$. 
$(a$, $b$, $c) \ne (1$, $1$, $1)$, since $f \not\in \cP_{4,3}^{c0+}$. 
Put ${\bf b} := (b \colon c \colon 1 \colon a) \in \Int(\P_+^3)$. 
Note that ${\bf a} \ne {\bf b}$. 
Then the line ${\bf a}{\bf b}$ is a bitangent line of the cubic 
surface $V_{\C}(f) \in \P_{\C}^3$. 
But a cubic surface has no bitangent line. A contradiction. 
\end{proof}

\bigskip

{\it Proof of \Tqah.}  
Let $\overline{B_0} 
  := \big\{ (0 \colon s \colon t \colon 1) \in \P_+^3$ $\big|$  
       $\hbox{$s$, $t \in \R_+$} \big\}$, 
and $\Omega := \{(1 \colon 1 \colon 1 \colon 1)\} \cup \overline{B_0}$. 
By \Tqbby, it is enough to show 
  $\cE(X_{4,3}^{c+}) \subset \Phi_{4,3}^c(\Omega)$. 
Take any ${\bf x} \in \cE(X_{4,3}^{c+})$. 
Then, there exists $D \in \Delta(X_{4,3}^{c+})$ such that ${\bf x} \in D$ 
and that $\cF_D$ is a face component. 
By the above proposition, 
$D \subset \partial X_{4,3}^{c+} \cup \Sing(X_{4,3}^{c+})$. 
If ${\bf x} \in \partial X_{4,3}^{c+}$, 
then ${\bf x} \in \Phi_{4,3}^c(\overline{B_0})$. 
If ${\bf x} \in \Sing(X_{4,3}^{c+})$, 
then ${\bf x} = \Phi_{4,3}^c(1 \colon 1 \colon 1 \colon 1)$ 
by \Tqccd.  \QED

\bigskip

For test set, we can prove the following by the same idea. 

\proclaim{Proposition 4.18} 
{\sl Assume that $f(x_1$,$\ldots$, $x_n) \in \cH_{n,3}$, 
and there exists ${\bf a} \in \Int(\P_+^{n-1})$ such that 
$f({\bf a}) = 0$ and $\displaystyle 
\frac{\partial}{\partial x_i}f({\bf a}) = 0$ for all $i=1$,$\ldots$, $n$. 
Then $f \in \cP_{n,3}^+$ if and only if $f({\bf b}) \geq 0$ 
for all ${\bf b} \in \partial \P_+^{n-1}$. }
\endproclaim

\Proof 
Assume that $f({\bf c}) < 0$ for a certain ${\bf c} \in \Int(\P_+^{n-1})$. 
We may assume that $f$ take a minimal value at ${\bf c}$. 
Put $g(t) := f((1-t){\bf a}+t{\bf c})$. 
Then, a cubic polynomial $g(t)$ takes minimal values at $t=0$ and $t=1$. 
A contradiction. 
\end{proof}

\removelastskip\penalty-400\vskip2.5em plus0.3em minus0.3em
{\bf Section 5. Philosophy of Semialgebraic Variety.} 
\par\penalty1000\vskip0.8em plus0.2em minus0.2em
{\bf 5.1. Real algebraic quasi-variety.}
\par\penalty1000\vskip0.4em plus0.1em minus0.1em
Till \S 4, we used the notion of (quasi-) semialgebraic varieties 
without exact definition. 
In this section, we shall discuss how its definition should be, 
at least for theory of PDS cones. 
Before to give it, we must discuss what a real algebraic variety is. 

Usually, we say $(X$, $\cO_X)$ is an {\it algebraic variety over $\R$} 
when $(X$, $\cO_X)$ is 
an integral separated scheme of finite type over $\R$. 
$X(\R)$ denotes the set of $\R$-rational points, 
and $X_{\C} := X \times_{\Spec \R} \Spec \C$. 
By this definition, $X$ and $X_{\C}$ are irreducible and reduced. 
To treat possibly reducible or non-reduced varieties, 
we shall call a separated scheme of finite type over $\Spec(\R)$ to 
be an {\it algebraic quasi-variety over $\R$}. 
This notion is not convenient for algebraic inequalities. 
For example, there exists infinitely many algebraic varieties $X$ over $\R$ 
such that $X(\R) = \R^2$. $X$ may not be affine even if $X(\R)=\R^2$. 

The definition of a real algebraic variety is given 
in \cite[\S 3.2]{RefBCR}. 
According to this definition, every real algebraic variety is reduced 
but may be reducible (i.e. not irreducible). 
To keep consistency with complex algebraic geometry, 
we shall add a restriction that 
real algebraic varieties must be irreducible and separated. 
To treat possibly non-reduced varieties, 
we shall give alternative definition of 
real algebraic quasi-varieties as the following: 

\def\Tpga{Definition 5.1} %
\proclaim{Definition 5.1}{\rm (Real algebraic quasi-variety)} 
{\rm (I) A locally ringed space $(X$, $\cR_X)$ is called 
a {\it real algebraic quasi-variety}, if there exists 
a separated scheme $(Y$, $\cO_Y)$ of 
finite type over $\Spec \R$ which satisfies the followings: \par
{\parindent=20pt
\Item{(1)} There exists an injective morphism $\iota \colon (X$, $\cR_X) 
 \lto (Y$, $\cO_Y)$ as locally ringed spaces, 
and $\iota$ induces a homeomorphism $X \to Y(\R)$ as topological spaces 
with respect to the Zariski topology and the Euclidean topology. 
\Item{(2)} Take any affine open subset $V \subset Y$. 
Let $\frn_P$ be the maximal ideal of $\cO_Y(V)$ corresponding 
to a closed point $P \in Y$. 
For an arbitral non-empty subset $U \subset V \cap \iota(X)$, we put 
\[S_U := \bigcap_{P \in U} \big(\cO_Y(V) - \frn_P\big).\]
If $U$ is an Euclidean open set, 
then $\iota^* : S_U^{-1} \cO_Y(V) \lto \cR_X(\iota^{-1}(U))$ is 
an isomorphism of $\R$-algebra. 
Thus, each maximal ideal $\frm \subset \cR_X(\iota^{-1}(V))$ 
corresponds to a point $P \in \iota^{-1}(V) \subset X$. 
\Item{(3)} Take an arbitral affine open subset $V \subset Y$. Then 
\[\big\{ f \in \cO_Y(V) \; \big| \; 
  \hbox{$f(P)=0$ for all $P \in V(\R)$}\big\}\]
is a nilpotent ideal of $\cO_Y(V)$. 

}
In this case, $Y$ is said to be a {\it $\R$-scheme which represents} $X$. 
If we can choose $Y$ such that $Y_{\C}$ is irreducible and reduced, 
then we shall call $X$ to be a {\it real algebraic variety} 
(See \cite[Notation 0.1]{RefKa}). 

$U \subset X$ is called an {\it affine open subset} of $X$, 
if there exists an affine open subset $U_Y \subset Y$ such 
that $U = \iota^{-1}(U_Y(\R))$. 
Zariski open (resp. closed) subsets are defied similarly. 
The {\it Euclidean topology} of $X$ is the topology 
induced from the analytic topology of $Y_{\C}$. 
$Y(\R)$ is also denoted as $Y_{\C}(\R)$. 
When $V \subset Y$ is an affine open subset and $B \subset V(\R)$ is 
a subset such that $\Cls_{Y(\R)}(\Int(B)) = \Cls_{Y(\R)}(B)$, we put 
\[S_B := \bigcap_{P \in B} \big(\cO_Y(V) - \frn_P\big),\]
and $\cR_X(\iota^{-1}(B)) := \iota^*\big(S_B^{-1} \cO_Y(V)\big)$. 
By this definition, $(X$, $\cR_X)$ can be also regarded 
as a locally ringed space with respect to the Zariski topology 
and the Euclidean topology. 
We usually omit to write $\iota$. For example, we write $X=Y(\R)$. }
\endproclaim

Note that if $(X$, $\cR_X)$ is a (possibly reducible) 
separated real algebraic variety in the sense of \cite{RefBCR}, 
there exists a reduced scheme $(Y$, $\cO_Y)$ which satisfies 
the above conditions. 
Contrary, if $(X$, $\cR_X)$ is a reduced real algebraic 
quasi-variety as \Tpga, then $(X$, $\cR_X)$ is a real algebraic variety 
in the sense of \cite{RefBCR}. 
\Tpga \ may not be so clear, the author wishes someone 
will give more nice definition. 

\bigbreak
{\bf 5.2. Semialgebraic quasi-variety.}\hfil\par\penalty 1000\vskip1em
\par\penalty1000\vskip0.4em plus0.1em minus0.1em
\def\Tpgb{Definition 5.2} %
\proclaim{Definition 5.2}{\rm (Semialgebraic quasi-variety)} 
{\rm A locally ringed space $(A$, $\cR_A)$ is called 
{\it semialgebraic quasi-variety}, if there exists 
a real algebraic quasi-variety $(X$, $\cR_X)$ and 
a finite affine open covering $\{V_i\}_{i=1}^r$ of $X$ which 
satisfies the followings: \par
{\parindent=20pt
\Item{(1)} There exists an injective morphism $\iota \colon (A$, $\cR_A) 
 \lto (X$, $\cR_X)$ as locally ringed spaces, 
and $\iota$ induces a homeomorphism $A \to \iota(A)$ as Euclidean spaces. 
Moreover, $\iota(A)$ is a semialgebraic subset of $X$, 
i.e. $\iota(A) \cap V_i$ is a semialgebraic subset of $V_i$ for 
each $i=1$,$\ldots$, $r$. 
\Item{(2)} $\Zar_X(A) = X$. 
\Item{(3)} Take an arbitral $i \in \{1$, $2$,$\ldots$, $r\}$, 
and take any Euclidean open subset $U \subset \iota^{-1}(V_i)$. 
Put $R_i := \cR_{V_i}(V_i)$. 
For a point $P \in \iota(U)$, let $\frm_P$ be the maximal ideal 
of $R_i$ corresponding to $P$, and let 
\[S_U := \bigcap_{P \in U} \big(R_i - \frm_P\big) \subset R_i.\]
Then $\iota^* : S_U^{-1} R_i \lto \cR_A(U)$ is an isomorphism 
of $\R$-algebra. 

}

Moreover, if $X$ is a real algebraic variety, 
then $A$ is said to be an {\it semialgebraic variety}. 
In this case, the field of fractions $Q(\cR_A(U_i))$ is 
called the {\it field of rational functions}, 
and is denoted by $\Rat(A) := Q\big(\cR_A(U_i)\big)$. 

The Zariski topology and the Euclidean topology on $A$ are 
defined naturally. 
A semialgebraic quasi-variety $A$ is called {\it irreducible} 
if it is irreducible with respect to the Zariski topology. 
$A$ is said to be {\it reduced} 
if $\cR_{A,P}$ has no nilpotent elements except $0$ for 
every $P \in A$. 
$\dim A$ is defined by 
$\displaystyle \dim A = \max_{P \in A} \Krulldim \cR_{A,P}$. 
$A$ is called {\it connected} if it is connected with respect to 
the Euclidean topology. 
Note that $A$ may not be connected even if $A$ is irreducible. 
$A$ is called {\it affine}, if we can choose $X$ to be isomorphic 
to a closed Zariski subset of $\R^n$ for a certain $n$. 

Notions about singularities of $A$ are defined using $\cR_{A,P}$. 
Note that if $Y$ is a $\R$-scheme which represents $X$, 
then $\cR_{A,P} \cong \cO_{Y,P}$. 
We denote 
\begin{align*}
 & \Sing(A) := \big\{ P \in A \; \big| \;
   \hbox{$\cR_{A,P}$ is not a regular local ring} \big\}, \\
 & \Reg(A) := \Int(A) - \Sing(A). 
\end{align*}

A {\it regular map} or {\it holomorphic map} (resp. {\it isomorphism}) between 
semialgebraic quasi-varieties is defined as a morphism 
(resp. isomorphism) of locally ringed space. 


We can choose a real algebraic quasi-variety $X$ 
and a separated scheme $Y$ of finite type over $\R$ 
so that $Y_{\C}$ is complete and $Y$ represents $X$. 
Then, we say $X$ is a {\it real envelope} of $A$, 
and $Y_{\C}$ is a {\it complex envelope} of $A$. }
\endproclaim


$X$ and $Y_{\C}$ are not unique for $A$, 
but it is easy to see that: 

\proclaim{Proposition 5.3} 
{\sl Let $A$ be a semialgebraic quasi-variety, 
$Y_{\C}$ and $Y_{\C}'$ be complex envelopes of $A$. 
Then $Y_{\C}$ and $Y_{\C}'$ are birational equivalent. 
If $A$ is a semialgebraic variety, 
then $\Rat(A) \otimes_{\R} \C = \Rat(Y_{\C})$. }
\endproclaim

This follows from \Tpgk \ given later. 

By this proposition, if $\nu(Y_{\C})$ is a certain birational invariant of 
complex algebraic varieties, 
then we can define $\nu(A) := \nu(Y_{\C})$ to be an invariant of $A$. 
Especially, when $A$ is non-singular semialgebraic variety, 
we can choose $Y$ to be non-singular projective, 
and we can define 
$h^i(A) := \dim_{\C} H^i(Y_{\C}$, $\cO_{Y_{\C}})$ and 
$P_m(A) := \dim_{\C} H^0(Y_{\C}$, $\cO_{Y_{\C}}(mK_{Y_{\C}}))$ for 
$m \geq 0$. 
Using $P_m(A)$, we can define the {\it Kodaira dimension} $\kappa(A)$, 

\proclaim{Remark 5.4} 
{\rm (1) $\Reg(A) \ne \emptyset$ if $A$ is reduced. 

(2) $\Reg(A)$ is not always dense in $A$ with respect 
to the Euclidean topology. 
For example, consider the case that $A$ has an isolated singularity 
as a connected component. 

(3) If $P \in \Reg(A) \cap \Int(A)$ and $\dim A = n$, 
then there exists an Euclidean open neighborhood $P \in U \subset A$ such 
that $U$ is homeomorphic to an open subset of $\R^n$. 

(4) By our definition, an isolated singular locus of $A$ is 
included in $\Int(A)$. 
But $\Sing(A)$ sometimes acts as if it is a boundary. 
So it will be safe to discuss $\Int(A) \cap \Reg(A)$. }
\endproclaim

In complex algebraic geometry, 
a subscheme is a closed subscheme of an open subscheme. 
But to define semialgebraic subvarieties, we must be careful. 
For example, any semialgebraic subset $B$ of a real algebraic variety $A$, 
must be able to be treated as semialgebraic quasi-subvariety of $A$. 

\def\Tpgc{Definition 5.5} %
\proclaim{Definition 5.5}{\rm (Image of a regular map)} 
{\rm Let $A$, $B$ be semialgebraic quasi-varieties, 
and $\varphi \colon A \to B$ be a regular map. 
Let $C := \varphi(A)$. 
By Tarski-Seidenberg theorem, $C$ is a semialgebraic subset of $B$. 
We define $\cR_C$ as the following: 

We may assume $A$ and $B$ are affine, 
since the definition of $\cR_C$ is local. 
Let $R_A := \cR_A(A)$, $R_B := \cR_B(B)$, and 
$\varphi^* \colon R_B \to R_A$ be the homomorphism induced by $\varphi$. 
We put $R := R_B/\Ker \varphi^*$. 
Note that $R$ defines $\Zar_B(C)$. 
For a point $P \in C$, there exists the unique maximal 
ideal $\frm_P \subset R$ corresponding to $P$. 
Put $\displaystyle S := \bigcap_{P \in C} (R - \frm_P)$, 
and $R_C := S^{-1}R$. 
Note that $R_C$ is a $R_B$-module. 
The structure sheaf of $C$ is defined by $\cR_C := \widetilde{R_C}$ 
which is the coherent $\cR_B$-module defined by $R_C$. 

$(C$, $\cR_C)$ is called the {\it image} of $\varphi$, 
and simply denoted by $C = \varphi(A)$. }
\endproclaim

\proclaim{Definition 5.6}{\rm (Semialgebraic quasi-subvariety)} 
{\rm Let $A$, $B$ be semialgebraic quasi-varieties. 
A morphism $\varphi : (B$, $\cR_B) \lto (A$, $\cR_A)$ is 
called an {\it immersion}, if $\varphi$ induces 
an isomorphism $B \lto \varphi(B)$. 

If $B$ is a semialgebraic subset of $A$, 
and the inclusion map $B \to A$ is an immersion, 
then $B$ is called a semialgebraic {\it quasi-subvariety} of $A$. 

If $A$ is a semialgebraic quasi-variety, 
and $B \subset A$ be a semialgebraic subset. 
Then, there exists a unique sheaf of rings $\cR_B$ 
such that $(B$, $\cR_B)$ is a semialgebraic quasi-subvariety 
of $(A$, $\cR_A)$ and $(B$, $\cR_B)$ is reduced. 
$\cR_B$ is called the {\it reduced structure} of $B \subset A$. }
\endproclaim


\proclaim{Definition 5.7}{\rm (Fibre product)} 
{\rm Let $A$, $B$, $C$ be semialgebraic quasi-varieties, 
and $f \colon A \to C$, $g \colon B \to C$ be regular maps. 
The {\it fiber product} $A \times_C B$ is a semialgebraic set 
\[A \times_C B = \big\{ (a,b) \in A \times B \; \big| \; 
 \hbox{$f(a) = g(b)$}\big\}\]
with a structure sheaf $\cR_A \otimes_{\cR_C} \cR_B$. }
\endproclaim

\proclaim{Definition 5.8}{\rm (Inverse image)} 
{\rm Let $A$, $B$ be semialgebraic quasi-varieties, 
and $\varphi \colon A \to B$ be a regular map. 
Let $C \subset B$ be a semialgebraic quasi-subvariety. 
The {\it inverse image} $\varphi^{-1}(C)$ is 
defined as the fiber product $\varphi^{-1}(C) := A \times_B C$. }
\endproclaim

\proclaim{Definition 5.9}{\rm (Birational map)} 
{\rm Let $A$, $B$ be semialgebraic quasi-varieties. 
If there exists Zariski open subsets $U \subset A$ and $W \subset B$ such 
that $\Zar_A(U)=A$, $\Zar_B(W)=B$ and 
there exists a regular map $\varphi \colon U \to W$, 
then we say that there exists a {\it rational map} $\varphi 
 \colon A \cdots\to B$. 
Moreover, if $\varphi \colon U \to W$ is an isomorphism, 
we say that $\varphi \colon A \cdots\to B$ is a {\it birational map}, 
and $A$ and $B$ are {\it birational equivalent}. }
\endproclaim

\def\Tpgk{Proposition 5.10} %
\proclaim{Proposition 5.10} 
{\sl Let $A$, $B$ be semialgebraic quasi-varieties, 
and let $X$, $Y$ be complex envelopes of $A$, $B$.} \par
{\parindent=20pt
\Item{(1)} {\sl If there exists a rational map $\varphi 
 \colon A \cdots\to B$, then there exists 
a rational map $\Phi \colon X_{\C} \cdots\to Y_{\C}$ 
such that $\Phi|_A = \varphi$. }
\Item{(2)} {\sl In (1), if $\varphi$ is a birational map, 
then $\Phi$ is a birational map. }

}
\endproclaim

\Proof 
(1) We may assume $\varphi$ is a regular map. 
Take a point $P \in \Int(A)$ such that $Q := \varphi(P) \in \Int(B)$, 
and take an affine open subset $W \subset Y$ such that $Q \subset W$. 

We can choose $f_1$,$\ldots$, $f_r \in \cR_{Y,Q}$ such that 
we can regard $f_i \in \cO_Y(W)$ 
and $\cO_Y(W) = \C[f_1,\ldots, f_r]$. 
Put $g_j := \varphi^*(f_j) \in \cR_{A,P}$. 
We can find an affine open subset $U \subset X_{\C}$ such 
that $g_1$,$\ldots$, $g_r$ are holomorphic (regular) on $U$, 
and that $U \cap X$ is dense in $X$ and $U \cap A$ is dense in $A$. 
Then, $\psi^* \colon \cR_B \to \cR_A$ induces 
$\Psi^* \colon \cO_Y(W)  \lto \cO_X(U)$. 
$\Psi^*$ induces a rational map $\Phi \colon X \cdots\to Y$. 

(2) is easy. 
\end{proof}

\bigbreak
{\bf 5.3. Some properties of semialgebraic quasi-varieties.}
\par\penalty1000\vskip0.4em plus0.1em minus0.1em
A notion of semialgebraic quasi-varieties brings some merits 
to Real Algebraic Geometry. 

\def\Tpgla{Theorem 5.11} %
\proclaim{Theorem 5.11} 
{\sl Every semialgebraic quasi-variety is affine. 
In other words, if $A$ is a semialgebraic quasi-variety, 
then there exists $m \in \N$ and 
an immersion $\iota \colon A \to \R^m$. }
\endproclaim

\Proof 
Let $(A$, $\cR_A)$ be a semialgebraic quasi-variety. 
We can take a real envelope $(X$, $\cR_X)$ of $A$, 
and a finite affine open covering $\{V_i\}_{i=1}^r$ of $X$ as \Tpgb. 
Fix an integer $1 \leq i \leq r$. 
Then the coordinate ring $\cR_{V_i}(V_i)$ can be represented as 
$\cR_{V_i}(V_i) \cong \R[x_1$,$\ldots$, $x_n]/I_i$ 
where $n$ is an integer which does not depend on $i$, 
and $I_i$ is an ideal of $\R[x_1$,$\ldots$, $x_n]$. 
We can regard $V_i$ as a closed real quasi-subvariety of $U_i = \R^n$, 
and $(x_1$,$\ldots$, $x_n)$ as the coordinate system of $U_i = \R^n$. 
We can consider $x_j$ as a regular function $x_j \colon V_i \to \R$. 
Let 
\[s(P) := \frac{1}{1+x_1(P)^2+\cdots+x_n(P)^2} \hskip 10pt \hbox{and}
 \hskip 10pt t_j(P) := \frac{x_j(P)}{1+x_1(P)^2+\cdots+x_n(P)^2}\]
($j=1$,$\ldots$, $n$) for $P \in U_i$. 

Consider the case $P \in X - U_i$. 
There exists a sequence of poins $\{P_k\}_{k=1}^{\infty}$ in 
$U_i$ such that $\displaystyle \lim_{k \to \infty} P_k = P$. 
Then, there exsists $1 \leq \exists l \leq n$ such 
that $\displaystyle \lim_{k \to \infty} x_l(P_i) = \infty 
 \in \P_{\R}^1 = \R \cup \{\infty\}$. 
Therefore $\displaystyle \lim_{k \to \infty} s(P_k) = 0$ 
and $\displaystyle \lim_{k \to \infty} t_j(P_k) = 0$. 
Thus we can define $s(P) = 0$ and $t_j(P) = 0$. 

In this way, we can extend $s \colon V_i \to \R$ and $t_j \colon V_i \to \R$ 
to algebraic functions $s \colon X \to \R$ and $t_j \colon X \to \R$. 
These are regular holomorphic functions on $X$. 

Let $(y_0$, $y_1$,$\ldots$, $y_n)$ be a coordinate system 
of $\R^{n+1}$. 
We define a map $\Phi_i \colon X \lto \R^{n+1}$ 
by $y_0 = s(P)$ and $y_j = t_j(P)$ ($j=1$,$\ldots$, $n$) for $P \in X$. 
This $\Phi_i$ is a regular map as semialgebraic quasi-varieties.

Note that $\Phi_i(X)$ is a semialgebraic quasi-variety 
by Tarski-Seidenberg theorem, 
but is not always an algebraic quasi-variety. 

Let $P \in V_i$ and $Q := \Phi_i(P)$. Note thet $y_0(Q) = s(P) \ne 0$. 
Then $x_j(P) = t_j(P)/s(P) = y_j(Q)/y_0(Q)$. 
Thus $\Phi_i|_{V_j} \colon V_j \to \Phi_i(V_j)$ is an isomorphism. 

We define a regular map $\Phi \colon X \to \R^{r(n+1)}$ 
by $\Phi(P) = \big(\Phi_1(P)$,$\ldots$, $\Phi_r(P)\big)$. 
Let $B := \Phi(A)$ and $Y = \Phi(X)$. 

Since $\Phi_i \colon V_i \lto \R^{n+1}$ is an immersion, 
we have $\Phi|_{V_i} \colon V_i \lto \R^{r(n+1)}$ is also an immersion. 
Thus $\Phi \colon X \lto \R^{2rn}$ is an immersion. 
Note that we can take $\cR_B$ and $\cR_Y$ 
so that $\cR_A = \Phi^* \cR_B$ and $\cR_X := \Phi^* \cR_Y$. 
So, $(X$, $\cR_X) \cong (Y$, $\cR_Y)$ and 
$(A$, $\cR_A) \cong (B$, $\cR_B)$. 
Thus $\Phi \colon A \to \R^m$ is an immersion. 
\end{proof}

\proclaim{Remark 5.12} 
{\rm A real algebraic variety is an affine semialgebraic variety, 
but is not always a real affine variety. 
For example, $\R^2 - \{(0,0)\}$ is not a real affine variety. }
\endproclaim

\def\Tpglc{Corollary 5.13} %
\proclaim{Corollary 5.13} 
{\sl Let $A$ be a semialgebraic quasi-variety 
(or a real algebraic quasi-variety) and 
$\cF$ be a quasi-coherent $\cR_A$-module. 
Put $R_A := \cR_A(A)$ and $F := \cF(F)$. 
Then, $\cR_A$ and $\cF$ are the sheaves obtained 
as $\widetilde{R_A}$ and $\widetilde{F}$.} 
\endproclaim

Note that $R_A$ is a Noetherian ring, 
but is not finitely generated over $\R$ if $\dim A \geq 1$. 
Each maximal ideal of $R_A$ corresponds to a certain point of $A$. 

\def\Tpgm{Corollary 5.14} %
\proclaim{Corollary 5.14} 
{\sl Let $A$ be a semialgebraic quasi-variety 
(or a real algebraic quasi-variety) and 
$\cF$ be a quasi-coherent $\cR_A$-module. 
Then, $H^i(A$, $\cF) = 0$ for all $i> 0$. }
\endproclaim

\Proof 
There exists an immersion $\iota \colon A \to \R^m$ by \Tpgla. 
Thus we may assume $A \subset \R^m$. 
Let $X$ be a real envelope of $A$ such that $A \subset X \subset \R^m$. 
Let $R_X := \cR_X(X)$ and $R_A := \cR_A(A)$. 
We can present as $R_A = S_A^{-1}R_X$ by a certain 
multiplicatively closed set $S_A \subset R_X$. 
Since $R_A$ is an $R_X$-module, $\cF$ is a 
quasi-coherent $\cR_X$-module. 

Put $M := \cF(A) = H^0(X$, $\cF)$. Then $\cF = \tilde{M}$. 
Let $0 \to M \to I^0 \to I^1 \to \cdots$ be an injective resolution 
of $M$ as a $R_X$-module. 
Then $0 \to \cF \to \tilde{I}^0 \to \tilde{I}^1 \to \cdots$ is 
an injective resolution of $\cF$ as an $\cR_X$-module. 

We can prove that $\tilde{I}^i$ are flasque on $X$ 
as \cite[Chap.II, Proposition 3.4]{RefH}. 

Thus we have 
\[H^i(A, \, \cF) \cong H^i\big(X, \, \cF\big) = 0\]
for $i>0$ (see \cite[Chap.III, Theorem 3.5]{RefH}). 
\end{proof}

Let $A := \R^2 - \{(0,0)\}$ and $A_{\C} := \C^2 - \{(0,0)\}$. 
As complex algebraic variety, $H^1(A_{\C}$, $\cO_{A_{\C}}) \ne 0$ 
(see \cite[Chap.III, Exercise 4.3]{RefH}). 
But $H^1(A$, $\cR_A) = 0$. 
Sheaves $\cR_A$ and $\cO_{A_{\C}}$ have very different properties. 

\smallskip

Similarly, birational geometries of complex and real algebraic varieties 
are very different. 
In a complete complex algebraic variety, 
exceptional subsets are special subsets. 
This is not true for complete real algebraic varieties. 

\def\Tpgn{Theorem 5.15} %
\proclaim{Theorem 5.15} 
{\sl Let $A$ be a semialgebraic variety, 
$E \subset A$ be a closed semialgebraic subset 
such that $E = \Zar_A(E) \subsetne A$. 
Then there exists a semialgebraic variety $B$ and 
a regular surjective morphism $\varphi \colon A \to B$ 
such that $P := \varphi(E)$ is a point 
and that $\varphi|_{A-E} : (A - E) \lto (B - P)$ is 
an isomorphism, i.e. $\varphi$ is a contraction of $E$ to a point $P$. }
\endproclaim

\Proof 
We may assume $A \subset \R^n$. 
Let $f_1$,$\ldots$, $f_r$ be defining polynomials of $\Zar_{\R^n}(E)$ in 
$\R[x_1$,$\ldots$, $x_n]$. 
Consider a map $\Phi \colon \R^n \to \R^{rn}$ defined by 
linear system with the basis $\big\{ x_i f_j$ $\big|$ $1 \leq i \leq n$, 
$1 \leq j \leq r \big\}$. 
$\Phi$ is a regular map. 
Put $B := \Phi(A)$ and $\varphi := \Phi|_A \colon A \to B$. 
Then, $B$ and $\varphi$ satisfy the conclusion of the Theorem. 
\end{proof}


\bigbreak


\begin{thebibliography}{20}

\def\refJ#1#2#3#4#5#6#7{%
{\rm #1}: {\rm #2}. {\rm #3} {\rm #4}, {\rm #6} (#5)}%
\def\refJb#1#2#3#4#5#6#7{%
{\rm #1}: {\rm #2}. {\rm #3} {\rm #4}, {\rm #6} (#5).}%

\def\refA#1#2#3{%
{\rm #1}: {\rm #2}. {\rm arXiv}:{\rm #3}.}%
\def\refAb#1#2#3{%
{\rm #1}: {\rm #2}. \hfill\break{\rm arXiv}:{\rm #3}.}%

\def\refAY#1#2#3#4{%
{\rm #1}: {\rm #2}. (#4) {\rm arXiv}:{\rm #3}.}%
\def\refAb#1#2#3{%
{\rm #1}: {\rm #2}. \hfill\break{\rm arXiv}:{\rm #3}.}%

\def\refB#1#2#3#4{
{\rm #1}: {\rm #2}, #3 (#4)}%

\bibitem{RefAa}
\refJ{T. Ando}
{Some homogeneous cyclic inequalities of three variables 
of degree three and four}
{Aust. J. Math. Anal. Appl.}{7 (2) Art. 11}{2010}{}{} 

\bibitem{RefAb}
\refJ{T. Ando}
{Cubic and quartic cyclic homogeneous inequalities of Three Variables}
{Math. Inequal. Appl.}{16}{2013}{127-142}{}

\bibitem{RefAc}
\refJ{T. Ando}
{Discriminants of Cyclic Homogeneous Inequalities of Three Variables}
{J. Alg.}{\bf 514}{2018}{384-441}{}

\bibitem{RefAd}
\refJ{T. Ando}
{ Extremal Cubic Inequalities of Three 
Variables}{Turkish J. Ineq.}{7}{2024}{1-34}{}

\bibitem{RefBHORS}
\refJ{G. Blekherman \& J. Hauenstein \& J. C. Ottem \& 
K. Ranestad \& B. Sturmfels}
{Algebraic Boundaries of Hilbert's SOS Cones}
{Compositio Math}{148}{2012}{1717-1735}{}

\bibitem{RefB}
\refB{G. Blekherman \& P. A. Parrilo \& R. Thomas, ed.}
{Semidefinite Optimization and Convex Algebraic Geometry}
{SIAM}{2013}

\bibitem{RefGKBR} 
\refJ{G. Blekherman \& C. Riener}
{Symmetric Nonnegative Forms and Sums of Squares}
{Disc. Comp. Geom.}{65}{2021}{764-799}{}

\bibitem{RefBCR}
\refB{J. Bochnak \& M. Coste \& M.F. Roy}
{Real Algebraic Geometry}{Springer}{1998} 

\bibitem{RefCL}
\refJ{M. D. Choi \& T. Y. Lam}
{Extremal Positive Semidefinite Forms}
{Math. Ann.}{231}{1977}{1-18}
{https://doi.org/10.1007/BF01360024}

\bibitem{RefCLRb}
\refJ{M.D. Choi \& T.Y. Lam \& B. Reznick}
{Real zeros of positive semidefinite forms I}
{Math. Z.}{171}{1980}{1-26}{}

\bibitem{RefCLR}
\refJ{M. D. Choi \& T. Y. Lam  \& B. Reznick}
{Even Symmetric Sextics}{Math. Z.}{195}{1987}{559-580}
{https://doi.org/10.1007/BF01166704}

\bibitem{RefCLRc}
\refJ{M. D. Choi \& T. Y. Lam  \& B. Reznick}
{Positive Sextics and Schur's Inequalities}{J. Alg.}{141}{1991}{36-77} 
{https://doi.org/10.1016/0021-8693(91)90203-K}

\bibitem{RefCa}
\refJ{V. C\^{\i}rtoaje}
{On the cyclic homogeneous polynomial inequalities of degree four}
{J. Inequal. Pure and Appl. Math.}{10 (3) Art. 67}{2009}{}{}

\bibitem{RefGKR} 
\refJ{C. Goel \& S. Kuhlmann \& B. Reznick}
{The analogue of Hilbert's 1888 theorem for even symmetric forms}
{J. Pure Appl. Alg.}{221}{2017}{1438-1448}
{https://doi.org/10.1016/j.jpaa.2016.10.003}

\bibitem{RefHs} 
\refJ{W. R. Harris}
{Real Even Symmetric Ternary Forms}{J. Alg.}{222}{1999}{204-245}{}

\bibitem{RefH}
\refB{R. Hartshorne}
{Algebraic Geometry}{Springer}{1977}

\bibitem{RefHil}
\refJ{D. Hilbert}
{Uber die Darstellung definiter Formen als Summe von Formenquadraten}
{Math. Ann.}{32}{1888}{342-350}
{https://doi.org/10.1007/BF01443605}

\bibitem{RefKa}
\refAY{J. Koll\'ar}
{Real Algebraic Surfaces}{alg-geom/9712003}{1997}

\bibitem{RefMM}
\refA{M. Milev \& N. Milev}
{The `core' of symmetric homogeneous polynomial inequalities 
of degree four of three real variables}{arXiv:1604.00993}

\bibitem{RefN}
\refJ{J. Nie}
{Discriminants and Nonnegative Polynomials}
{J. Symbolic Comput}{47}{2012}{167-191}{}

\bibitem{RefR}
\refJ{B. Reznick}
{Forms Derived from the Arithmetic-Geometric Inequality}
{Math. Ann.}{238}{1989}{431-464}{}

\bibitem{RefRie}
\refJ{C. Riener}
{ On the degree and half-degree principle for symmetric polynomials}
{J. Pure Appl. Alg.}{216.4}{2012}{850-856}{}

\bibitem{RefTa}
\refJ{V. Timofte}
{On the positivity of symmetric polynomial functions.
Part I: General results}{J. Math. Anal. Appl.}{284}{2003}{174-190}{}

\bibitem{RefTb}
\refJ{V. Timofte}
{On the positivity of symmetric polynomial functions. Part II: 
Lattice general results and positivity criteria for degrees 4 and 5}
{J. Math. Anal. Appl.}{304}{2005}{652-667}{}

\bibitem{RefTc}
\refJ{V. Timofte}
{On the positivity of symmetric polynomial functions. 
Part III: Extremal polynomials of degrees 4}
{J. Math. Anal. Appl.}{307}{2005}{565-578}{}

\end{thebibliography}
\end{document}